\documentclass[12pt,oneside]{book}
  \usepackage{amsfonts,amssymb,mathrsfs,pb-diagram,makeidx}
  \usepackage[all]{xy}
  \makeindex

\def \1{{\bf 1}}
\def \a{{{\mathfrak a}}}

\def \al{\alpha}
\def \A{{\mathbb A}}
\def \Ad{{\rm Ad}}

\def \b{{{\mathfrak b}}}
\def \be{\beta}
\def \bs{\backslash}

\def \bwedge{\left.\bigwedge\!\right.}

\def \card{{\rm card}}

\def \C{{\mathbb C}}

\def \CC{{\mathscr C}}
\def \CCC{{\cal C}}
\def \CCD{{\cal D}}
\def \CD{{\mathscr D}}
\def \CE{{\mathscr E}}
\def \CEE{{\mathscr E}}
\def \CF{{\mathscr F}}
\def \CG{{\mathscr G}}

\def \CCO{{\mathscr O}}
\def \CO{{\cal O}}
\def \cos{{\rm cos \hspace{2pt}}}
\def \coth{{\rm coth}}
\def \CP{{\mathscr P}}

\def \d{\delta}
\def \D{\Delta}
\def \det{{\rm det \hspace{2pt}}}
\def \diag{{\rm diag}}

\def \eac{{\acute{\rm e}}}
\def \eacit{{\acute{\it e}}}
\def \ell{{\rm ell}}
\def \ep{\varepsilon}
\def \exp{{\rm exp}}

\def \fin{{\rm fin}}
\def \F{{\mathbb F}}
\def \CF{{\mathfrak F}}

\def \g{{{\mathfrak g}}}
\def \ga{\gamma}
\def \Ga{\Gamma}
\def \gen{{\rm gen}}
\def \GL{{\rm GL}}
\def \Gr{{\rm Gr \hspace{2pt}}}
\def \h{{{\mathfrak h}}}
\def \H{{\mathbb H}}

\def \hra{\hookrightarrow}

\def \Im{{\rm Im \hspace{2pt}}}
\def \ind{{\rm ind \hspace{2pt}}}
\def \Ind{{\rm Ind \hspace{2pt}}}
\def \Inv{{\rm Inv \hspace{2pt}}}
\def \k{{{\mathfrak k}}}

\def \ka{\kappa}

\def \la{\lambda}

\def \La{\Lambda}
\def \li{{\rm li}}
\def \Lie{{\rm Lie}}
\def \Lin{{\rm Lin}}
\def \m{{{\mathfrak m}}}
\def \Mat{{\rm Mat}}
\def \mod{{\rm mod}}
\def \n{{{\mathfrak n}}}

\def \N{\mathbb N}
\def \nreg{\rm nreg}

\def \O{{\rm O}}
\def \om{\omega}
\def \Om{\Omega}
\def \ox{\otimes}
\def \p{{{\mathfrak p}}}

\def \prf{\noindent{\bf Proof: }}

\def \PGL{{\rm PGL}}

\def \qed{\ifmmode\eqno $\square$ \else\noproof\vskip 12pt plus 3pt minus 9pt \fi}
 \def\noproof{{\unskip\nobreak\hfill\penalty50\hskip2em\hbox{}%
     \nobreak\hfill $\square$ \parfillskip=0pt%
     \finalhyphendemerits=0\par}}
\def \Q{\mathbb Q}

\def \R{{\mathbb R}}
\def \Re{{\rm Re \hspace{1pt}}}

\def \ra{\rightarrow}
\def \Ra{\Rightarrow}

\def \reg{{\rm reg}}

\def \rw{{\rm rw}}
\def \rnw{{\rm rnw}}
\def \si{\sigma}
\def \Si{\Sigma}
\def \sin{{\rm sin \hspace{2pt}}}
\def \SL{{\rm SL}}
\def \sl{{\rm sl}}
\def \SO{{\rm SO}}

\def \SS{{\rm S}}

\def \supp{{\rm supp}}

\def \tanh{{\rm tanh}}
\def \Th{\Theta}
\def \th{\theta}
\def \tr{{\hspace{1pt}\rm tr\hspace{2pt}}}
\def \Upp{{\rm Upp}}
\def \ut{\tilde{u}}

\def \vol{{\rm vol}}

\def \x{\times}
\def \z{{\mathfrak z}}
\def \ze{\zeta}
\def \Z{\mathbb Z}
\def \={\ =\ }

\newcommand{\rez}[1]{\frac{1}{#1}}

\renewcommand{\matrix}[4]{\left( \begin{array}{cc}#1 & #2 \\ #3 & #4 \end{array}
            \right)}
\newcommand{\matrixfour}[4]
	   {\left( \begin{array}{cccc} 
	       #1 & \  & \  & \  \\ 
	       \  & #2 & \  & \  \\  
	       \  & \  & #3 & \  \\  
	       \  & \  & \  & #4 \end{array}\right)}
\newcommand{\matrixtwo}[2]{\matrix {#1}{}{}{#2}}
\newtheorem{theorem}{Theorem}[section]
\newtheorem{conjecture}[theorem]{Conjecture}
\newtheorem{lemma}[theorem]{Lemma}

\newtheorem{proposition}[theorem]{Proposition}

\begin{document}


\thispagestyle{empty}
\begin{center}

\ 
\bigskip
\bigskip
\bigskip
\bigskip
\bigskip
\bigskip
\bigskip
\bigskip
\bigskip
\bigskip

\Huge \textbf{Class Numbers of Orders}

\bigskip

\textbf{in}

\bigskip

\textbf{Quartic Fields}

\bigskip
\bigskip
\bigskip
\bigskip
\bigskip
\bigskip

\Large Dissertation

\bigskip
\bigskip

\normalsize der Fakult\"at f\"ur Mathematik und Physik

\bigskip

der Eberhard-Karls-Universit\"at T\"ubingen

\bigskip

zur Erlangung des Grades eines Doktors der Naturwissenschaften

\bigskip
\bigskip
\bigskip
\bigskip

vorgelegt von

\bigskip

\textsc{Mark Pavey}

\bigskip

aus Cheltenham Spa, UK

\bigskip
\bigskip
\bigskip
\bigskip
\bigskip
\bigskip
\bigskip
\bigskip
\bigskip
\bigskip

2006

\end{center}

\newpage
\thispagestyle{empty}

\vspace*{\stretch{1}}

Mundliche Pr\"ufung: 11.05.2006

\bigskip

\begin{tabular}{@{}ll}
Dekan: & Prof.Dr.P.Schmid \\
1.Berichterstatter: & Prof.Dr.A.Deitmar \\
2.Berichterstatter: & Prof.Dr.J.Hausen
\end{tabular}

\newpage
\thispagestyle{empty}

\vspace*{\stretch{1}}

\begin{center}

Dedicated to my parents:

\bigskip

\textsc{Deryk and Frances Pavey},

\bigskip

who have always supported me.

\end{center}

\vspace*{\stretch{2}}

\newpage
\thispagestyle{empty}
\ 
\newpage

\chapter*{Acknowledgments}

First of all thanks must go to my doctoral supervisor Professor Anton Deitmar, without whose enthusiasm and unfailing generosity with his time and expertise this project could not have been completed.

The first two years of research for this project were undertaken at Exeter University, UK, with funding from the Engineering and Physical Sciences Research Council of Great Britain (EPSRC).  I would like to thank the members of the Exeter Maths Department and EPSRC for their support.  Likewise I thank the members of the Mathematisches Institut der Universit\"at T\"ubingen, where the work was completed.

Finally, I wish to thank all the family and friends whose friendship and encouragement have been invaluable to me over the last four years.  I have to mention particularly the Palmers, the Haywards and my brother Phill in Seaton, and those who studied and drank alongside me: Jon, Pete, Dave, Ralph, Thomas and the rest.  Special thanks to Paul Smith for the Online Chess Club, without which the whole thing might have got done more quickly, but it wouldn't have been as much fun.

\newpage
\thispagestyle{empty}
\ 
\newpage

\frontmatter

\tableofcontents

  \chapter{Introduction}
    \label{ch:Intro}
    \markright{\textnormal{Introduction}}
    In this thesis we present two main results.  The first is a Prime Geodesic Theorem for compact symmetric spaces formed as a quotient of the Lie group $\SL_4(\R)$.  The second is an application of the Prime Geodesic Theorem to prove an asymptotic formula for class numbers of orders in totally complex quartic fields with no real quadratic subfield.  Before stating our results we give some background.

Let $\CD$ be the set of all natural numbers $D\equiv 0,1\ \mod\,4$ with $D$ not a square.  Then $\CD$ is the set of all discriminants of orders in real quadratic fields.  For $D\in\CD$ the set
$$
\CO_D=\left\{\frac{x+y\sqrt{D}}{2}:x\equiv yD\ \mod\,2\right\}
$$
is an order in the real quadratic field $\Q(\sqrt{D})$ with discriminant $D$.  As $D$ varies, $\CO_D$ runs through the set of all orders of real quadratic fields.  For $D\in\CD$ let $h(\CO_D)$ denote the class number and $R(\CO_D)$ the regulator of the order $\CO_D$.  It was conjectured by Gauss (\cite{Gauss86}) and proved by Siegel (\cite{Siegel44}) that, as $x$ tends to infinity
$$
\sum_{{D\in\CD}\atop{D\leq x}}h(\CO_D)R(\CO_D)=\frac{\pi^2 x^{3/2}}{18\ze(3)}+O(x\log x),
$$
where $\ze$ is the Riemann zeta function.

For a long time it was believed to be impossible to separate the class number and the regulator in the summation.  However, in \cite{Sarnak82}, Theorem 3.1, Sarnak showed, using the Selberg trace formula, that as $x\ra\infty$ we have
$$
\sum_{{D\in\CD}\atop{R(\CO_D)\leq x}}h(\CO_D)\sim\frac{e^{2x}}{2x}.
$$
More sharply,
$$
\sum_{{D\in\CD}\atop{R(\CO_D)\leq x}}h(\CO_D)=L(2x)+O\left(e^{3x/2}x^2\right)
$$
as $x\ra\infty$, where $L(x)$ is the function
$$
L(x)=\int_1^x\frac{e^t}{t}\,dt.
$$
Sarnak established this result by identifying the regulators with lengths of closed geodesics of the modular curve $\SL_2(\Z)\bs\H$, where $\H$ denotes the upper half-plane, and using the Prime Geodesic Theorem for this Riemannian surface.  Actually Sarnak proved not this result but the analogue where $h(\CO)$ is replaced by the class number in the narrower sense and $R(\CO)$ by a ``regulator in the narrower sense".  But in Sarnak's proof the group $\SL_2(\Z)$ can be replaced by $\PGL_2(\Z)$ giving the above result.

The Prime Geodesic Theorem in this context gives an asymptotic formula for the number of closed geodesics on the surface $\SL_2(\Z)\bs\H$ with length less than or equal to $x>0$.  This formula is analogous to the asymptotic formula for the number of primes less than $x$ given in the Prime Number Theorem.  The Selberg zeta function (see \cite{Selberg56}) is used in the proof of the Prime Geodesic Theorem in a way analogous to the way the Riemann zeta function is used in the proof of the Prime Number Thoerem (see \cite{Chandrasekharan68}).  The required properties of the Selberg zeta function are deduced from the Selberg trace formula (\cite{Selberg56}).

It seems that following Sarnak's result no asymptotic results for class numbers in fields of degree greater than two were proven until in \cite{Deitmar02}, Theorem 1.1, Deitmar proved an asymptotic formula for class numbers of orders in complex cubic fields, that is, cubic fields with one real embedding and one pair of complex conjugate embeddings.  Deitmar's result can be stated as follows.

Let $S$ be a finite set of prime numbers containing at least two elements and let $C(S)$ be the set of all complex cubic fields $F$ such that all primes $p\in S$ are non-decomposed in $F$.  For $F\in C(S)$ let $O_F(S)$ be the set of all isomorphism classes of orders in $F$ which are maximal at all $p\in S$, ie. are such that the completion $\CO_p=\CO\ox\Z_p$ is the maximal order of the field $F_p=F\ox\Q_p$ for all $p\in S$.  Let $O(S)$ be the union of all $O_F(S)$, where $F$ ranges over $C(S)$.  For a field $F\in C(S)$ and an order $\CO\in O_F(S)$ define
$$
\la_S(\CO)=\la_S(F)=\prod_{p\in S}f_p(F),
$$
where $f_p(F)$ is the inertia degree of $p$ in $F$.  Let $R(\CO)$ denote the regulator and $h(\CO)$ the class number of the order $\CO$.

For $x>0$ we define
$$
\pi_S(x)=\sum_{{\CO\in O(S)}\atop{R(\CO)\leq x}}\la_S(\CO)h(\CO).
$$
Then as $x\ra\infty$ we have
$$
\pi_S(x)\sim\frac{e^{3x}}{3x}.
$$
More sharply,
$$
\pi_S(x)=L(3x)+O\left(\frac{e^{9x/4}}{x}\right)
$$
as $x\ra\infty$.

This result was again proved by means of a Prime Geodesic Theorem, this time for symmetric spaces formed as a compact quotient of the group $\SL_3(\R)$.  There exists a Selberg trace formula for such spaces (\cite{Wallach76}) by means of which the required properties of a generalised Selberg zeta function can be deduced (\cite{Deitmar00}) in order to prove the Prime Geodesic Theorem.  The class number formula is then deduced by means of a correspondence between primitive closed geodesics and orders in complex cubic number fields, under which the lengths of the geodesics correspond to the regulators of the number fields.

In this thesis we follow the methods of \cite{Deitmar02}.  In the first three chapters we prove a Prime Geodesic Theorem for compact quotients of $\SL_4(\R)$.  In the final two chapters we give an arithmetic interpretation in terms of class numbers.  Our main results can be stated as follows.

Let $G=\SL_4(\R)$ and let $K$ be the maximal compact subgroup $\SO(4)$.  Let $\Ga\subset G$ be discrete and cocompact.  We then have a one to one correspondence between conjugacy classes in $\Ga$ and free homotopy classes of closed geodesics on the symmetric space $X_{\Ga}=\Ga\bs G/K$.  Let
$$
A^-=\left\{ \matrixfour{a}{a}{a^{-1}}{a^{-1}}:0<a<1\right\}
$$
and
$$
B=\matrixtwo {\SO(2)}{\SO(2)}.
$$
Let $\CE(\Ga)$ be the set of primitive conjugacy classes $[\ga]$ in $\Ga$ such that $\ga$ is conjugate in $G$ to an element $a_{\ga}b_{\ga}$ of $A^-B$.  For $\ga\in\Ga$ we write $a_{\ga}$ also for the top left entry in the matrix $a_{\ga}$ and define the length $l_{\ga}$ of $\ga$ to be $8\log a_{\ga}$.  Let $G_{\ga},\Ga_{\ga}$ be the centralisers of $\ga$ in $G$ and $\Ga$ respectively and let $K_{\ga}=K\cap G_{\ga}$.  For $x>0$ define the function
$$
\pi(x)=\sum_{{[\ga]\in\CE(\Ga)}\atop{e^{l_{\ga}}\leq x}}\chi_1(\Ga_{\ga}),
$$
where $\chi_1(\Ga_{\ga})$ is the first higher Euler characteristic of the symmetric space $X_{\Ga_{\ga}}=\Ga_{\ga}\bs G_{\ga}/K_{\ga}$.

\begin{theorem}\textnormal{(Prime Geodesic Theorem)}
For $x\rightarrow\infty$ we have
$$
\pi(x)\sim\frac{2x}{\log x}.
$$
More sharply,
$$
\pi(x)=2\,\li(x)+O\left(\frac{x^{3/4}}{\log x}\right)
$$
as $x\rightarrow\infty$, where $\li(x)=\int_2^x \rez{\log t} dt$ is the integral logarithm.
\end{theorem}

We now state our result on class numbers.  Let $S$ be a finite, non-empty set of prime numbers containing an even number of elements.  We define the sets $C(S)$ and $O(S)$ and the constants $\la_S(\CO)$ as above, replacing complex cubic fields with totally complex quartic fields in the definitions.  A totally complex quartic field has at most one real quadratic subfield, as can be seen by comparing numbers of fundamental units.  Let $C^c(S)\subset C(S)$ be the subset of fields with no real quadratic subfield and let $O^c(S)\subset O(S)$ be the subset of isomorphy classes of orders in fields in $C^c(S)$.  Let $R(\CO)$ and $h(\CO)$ once again denote respectively the regulator and class number of the order $\CO$.

\begin{theorem}\textnormal{(Main Theorem)}
\label{thm:main1}
For $x>0$ we define
$$
\pi_S(x)=\sum_{{\CO\in O^c(S)}\atop{R(\CO)\leq x}}\la_S(\CO)h(\CO).
$$
Then as $x\ra\infty$ we have
$$
\pi_S(x)\sim\frac{e^{4x}}{8x}.
$$
\end{theorem}

Transfering the methods of \cite{Deitmar02} to the case of totally complex quartic fields does not proceed entirely smoothly, there are various technical difficulties to be overcome.  These have been successfully overcome to the point of proving the Prime Geodesic Theorem.  However, the correspondence between primitive closed geodesics and orders is considerably more complicated than in either the real quadratic or the complex cubic case.  As a result of this extra complexity our result is weaker than that obtained in the complex cubic case.  We have had to impose the extra condition on the finite set $S$ of prime numbers that it must contain an even number of elements and we have been unable to provide an error term in our final asymptotic.

Also we have had to restrict ourselves to counting the orders in fields without a real quadratic subfield, as the fields with a real quadratic subfield cannot be counted using our method.  Indeed, any real quadratic field which occurs as a subfield of a totally complex quartic field may in fact occur as a subfield of infinitely many such fields.  The fundamental units in the quadratic field are powers of fundamental units in the quartic fields, so the (possibly infinitely many) quartic fields all have regulator less than or equal to that of the quadratic field.  Hence an asymptotic formula for a sum of class numbers of orders bounded by the regulators is not even possible in the case of totally complex quartic fields with a real subfield.  

In comparison with the real quadratic and complex cubic cases we might have expected to get the asymptotic
$$
\pi_S(x)\sim\frac{e^{4x}}{4x}.
$$
In the result that we do in fact get there is an extra factor of one half.  The fact that we have had to restrict the fields over which we are counting gives an explanation for this discrepancy.

We also prove the following asymptotic in which we introduce an extra factor into the summands.  For an order $\CO\in O(S)$ this extra factor is defined in terms of the arguments of the fundamental units of $\CO$ under the embeddings of $F$ into $\C$ as follows.  If $\ep$ is a fundamental unit in $\CO$ then $\ep^{-1},\ze\ep$ and $\ze\ep^{-1}$ are also fundamental units in $\CO$, where $\ze$ is a root of unity contained in $\CO$.  Let
$$
\nu(\CO)=\rez{2\mu_{\CO}}\sum_{\ep}\prod_{\al}\left(1-\frac{\al(\ep)}{|\al(\ep)|}\right),\index{$\nu(\CO)$}
$$
where $\mu_{\CO}$ is the number of roots of unity in $\CO$, the sum is over the $2\mu_{\CO}$ different fundamental units of $\CO$ and the product is over the embeddings of $\CO$ into $\C$.  We have:

\begin{theorem}
\label{thm:main2}
For $x>0$ we define
$$
\tilde{\pi}_S(x)=\sum_{{\CO\in O^c(S)}\atop{R(\CO)\leq x}}\nu(\CO)\la_S(\CO)h(\CO).
$$
Then as $x\ra\infty$ we have
$$
\tilde{\pi}_S(x)\sim\frac{e^{4x}}{2x}.
$$
\end{theorem}

An element $\ga$ in $\Ga$ is called \emph{regular} if it centraliser in $G$ is a torus and \emph{non-regular} otherwise.  The factor $\nu(\CO)$ was originally introduced in order to separate the contribution of non-regular elements from that of the regular elements in the Prime Geodesic Theorem.  As it turned out, the complexity of the correspondence between geodesics and orders meant that we did not really gain anything from this approach.  However, the factor $\nu(\CO)$ contains information about the arguments of the fundamental units in the order $\CO$ which is interesting in its own right.  In comparison with Theorem \ref{thm:main1} we can see that ``on average" the value of $\nu(\CO)$ as $R(\CO)$ goes to infinity is 4.  If $\ep$ is a fundamental unit in $\CO$ with arguments $\th,-\th,\phi,-\phi$ under the four embeddings of $\CO$ into $\C$ then a simple calculation shows that
$$
\prod_{\al}\left(1-\frac{\al(\ep)}{|\al(\ep)|}\right)=4(1-\cos 2\th)(1-\cos 2\phi).
$$
Since this takes ``on average" the value 4 we can see that there is a sense in which we can say that the arguments of the fundamental units in the orders $\CO$ are evenly distributed about $\pm\pi/2$ as $R(\CO)$ goes to infinity.

In what follows we give a chapter by chapter summary of the arguments and techniques used in the proof of our main results.  In particular we shall point out the difficulties that arise in applying the methods of \cite{Deitmar02} to our situation.  Chapter \ref{ch:EulerChar} provides some necessary background and preliminary results.  In Chapter \ref{ch:Ruelle} we introduce the zeta functions we shall be studying and prove some of their analytic properties.  The results of Chapter \ref{ch:EulerChar} are made use of here in the definitions of the zeta functions and proofs of their analytic properties.  In Chapter \ref{ch:PGT} we apply standard techniques from analytic number theory to the results of Chapter \ref{ch:Ruelle} in order to prove the Prime Geodesic Theorem required in our context.  In Chapter \ref{ch:DivAlg} we show that the fields we are interested in can be obtained as maximal subfields of a suitably chosen division algebra over $\Q$ and count the number of embeddings of the fields into the division algebra.  Finally, in Chapter \ref{ch:Compare} we apply the results of Chapter \ref{ch:DivAlg} to give the arithmetic interpretation of the Prime Geodesic Theorem in terms of class numbers.

The correspondence between closed primitive geodesics in the Prime Geo-desic Theorem and orders in totally complex quartic fields goes via the choice of a division algebra $M(\Q)$ of degree four over $\Q$ whose maximal subfields include the ones we are interested in, ie. the fields in the set $C^c(S)$.  We are in fact able to choose $M(\Q)$ so that as well as having the properties described above it also satisfies $M(\Q)\ox\R\cong\Mat_4(\R)$.  The map $M(\Q)\ox\R\ra\R$ induced from the reduced norm on $M(\Q)$ agrees with the determinant map on $\Mat_4(\R)$.  We define
$$
\CG(\R)=\{X\in M(\Q)\ox\R:\det X=1\}\cong\SL_4(\R).
$$
If we let $M(\Z)$ denote the maximal order of $M(\Q)$ (\cite{Reiner75}) then
$$
\Ga=\left(M(\Z)\ox 1\right)\cap\CG(\R)
$$
is a discrete, cocompact subgroup of $\SL_4(\R)$ (\cite{BorelHarder78}).  Note that $\Ga$ is not necessarily torsion free.

The group $\Ga$ forms the link between geodesics and class numbers in the following sense.  Firstly, there is a one-to-one correspondence between conjugacy classes in $\Ga$ and free homotopy classes of closed geodesics on the symmetric space $X_{\Ga}=\Ga\bs\SL_4(\R)/\SO(4)$.  Under this correspondence primitive elements of $\Ga$ correspond to primitive geodesics.  Secondly, the primitive elements of $\Ga$ viewed as a subset of the maximal order $M(\Z)$ in $M(\Q)$ correspond to fundamental units of orders of subfields of $M(\Q)$.  We shall be interested in the primitive conjugacy classes in $\Ga$ which correspond to orders in totally complex quartic fields with no real quadratic subfield.

Our strategy is to use a suitably defined generalised Selberg zeta function to get information about the distribution of the primitive, closed geodesics on $X_{\Ga}$ which correspond to orders in $O^c(S)$.  The definition of such a function is given as a product over the relevant conjugacy classes in $\Ga$.  We are able to study our chosen zeta function by choosing a suitable test function for the Selberg trace formula on $X_{\Ga}$ so that the geometric side gives a higher derivative of the zeta function.  We do this using the theory developed in \cite{Deitmar00}.

It is at this point that we encounter the first obstacle.  We say that an element $g\in G$ is \emph{weakly neat} if the adjoint $\Ad(g)$ has no non-trivial roots of unity as eigenvalues.  A subgroup of $G$ is weakly neat if every element is.  The results of \cite{Deitmar00} use the assumptions that the group $G$ has trivial centre and the group $\Ga$ is weakly neat.  Note that, if $G$ has trivial centre then $\Ga$ weakly neat implies $\Ga$ torsion free, since any non weakly neat torsion elements must be central.  To generalise the results of \cite{Deitmar00} to the case where $\Ga$ is not weakly neat requires two things.  Firstly, we need to modify the definition of the zeta function slightly from that given in \cite{Deitmar00}.  Secondly, the definition of the zeta function includes the first higher Euler characteristic of the spaces $X_{\Ga_{\ga}}$ for $\ga\in\CE(\Ga)$.  In \cite{Deitmar00} these are only defined for $\Ga_{\ga}$ torsion free so we need to broaden the definition.

In Chapter \ref{ch:EulerChar} we give a definition of first higher Euler characteristics which is general enough for our application.  We further prove that its value is always positive in the cases we consider.  This fact is needed in the proof of the Prime Geodesic Theorem.  In \cite{Deitmar00} it is shown that the position of the poles and zeros of the generalised Selberg zeta function depend on the Lie algebra cohomology of the irreducible, unitary representations of $\SL_4(\R)$.  Using a result of Hecht and Schmid (\cite{HechtSchmid83}) it suffices to look at the infinitesimal characters of the irreducible, unitary representations of $\SL_4(\R)$.  For this purpose, in Chapter \ref{ch:EulerChar} we also describe the unitary dual of $\SL_4(\R)$ using a result of Speh (\cite{Speh81}) and give a result about the infinitesimal characters of certain elements of this set.

In Chapter \ref{ch:Ruelle} we define the generalised Selberg zeta function and use the trace formula to deduce its analytic properties.  In particular, we give a formula for its vanishing order at a given point and prove a functional equation, from which we can deduce that it is of finite order.  We then define the generalised Ruelle zeta function and prove that it is a finite quotient of generalised Selberg zeta functions.  In particular, in the cases we are interested in it has a zero at $s=1$ and all other poles and zeros in the half plane $\Re s\leq\frac{3}{4}$, and is furthermore of finite order.  We introduce the Ruelle zeta function as its logarithmic derivative is a Dirichlet series from whose properties we can prove the Prime Geodesic Theorem.

The definitions of the generalised Selberg and Ruelle zeta functions depend on the choice of a finite dimensional virtual representation of a particular closed subgroup $M$ of $\SL_4(\R)$.  The trace of this representation at elements of $\Ga$ appears in the Dirichlet series arising as the logarithmic derivative of the Ruelle zeta function.  We can choose such a representation $\tilde{\si}$ so that its trace is zero for all non-regular elements of $\Ga$.  The factors $\nu(\CO)$ in Theorem~\ref{thm:main2} come from the traces of the representation $\tilde{\si}$.

In Chapter \ref{ch:PGT} we apply the methods of \cite{Randol77} together with standard techniques of analytic number theory to prove a Prime Geodesic Theorem.  In carrying out the application to class numbers it is necessary to show that certain subsets of the summands contribute negligibly to the asymptotic.  This is done using a ``sandwiching" argument.  To carry out this ``sandwiching" argument we define a sequence of Dirichlet series, which are not connected to any zeta function, but whose analytic properties can also be deduced from the Selberg trace formula.  We use a version of the Wiener-Ikehara theorem to prove an asymptotic result for an increasing sequence of functions derived from these Dirichlet series, which we use to ``sandwich" the product over the elements  we are interested in against the sum over all elements which comes from the Ruelle zeta function.  Unfortunately the asymptotics we are able to derive from the Wiener-Ikehara theorem do not provide an error term like those we can deduce from the Ruelle zeta function using the methods of \cite{Randol77}.  This is why we lose the error term in Theorem~\ref{thm:main1}.

In Chapter \ref{ch:DivAlg} we use a classification of the set of equivalence classes of division algebras over $\Q$ by means of a description of the Brauer group of $\Q$ (\cite{Pierce82}, Theorem 18.5) to show the existence of a division algebra $M(\Q)$ whose maximal subfields include the fields in the set $C^c(S)$ and such that $M(\Q)\ox\R\cong\Mat_4(\R)$.  In fact the maximal subfields of $M(\Q)$ are all quartic extensions of $\Q$ and the set of totally complex maximal subfields of $M(\Q)$ coincides with $C(S)$.  We further show that for an order $\CO$ in $O(S)$ the number of embeddings of $\CO$ into $M(\Q)$, up to conjugation by the unit group of the maximal order $M(\Z)$, is $\la_S(\CO)h(\CO)$.  The restriction on the set $S$ of prime numbers that it has to contain an even number of elements is a consequence of the classification of division algebras over $\Q$.

In Chapter \ref{ch:Compare} we prove a correspondence between the primitive, closed geodesics in the Prime Geodesic Theorem and the orders of totally complex quartic fields.  Under this correspondence the lengths of the geodesics correspond to the regulators of the orders.  It turns out that it is the regular, weakly neat elements of $\CE(\Ga)$ which correspond to orders in totally complex quartic fields with no real quadratic subfield.  We use the above mentioned ``sandwiching" argument to isolate these elements in the Prime Geodesic Theorem.

We finish this introduction with a few remarks about the limitations of the method used here in terms of further applications and mention a couple of other recent results in the same direction.  In order to be able to make use of the trace formula for compact spaces we have had to limit our sum over class numbers by means of the choice of a finite set of primes, as described above.  In \cite{DeitmarHoffmann05} Deitmar and Hoffmann have been able to use a different trace formula to prove that as $x\ra\infty$
$$
\sum_{R(\CO)\leq x}h(\CO)\sim\frac{e^{3x}}{3x},
$$
where the sum is over all isomorphism classes of orders in complex cubic fields.

In order to get the error term in the Prime Geodesic Theorem we have made use of the classification of the unitary dual of $\SL_4(\R)$.  At present the unitary dual is not known for any higher dimensional groups so a Prime Geodesic Theorem with error term is not possible using our methods.  Finally we mention that the correspondence between geodesics and orders actually works by identifying primitive geodesics with fundamental units in orders.  By Dirichlet's unit theorem, an order in a number field $F$ has a unique fundamental unit (up to inversion and multiplication by a root of unity) only if $F$ is real quadratic, complex cubic or totally complex quartic.  Hence an asymptotic of our form can be proven only in these three cases.  In \cite{Deitmar04} Deitmar has proved a Prime Geodesic Theorem for higher rank spaces, from which he deduces an asymptotic formula for class numbers of orders in totally real number fields of prime degree.

\mainmatter

  \chapter{Euler Characteristics and Infinitesimal Characters}
    \label{ch:EulerChar}
    \markright{\textnormal{\thechapter{. Euler Characteristics and Infinitesimal Characters}}}
    In this chapter we introduce some concepts and prove some results which will be needed for our consideration of the zeta functions in the next chapter.

\section{The unitary duals of $\SL_2(\R)$ and $\SL_4(\R)$}

For a locally compact group $G$ (or more generally a topological group) a \emph{representation} \index{representation} $(\pi,V_{\pi})$ of $G$ on a complex Hilbert space $V_{\pi}\neq 0$ is a homomorphism of $G$ into the group of bounded linear operatots on $V_{\pi}$ with bounded inverses, such that the resulting map of $G\x V_{\pi}$ into $V_{\pi}$ is continuous.  A representation will also be denoted simply as $\pi$.

An \emph{invariant subspace} \index{invariant subspace} for such a $(\pi,V_{\pi})$ is a vector subspace $U\subset V_{\pi}$ such that $\pi(g)U\subset U$ for all $g\in G$.  The representation is \emph{irreducible} \index{representation!irreducible} if it has no closed invariant subspaces other than $0$ and $V$.  The representation $\pi$ is \emph{unitary} \index{representation!unitary} if $\pi(g)$ is unitary for all $g\in G$.  Two representations $(\pi,V_{\pi})$ and $(\si,V_{\si})$ are \emph{equivalent} \index{representation!equivalent} if there is a bounded linear map $E:V_{\pi}\ra V_{\si}$ with a bounded inverse such that $\si(g)E=E\pi(g)$ for all $g\in G$.  If $\pi$ and $\si$ are unitary, they are \emph{unitarily equivalent} \index{representation!unitarily equivalent} if they are equivalent via an operator $E$ that is unitary.

The \emph{unitary dual} \index{unitary dual} $\hat{G}$ \index{$\hat{G}$} of $G$ is the set of all equivalence classes of irreducible unitary representations of $G$.  If $\pi$ is an irreducible, unitary representation of $G$ then we shall write, by slight abuse of notation, $\pi\in\hat{G}$.

For a real Lie Algebra $\g_{\R}$ a \emph{Lie algebra representation} \index{representation!Lie algebra} $(\pi,V_{\pi})$ of $\g_{\R}$ on a complex vector space $V_{\pi}\neq 0$ is a homomorphism $\pi$ of $\g_{\R}$ into the Lie algebra of all automorphisms of $V_{\pi}$.  Invariant subspaces, irreducibility and equivalence are defined in an analogous way as for group representations.

Let $G_1=\SL_2(\R)$.  We define the following subgroups: let $K_1$ be the maximal compact subgroup $\SO(2)$; let $M_1=\{\pm I_2\}$;
$$
A_1=\left\{\matrixtwo{a}{a^{-1}}:a>0\right\};
$$
and
$$
N_1=\left\{\matrix{1}{x}{}{1}:x\in\R\right\}.
$$
Let $\g_1$, $\k_1$, $\a_1$, $\n_1$ be the respective complexified Lie algebras.  Let $\rho_1\in\a_1^*$ be defined by
$$
\rho_1\matrixtwo{a}{-a}=a,
$$
and let $P_1=M_1 A_1 N_1$ be the parabolic subgroup of $G_1$ with split torus $A_1$ and unipotent radical $N_1$.

For $n\geq 2$ we define the discrete series representations \index{representation!discrete series} $\CD_n^+$ \index{$\CD_n^\pm$} and $\CD_n^-$ of $G_1$.  The Hilbert space for $\CD_n^+$ is the space of analytic functions $f$ on the upper half plane such that
$$
\|f\|^2=\int_0^{\infty}\int_{-\infty}^{\infty}|f(x+iy)|^2 y^{n-2}\ dx\,dy <\infty.
$$
The inner product is given by
$$
\langle f,g\rangle=\int_0^{\infty}\int_{-\infty}^{\infty}f(x+iy)\overline{g(x+iy)}y^{n-2}\ dx\,dy
$$
and the group action is
$$
\CD_n^+\matrix{a}{b}{c}{d}f(z)=(-bz+d)^{-n}f\left(\frac{az-c}{-bz+d}\right).
$$
The representations $\CD_n^-$ are defined analogously for functions on the lower half plane.  We also define the two limit of discrete series representations \index{representation!limit of discrete series} $\CD_1^+$ and $\CD_1^-$ on the spaces of analytic functions $f$ on the upper (respectively, lower) half plane such that
$$
\|f\|^2=\sup_{y>0}\int_{-\infty}^{\infty}|f(x+iy)|^2\ dx<\infty.
$$
The inner product is given by
$$
\langle f,g\rangle=\sup_{y>0}\int_{-\infty}^{\infty}f(x+iy)\overline{g(x+iy)}\ dx
$$
and the group action is as for the discrete series representations.  The representations $\CD_n^+$ and $\CD_n^-$ are irreducible unitary representations of $G_1$ for all $n\geq 1$.

Given an irreducible, unitary representation $\tau$ of $M_1$ and an element $\nu$ of $\a^*$ we define the \emph{induced representation} \index{representation!induced} $\Ind^{G_1}_{P_1}(\tau\otimes\nu)$ \index{$\Ind_P^G$} of $G_1=\SL_2(\R)$.  Consider the space
$$
\{f:G_1\ra V_{\tau}\ |\ f\ {\rm continuous},\ f(xman)=e^{-(\nu+\rho_1)\log a}\tau(m)^{-1}f(x)\}
$$
with the inner product
$$
\langle f,g\rangle=\int_{K_1} \langle f(k),g(k)\rangle_{\tau}\ dk.
$$
$G$ acts on the completion of this space via $\Ind^{G_1}_{P_1}(\tau\otimes\nu)f(x)=f(g^{-1}x)$, to give a representation, which is irreducible.

Let $\tau^+$ be the trivial representation on $M_1$ and $\tau^-$ the representation defined by $\tau^-(\pm I_2)=\pm 1$.  For $\nu\in\a_1^*$ we define
$$
\CP^{\pm,\nu}=\Ind_{P_1}^{G_1}(\tau^{\pm}\ox\nu).
$$
For $x\in\R$ define the \emph{principal series representations}\index{representation!principal series}
$$
\CP^{\pm,ix}=\CP^{\pm,ix\rho_1}.
$$
Then the representations $\CP^{\pm,ix}$ \index{$\CP^{\pm,ix}$} are unitary and are all irreducible except for $\CP^{-,0}\cong\CD_1^+\oplus\CD_1^-$.  The only equivalences among the representations $\CP^{\pm,\nu}$ are that $\CP^{+,ix}$ and $\CP^{-,ix}$ are equivalent to $\CP^{+,-ix}$ and $\CP^{-,-ix}$ respectively for all $x\in\R$.

We also define the \emph{complementary series representations}\index{representation!complementary series}
$$
\CC^x=\CP^{+,x\rho_1},
$$
for $0<x<1$.  Instead of the usual inner product we give $\CC^x$ \index{$\CC^x$} the inner product:
$$
\langle f,g\rangle=\int_{-\infty}^{\infty}\int_{-\infty}^{\infty}\frac{\langle f(x),g(y)\rangle_{\tau^+}}{|x-y|^{1-x}}\ dx\,dy.
$$
With this inner product the representations $\CC^x$ are unitary for all $0<x<1$.

We have the following classification theorem:

\begin{theorem}
\label{thm:SL2Dual}
\index{unitary dual!of $\SL_2(\R)$}
The unitary dual of $\SL_2(\R)$ consists of

a) the trivial representation;

b) the principal series representations $\CP^{+,ix}$ for $x\in\R$ and $\CP^{-,ix}$ for $x\in\R\smallsetminus\{0\}$;

c) the complementary series $\CC^x$ for $0<x<1$;

d) the discrete series $\CD_n^+$ and $\CD_n^-$ for $n\geq 2$ and limits of discrete series $\CD_1^+$ and $\CD_1^-$.

The only equivalences among these representations are $\CP^{\pm,ix}\cong\CP^{\pm,-ix}$ for all $x\in\R$.

\end{theorem}
\prf
\cite{Knapp86}, Theorem 16.3.
\qed

Now let $G=\SL_4(\R)$ and $K=\SO(4)$, so $K$ is a maximal compact subgroup of $G$.  Let $P'=M'A'N'$ be a parabolic subgroup of $G$ with split component $A'$ and unipotent radical $N'$.  Let $\g$ and $\a'$ be the complexified Lie algebras of $G$ and $A'$ and let $\a'^*$ be the complex dual of $\a'$.  Let $\rho'$ be the half sum of the positive roots of the system $(\g',\a')$.  Then we can define induced representations in an entirely analogous way to that used for for $G_1=\SL_2(\R)$ above.  For an irreducible, unitary representation $\tau$ of $M'$ and for $\nu\in\a'^*$ we write the corresponding induced representation of $G$ as $\Ind_{P'}^G(\tau\ox\nu)$.

If $\nu=ix\rho'$ for some $x\in\R$ then the induced representation $\Ind_{P'}^G(\tau\ox\nu)$ is unitary with respect to the inner product
$$
\langle f,g\rangle=\int_{K} \langle f(k),g(k)\rangle_{\tau}\ dk.
$$
In this case we call $\Ind_{P'}^G(\tau\ox\nu)$ a principal series representation.  For other choices is $\nu\in\a'^*$ it may be possible to make $\Ind_{P'}^G(\tau\ox\nu)$ unitary with respect to a different inner product, in which case we call $\Ind_{P'}^G(\tau\ox\nu)$ a complementary series representation.  There are also certain irreducible, unitary subrepresentations of induced representations (see \cite{Speh81}, p121) which are called \emph{limit of complementary series representations}\index{representation!limit of complementary series}.  These limit of complementary series representations can often be realised as induced representations from a parabolic subgroup $P''\supset P'$.

We define the following subgroups of $G$.  Let
\begin{eqnarray*}
M & = & \SS\matrixtwo{\SL_2^{\pm}(\R)}{\SL_2^{\pm}(\R)} \\
  & \cong & \left\{(x,y)\in \Mat_2(\R)\times \Mat_2(\R)|\begin{array}{c}\det(x),\det(y)=\pm 1 \\ \det(x)\det(y)=1 \end{array}\right\},
\index{$M$}
\end{eqnarray*}
$$
A=\left\{ \matrixfour{a}{a}{a^{-1}}{a^{-1}}:a>0\right\},\index{$A$}
$$
and
$$
N=\matrix{I_2}{\Mat_2(\R)}{0}{I_2}.\index{$N$}
$$
Let $P$ be the parabolic subgroup of $G$ with Langlands decomposition $P=MAN$.  For $m_1, m_2\in\N$, we denote by $\bar{\pi}_{m_1,m_2}$ the representation of $M$ induced from the representation $\CD_{m_1}^+\ox\CD_{m_2}^+$ of $\SL_2(\R)\x\SL_2(\R)$.  For $m\in\N$, let $\bar{\pi}_m=\bar{\pi}_{m,m}$ and let $I_m =\Ind_P^G \left(\bar{\pi}_m \otimes \frac{1}{2} \rho_P \right)$.  The representations $I_m$ each have a unique irreducible quotient known as the Langlands quotient (see \cite{Knapp86}, Theorem 7.24).  We denote the Langlands quotient of $I_m$ by $\pi_m$.\index{$\pi_m$}

We have the following classification theorem:

\begin{theorem}
\label{thm:SL4Dual}
\index{unitary dual!of $\SL_4(\R)$}
The unitary dual of $\SL_4(\R)$ consists of

a) the trivial representation;

b) principal series representations;

c) complementary series representations $\Ind^G_P\left(\bar{\pi}_m\ox t\rho_P\right)$, for $m\in\N$ and $0<t<\rez{2}$;

d) complementary series representations induced from parabolics other than $P=MAN$;

e) limit of complementary series representations, which are irreducible unitary subrepresentations of $I_m$, for $m\in\N$;

f) the family of representations $\pi_m$, indexed by $m\in\N$.
\end{theorem}
\prf
This follows from \cite{Speh81}, Theorem 5.1, where the unitary dual of $\SL_4^{\pm}(\R)=\{X\in\Mat_4(\R):\det X=\pm 1\}$ is given.
\qed

\section{Euler characteristics}
\label{sec:EulerChars}
Let $G$ be a real reductive group and suppose that there is a finite subgroup $E$ of the centre of $G$ and a reductive and Zariski-connected linear group $\CG$ such that $G/E$ is isomorphic to a subgroup of $\CG(\R)$ of finite index.  Note these conditions are satisfied whenever $G$ is a Levi component of a connected semisimple group with finite centre.    Let $K$ be a maximal compact subgroup of $G$ and let $\Ga$ be a discrete, cocompact subgroup of $G$.  Let $X_{\Ga}$ be the locally symmetric space $\Ga\bs G/K$.  If $\Ga$ is torsion-free we define the \emph{first higher Euler characteristic} \index{first higher Euler characteristic}\index{Euler characteristic!first higher} of $\Ga$ to be
$$
\chi_1(X_{\Ga}) = \chi_1(\Ga) = \sum_{j=0}^{\dim X_{\Ga}} (-1)^{j+1}jh^j(X_{\Ga}),\index{$\chi_1(\Ga)$}
$$
where $h^j(X_{\Ga})$ is the $j$th Betti number of $X_{\Ga}$.  We want to generalise this definition to cases where the group $\Ga$ is not necessarily torsion free.  We prove below a proposition which allows us broaden the definition to all cases required by our applications.

Let $\th$ be the Cartan involution fixing $K$ pointwise.  There exists a $\th$-stable Cartan subgroup $H=AB$ of $G$, where $A$ is a connected split torus and $B\subset K$ is a Cartan of $K$.  We assume that $A$ is central in $G$.  Let $C$ denote the centre of $G$.  Then $C\subset H$ and we write $B_C$, $\Ga_C$ for $B\cap C$ and $\Ga\cap C$ respectively.  Let $G^1$ be the derived group of $G$ and let $\Ga^1=G^1\cap\Ga C$.  We note in particular that, since $G=G^1C$, we have $\Ga\subset\Ga^1C\subset\Ga^1AB$.  Let $\Ga_A=A\cap\Ga_C B_C$ \index{$\Ga_A$} be the projection of $\Ga_C$ to $A$.  Then $\Ga_A$ is discrete and cocompact in $A$ (see \cite{Wolf62}, Lemma 3.3).

Let $\g_{\R}$ be the real Lie algebra Lie($G$) with polar decomposition $\g_{\R}=\k_{\R}\oplus\p_{\R}$.  Let $b$ be a fixed nondegenerate invariant bilinear form on the Lie algebra $\g_{\R}$ which is negative definite on $\k_{\R}$ and positive definite on $\p_{\R}$.  Then the form $-b(X,\th(Y))$ is positive definite and thus defines a left invariant metric on $G$.  Unless otherwise stated all Haar measures will be normalised according to the Harish-Chandra normalisation given in \cite{HarishChandra75}, \S 7.  Note that this normalisation depends on the choice of the bilinear form $b$ on $\g$.

On the space $G/H$ we have a pseudo-Riemannian structure given by the form $b$.  The Gauss-Bonnet construction generalises to pseudo-Riemannian structures to give an Euler-Poincar$\eac$ measure $\eta$ on $G/H$.  Define a (signed) Haar measure on $G$ by
$$
\mu_{EP}=\eta\ox(\textnormal{Haar measure on }H).
$$
The \emph{Weyl group} \index{Weyl group} $W=W(G,H)$ is defined to be the quotient of the normaliser of $H$ in $G$ by the centraliser.  It is a finite group generated by elements of order two.  We define the \emph{generic Euler characteristic} \index{generic Euler characteristic}\index{Euler characteristic!generic} by
$$
\chi_{\gen}(\Ga)=\chi_{\gen}(X_{\Ga})=\frac{\mu_{EP}(\Ga\bs G)}{|W|}.\index{$\chi_{\gen}(\Ga)$}
$$

\begin{proposition}
\label{pro:ECharWellDef}
Assume $\Ga$ is torsion free, $A$ is central in $G$ of dimension one and $\Ga'\subset\Ga$ is of finite index in $\Ga$.  Then $A/\Ga_A$ acts freely on $X_{\Ga}$ and $\chi_{\gen}(X_{\Ga})=\chi_1(\Ga)\vol(A/\Ga_A)$.  It follows that
$$
\chi_1(\Ga)=\chi_1(\Ga')\frac{\left[\Ga_A:\Ga'_A\right]}{\left[\Ga:\Ga'\right]}.
$$
\end{proposition}
\prf
The group $A_{\Ga}=A/\Ga_A$ acts on $\Ga\bs G/B$ by multiplication from the right.  We claim that this action is free, i.e., that it defines a fibre bundle
$$
A_{\Ga}\ra\Ga\bs G/B\ra\Ga\bs G/H.
$$
To see this let $\Ga xaB=\Ga xB$ for some $a\in A$ and $x\in G$.  Then $a=x^{-1}\ga xb$ for some $\ga\in \Ga$ and $b\in B$.  Since $\Ga\subset\Ga^1C\subset\Ga^1AB$ we can write $\ga$ as $\ga^1 a_{\ga}b_{\ga}$, with $\ga^1\in\Ga^1$ and $a_{\ga}\in A$ and $b_{\ga}\in B_C=B\cap C$.  It follows that $a_{\ga}\in\Ga_A$.  Since $A$ is central in $G$, we can write $\ga^1=aa_{\ga}^{-1}xb^{-1}x^{-1}b_{\ga}^{-1}$.  Since $\ga^1\in G^1$ and $aa_{\ga}^{-1}\in A\subset C$, we must have $aa_{\ga}^{-1}=1$, so $a=a_{\ga}\in \Ga_{A}$, which implies the claim.

In the same way we see that we get a fibre bundle
\begin{equation}
\label{eqn:TBundle}
A_{\Ga}\ra\Ga\bs G/K\ra A\Ga\bs G/K.
\end{equation}

Let $\chi$ be the usual Euler characteristic.  From \cite{HopfSamelson41} we take the equation
$\chi(K/B)=|W|$.  Using multiplicativity of Euler numbers in the fibre bundle
$$
K/B\ra A\Ga\bs G/B\ra A\Ga\bs G/K
$$
we get
\begin{eqnarray*}
\chi_{\gen}(X_{\Ga}) & = & \vol(A/\Ga_A)\frac{\eta(\Ga\bs G/H)}{|W|} \\
  &   & \\
  & = & \vol(A/\Ga_A)\frac{\chi(\Ga\bs G/H)}{\chi(K/B)} \\
  &   & \\
  & = & \vol(A/\Ga_A)\frac{\chi(A\Ga\bs G/B)}{\chi(K/B)} \\
  &   & \\
  & = & \vol(A/\Ga_A)\chi(A\Ga\bs G/K) \\
  & = & \vol(A/\Ga_A)\chi(A\bs X_{\Ga}).
\end{eqnarray*}
It remains to show that $\chi(A\bs X_{\Ga})=\chi_1(\Ga)$.

Let $\a_{\R}$ and $\g^1_{\R}$ be the Lie algebras of $A$ and $G^1$ respectively.  Then we can write
\begin{equation}
\label{eqn:gDecomp}
\g_{\R}=\a_{\R}\oplus\g^1_{\R}\oplus\z_{\R},
\end{equation}
where $\z_{\R}$ is central in $\g_{\R}$.  Let $X$ be the bi-invariant vector field on $\Ga\bs G/K$ generating the $A_{\Ga}$ action.  Then we can consider $X$ as an element of $\a_{\R}$ under the decomposition (\ref{eqn:gDecomp}).  Considering the dual of (\ref{eqn:gDecomp}), we can identify $\a_{\R}^*$ with a subspace of $\g_{\R}^*$.  Since $\a$ is central in $\g$, a non-zero element of $\a_{\R}^*$ gives us an $A_{\Ga}$-invariant, closed 1-form $\om$ on $\Ga\bs G/K$ such that for every $\Ga gK\in\Ga\bs G/K$ we have $\om(\Ga gK)(X)\neq 0$.  Since $A_{\Ga}\cong\R/\Z$ is connected and compact, the cohomology of the de Rham complex $\Ga\bs G/K$ coincides with the cohomology of the subcomplex of $A_{\Ga}$-invariants $\Om(\Ga\bs G/K)^{A_{\Ga}}$.  Using local triviality of the bundle one sees that
$$
\Om(\Ga\bs G/K)^{A_{\Ga}}=\pi^*\Om(A\Ga\bs G/K)\oplus\pi^*\Om(A\Ga\bs G/K)\wedge\om,
$$
where $\pi^*$ denotes the projection map.  Set $C_0=\pi^*\Om(A\Ga\bs G/K)$ and $C_1=C_0\oplus C_0\wedge\om=\Om(\Ga\bs G/K)^{A_{\Ga}}$.  Then
$$
H^p(C_1)=H^p(C_0)\oplus H^{p-1}(C_0)
$$
and so
\begin{eqnarray}
\label{eqn:ECharEq}
\chi_1(C_1) & = & \sum_{p\geq 0} (-1)^{p+1} p\dim H^p(C_1) \nonumber \\
  & = & \sum_{p\geq 0} (-1)^{p+1} p\left(\dim H^p(C_0)+\dim H^{p-1}(C_0)\right) \nonumber \\
  & = & \sum_{p\geq 0} (-1)^{p+1} (p-(p+1))\dim H^p(C_0) \nonumber \\
  & = & \sum_{p\geq 0} (-1)^p \dim H^p(C_0) \nonumber \\
  & = & \chi(C_0).
\end{eqnarray}
This gives the required result.
\qed

Let $\Ga$ be a discrete, cocompact subgroup of $G$.  It is known that every arithmetic subgroup of $G$ has a torsion free subgroup of finite index (\cite{Borel69}, Proposition 17.6).  Let $\Ga'\subset\Ga$ be such a subgroup.  Suppose further that the torus $A$ is central in $G$ and of dimension one.  Define $\Ga_A$ and $\Ga'_A$ as above.

We define the first higher Euler characteristic \index{first higher Euler characteristic}\index{Euler characteristic!first higher} of $\Ga$ as
\begin{equation}
\label{eqn:fhEC}
\chi_1(X_{\Ga})=\chi_1(\Ga)=\chi_1(\Ga')\frac{\left[\Ga_A:\Ga'_A\right]}{\left[\Ga:\Ga'\right]}.\index{$\chi_1(\Ga)$}
\end{equation}
Propostition \ref{pro:ECharWellDef} shows that this is well-defined.  We note that in the case that $\Ga$ itself is torsion free, this definition of first higher Euler characteristic agrees with the one given above.

\section{Euler-Poincar\'e functions}
For the next definition we assume $G$ to be a semisimple real reductive group of inner type and we fix a maximal compact subgroup $K$.  We further assume that $G$ contains a compact Cartan subgroup.  Let $\g_{\R}$ be the real Lie algebra Lie($G$) with polar decomposition $\g_{\R}=\k_{\R}\oplus\p_{\R}$, and write $\g=\k\oplus\p$ for its complexification.  Recall that we are using Harish-Chandra's Haar measure normalisation as given in \cite{HarishChandra75}, \S 7 and this normalisation depends on the choice of an invariant bilinear form $b$ on $\g$, which for our purposes in this chapter we shall leave arbitrary.

A \emph{virtual representation} \index{representation!virtual} $\si$ of a group is a formal difference of two representations $\si=\si^+-\si^-$, which is called finite dimensional if both $\si^+$ and $\si^-$ are.  Two virtual representations $\si=\si^+-\si^-$ and $\tau=\tau^+-\tau^-$ of a group are said to be isomorphic if there is an isomorphism $\si^+\oplus\tau^-\cong\tau^+\oplus\si^-$.  The trace and determinant of a virtual representation $\si=\si^+-\si^-$ are defined by $\tr\si=\tr\si^+-\tr\si^-$ and $\det\si=\det\si^+/\det\si^-$.  The dimension of $\si$ is defined as $\dim\si=\dim\si^+-\dim\si^-$.

If $V$ is a representation space with $\Z$-grading then we shall consider it naturally as a virtual representation space by $V^+=V_{{\rm even}}$ and $V^-=V_{{\rm odd}}$.  In particular, if $V$ is a subspace of $\g$ we shall always consider the exterior product $\bwedge^*V$ as a virtual representation $\bwedge^*V=\bwedge^{{\rm even}}V-\bwedge^{{\rm odd}}V$ with respect to the adjoint representation.  We consider symmetric powers and cohomology spaces similarly.

For a smooth function $f$ on $G$ of compact support and an irreducible unitary representation $(\pi,V_{\pi})\in\hat{G}$ define the operator
$$
\pi(f)=\int_G \pi(g)f(g) dg\index{$\pi(f)$}
$$
on $V_{\pi}$.  Since $f$ is smooth and has compact support, $\pi(f)$ is of trace class.

Let $(\si,V_{\si})$ be a finite dimensional virtual representation of $G$.  An \emph{Euler-Poincar$\eacit$ function} \index{Euler-Poincar$\eac$ function} $f_{\si}$ \index{$f_{\si}$} for $\si$ is a compactly supported, smooth function on $G$ such that $f_{\si}\left(kxk^{-1}\right)=f_{\si}\left(x\right)$ for all $k\in K$ and for every irreducible unitary representation $(\pi,V_{\pi})$ of $G$ the following identity holds:
\begin{equation}
\label{eqn:EPFnId}
\tr \pi\left(f_{\si}\right)=\sum_{p=0}^{\dim(\p)} (-1)^p \dim\left( V_{\pi}\ox\bwedge^p\p\ox V_{\si}\right)^K,
\end{equation}
where the superscript $K$ denotes the subspace of $K$ invariants.  By \cite{Deitmar00}, Theorem 1.1 such functions exist.  We note that the value of such functions depends on the choice of Haar measure.  We shall assume all Euler-Poincar$\eac$ functions are given with respect to the Harish-Chandra normalisation.  In the following lemmas we prove some of their properties.

\begin{lemma}
\label{lem:EPfnrestrict}
Let $G$ denote a semisimple real reductive group of inner type, with connected component $G^0$, maximal compact subgroup $K$ and compact Cartan subgroup $T\subset K$.  Let $G^+ =TG^0$.  Further let $\si$ be a finite dimensional representation of $G$, $\si^+ =\si|_{G^+}$ and $f_{\si}$ an Euler-Poincar$\eacit$ function for $\si$ on $G$.

Then $f_{\si}|_{G^+}$ is an Euler-Poincar$\eacit$ function for $\si^+$ on $G^+$.
\end{lemma}

\prf
This is Lemma 1.5 of \cite{Deitmar00}.
\qed

\begin{lemma}
\label{lem:EPfnproduct}
Let $H$, $H_1$, $H_2$ be real reductive groups of inner type such that $H=H_1\x H_2$.  Let $\si$ be an irreducible representation of $H$.  There exist irreducible representations $\si_1, \si_2$ of $H_1, H_2$ respectively such that $\si\cong \si_1\ox\si_2$ and let $f_{\si_i}$ be an Euler-Poincar$\eacit$ function for $\si_i$ on $H_i$.

Then $f_{\si}(h_1,h_2)=f_{\si_1}(h_1)f_{\si_2}(h_2)$ is an Euler-Poincar$\eacit$ function for $\si$ on $H$.
\end{lemma}
\prf
Let $\pi\in\hat{H}$.  Then $\pi=\pi_1\ox\pi_2$ for some $\pi_i\in\hat{H}_i$.  Let $K_i$ be a maximal compact subgroup of $H_i$ and let $K=K_1\x K_2$ be a maximal compact subgroup of $H$.  Let $\h_{i,\R}=\k_{i,\R}\oplus\p_{i,\R}$ be the polar decomposition of the real Lie algebra $\h_{i,\R}$ of $H_i$ and write $\h_i=\k_i\oplus\p_i$ for its complexification.  Let $\p=\p_1\oplus\p_2$.  We note in particular that
$$
\bwedge^*\p\cong\bigoplus_{p,q}\bwedge^p\p_1\ox\bwedge^q\p_2.
$$
Then
\begin{eqnarray*}
\tr\pi(f_{\si}) & = & \tr\left(\int_{H_1}\int_{H_2}f_{\si_1}(h_1)f_{\si_2}(h_2)\left(\pi_1(h_1)\ox\pi_2(h_2)\right) dh_1 dh_2 \right) \\
  & = & \tr\left[\left(\int_{H_1}f_{\si_1}(h_1)\pi_1(h_1)dh_1\right)\ox\left(\int_{H_2}f_{\si_2}(h_2)\pi_2(h_2)dh_2\right)\right] \\
  & = & \left(\tr\pi_1(f_{\si_1})\right)\left(\tr\pi_2(f_{\si_2})\right) \\
  & = & \prod_{i=1,2}\sum_{p=0}^{\dim\p_i}(-1)^p\dim\left(V_{\pi_i}\ox\bwedge^p\p_i\ox V_{\si_i}\right)^{K_i} \\
  & = & \sum_{p=0}^{\dim\p_1}\sum_{q=0}^{\dim\p_2}(-1)^{p+q}\dim\left(V_{\pi}\ox\left(\bwedge^p\p_1\ox\bwedge^q\p_2\right)\ox V_{\si}\right)^K \\
  & = & \sum_{j=0}^{\dim\p}(-1)^j\dim\left(V_{\pi}\ox\bwedge^j\p\ox V_{\si}\right)^K,
\end{eqnarray*}
which is what was required to show.
\qed

\begin{lemma}
\label{lem:EPCentral}
Let $G$ denote a semisimple real reductive group of inner type, with maximal compact subgroup $K$ and compact Cartan subgroup $T\subset K$.  Let $g\in G$ be central.  Let $\si$ be a finite dimensional representation of $G$ and let $f_{\si},h_{\si}$ be Euler-Poincar$\eacit$ functions for $\si$ on $G$.  Then $f_{\si}(g)=h_{\si}(g)$. 
\end{lemma}

\prf
Since $g$ is central, the orbital integral can be written as
$$
\CO_g(f_{\si})=\int_{G_g\bs G}f_{\si}(x^{-1}gx)\,dx=f_{\si}(g),
$$
where $G_g$ denotes the centraliser of $g$ in $G$ (which in this case, since $g$ is central, is the whole of $G$).  Similarly $\CO_g(h_{\si})=h_{\si}(g)$.  By \cite{Deitmar00}, Proposition 1.4, the value of the orbital integral $\CO_g(f_{\si})$ of an Euler-Poincar$\eac$ function for $\si$ on $G$ depends only on $g$, not on the particular Euler-Poincar$\eac$ function chosen.  Hence, $f_{\si}(g)=\CO_g(f_{\si})=\CO_g(h_{\si})=h_{\si}(g)$, as claimed.
\qed

\begin{lemma}
\label{lem:EPfnvalue}
Let $g_1$ be an Euler-Poincar$\eacit$ function for the trivial representation on $\SL_2(\R)$.  Then $g_1(1)=g_1(-1)\in\R$.
\end{lemma}
\prf
For $v\in\R$ let $\CP^{+,iv}$ and $\CP^{-,iv}$ be principal series representations on $\SL_2(\R)$.  For $n\in\N, n\geq 2$ let $\CD^+_n$ and $\CD^-_n$ be discrete series representations on $\SL_2(\R)$ and write $\CD^{\pm}_n$ for $\CD^+_n\oplus\CD^-_n$.  By \cite{Knapp86}, Theorem 11.6 there is a constant $M>0$ such that for any compactly supported, smooth function $f$ on $\SL_2(\R)$ we have
\begin{eqnarray*}
f(1) & = & M\sum_{n=2}^{\infty}(n-1)\tr\CD^{\pm}_n(f) \\
  & & + \frac{M}{4}\int_{-\infty}^{\infty}\tr\CP^{+,iv}(f)v\tanh\left(\frac{\pi v}{2}\right)+\tr\CP^{-,iv}(f)v\coth\left(\frac{\pi v}{2}\right)dv.
\end{eqnarray*}
Lemma 1.3 of \cite{Deitmar00} tells us that $\tr\CP^{\pm,iv}(g_1)=0$ for all $v\in\R$, so we have
\begin{equation}
\label{eqn:EPSL2}
g_1(1) = M\sum_{n=2}^{\infty}(n-1)\tr\CD^{\pm}_n(g_1).
\end{equation}
By the definition of an Euler-Poincar$\eac$ function, for all $n\geq 2$
$$
\tr\CD^{\pm}_n(g_1)=\sum_{p=0}^2(-1)^p\dim\left(\CD^{\pm}_n\ox\bwedge^p\p\right)^{\SO(2)},
$$
where $\p$ is the complex Lie algebra
$$
\p=\left\{\matrix{a}{b}{b}{-a}:a,b\in\C\right\}.
$$
Hence $g_1(1)\in\R$.

We want to know for which values of $n$ the trace $\tr\CD^{\pm}_n(g_1)$ is non-zero.  For this we need to know the $\SO(2)$-types of $\CD^{\pm}_n$ and $\p$.

\begin{lemma}
\label{lem:SO2Dual}
For $l\in\Z$, define the one dimensional representation $\ep_l$ of $\SO(2)$ by
$$
\ep_l R(\th)=e^{il\th},
$$
where
$$
R(\th)=\matrix{\cos\th}{-\sin\th}{\sin\th}{\cos\th}\in \SO(2).
$$
Note that $\ep_0$ is the trivial representation, which we shall denote by $triv$.

The unitary dual of $\SO(2)$ is the set $\{\ep_l:l\in\Z\}$.
\end{lemma}
\prf
Since $\SO(2)$ is abelian all its irreducible representations are one dimensional by Schur's Lemma (\cite{Knapp86}, Proposition 1.5).  By unitarity and the fact that $I=R(0)=R(2\pi)$ we have that $\widehat{\SO(2)}=\{\ep_l:l\in\Z\}$.
\qed

By computing the adjoint action of $\SO(2)$ on $\p$ we get:
\begin{lemma}
We have the following isomorphisms of $\SO(2)$-modules:
\begin{eqnarray*}
\bwedge^0 \p & = & triv \\
\bwedge^1 \p & = & \ep_2\oplus\ep_{-2} \\
\bwedge^2 \p & = & triv.
\end{eqnarray*}
\end{lemma}
\qed

\begin{lemma}
\label{lem:DiscSer}
For $n\in\N$ we have the following isomorphism of $\SO(2)$-modules:
$$
\CD^{\pm}_n \cong \bigoplus_{{|j|\geq n}\atop{j\equiv n\,\mod 2}}\ep_j.
$$
\end{lemma}
\prf
Let $\tau_n$ be the unique $n$-dimensional representation of $\SL_2(\R)$ (see \cite{Knapp86}, Chapter II $\S$1).  It follows from the definition of $\tau_n$ that the following isomorphism of $\SO(2)$-modules holds:
$$
\tau_{n} \cong \bigoplus_{{|j|\leq n-1}\atop{j\equiv n-1\ \mod(2)}}\ep_j.
$$

Let $P_1=M_1 A_1 N_1$ be the minimal parabolic of $\SL_2(\R)$.  Then the unitary dual of $M_1\cong \{1,-1\}$ consists of two one dimensional representations, which we denote by $1=triv$ and $-1$.  We denote by $\rho_1$ the character of $A_1$
$$
\rho_1\matrixtwo{a}{a^{-1}}=a,
$$
and write $\pi_{\pm 1,n-1}$ for $\Ind_{\bar{P}}\pm 1\otimes (n-1)\rho_{\bar{P}}$.  Then we have the following exact sequences of $\SL_2(\R)$-modules (see \cite{Knapp86}, Chapter II \S5):
$$
0\rightarrow\CD^{\pm}_n\rightarrow\pi_{1,n-1}\rightarrow\tau_{n-1}\rightarrow 0,\ \ \ \ n\in\N,\ n\ {\rm even}
$$
and
$$
0\rightarrow\CD^{\pm}_n\rightarrow\pi_{-1,n-1}\rightarrow\tau_{n-1}\rightarrow 0,\ \ \ \ n\in\N,\ n\ {\rm odd}, n\neq 1
$$
and the isomorphism of $\SL_2(\R)$-modules
$$
\CD_1^{\pm}\cong\pi_{-1,0}.
$$
From the so-called \emph{compact picture} of induced representations given in \cite{Knapp86}, Chapter VII \S1, it is fairly straightforward to compute the following isomorphisms of $\SO(2)$-modules, where the sums are over all integers with the same parity:
\begin{eqnarray*}
\pi_{1,n-1} & \cong & \bigoplus_{j\ {\rm even}}\ep_j \\
\pi_{-1,n-1} & \cong & \bigoplus_{j\ {\rm odd}}\ep_j.
\end{eqnarray*}
The lemma then follows easily by combining the various elements of the proof.
\qed

From the previous two lemmas we can see that $\tr\CD^{\pm}_n(g_1)$ is non-zero only if $n=2$ and in that case $\tr\CD^{\pm}_2(g_1)=-2$, so $g_1(1)=-2M$.

It remains to show that $g_1(-1)=g_1(1)$.  Let $R_z$ be the right multiplication operator of $\SL_2(\R)$ on the space $C^{\infty}_{c}\left(\SL_2(\R)\right)$ of smooth, compactly supported functions on $\SL_2(\R)$.  That is, $R_zg(x)=g(xz)$ for all $x$, $z\in\SL_2(\R)$ and $g\in C^{\infty}_{c}\left(\SL_2(\R)\right)$.  Let $\pi$ be an irreducible, unitary representation of $\SL_2(\R)$.  The matrix $-1=-I_2$ is central in $\SL_2(\R)$ so $\pi(-1)$ commutes with $\pi(x)$ for all $x\in\SL_2(\R)$ and hence, by Schur's Lemma (\cite{Knapp86}, Proposition 1.5) is scalar.  This means that for any irreducible, unitary representation $\pi$ on $\SL_2(\R)$ we have $\tr\pi\left(R_{-1}g_1\right)=\pi(-1)\tr\pi(g_1)$.  Thus we get
\begin{eqnarray*}
g_1(-1) & = & R_{-1}g_1(1) \\
  & = & M\sum_{n=2}^{\infty}(n-1)\CD^{\pm}_n(-1)\tr\CD^{\pm}_n(g_1) \\
  & = & \CD^{\pm}_2(-1)g_1(1),
\end{eqnarray*}
where $\CD^{\pm}_2(-1)$ is a scalar, which we now compute.  For $g\in\CD^+_2$ and $z\in\C$ we have the action
$$
\CD^+_2\matrix{a}{b}{c}{d}g(z)=(-bz+d)^{-2}g\left(\frac{az-c}{-bz+d}\right).
$$
Hence
$$
\CD^+_2\matrixtwo{-1}{-1}g(z)=(-1)^{-2}g\left(\frac{-z}{-1}\right)=g(z).
$$
Similarly we get
$$
\CD^-_2\matrixtwo{-1}{-1}g(z)=g(z)
$$
so $\CD^{\pm}_2(-1)=1$ and the lemma is proved.
\qed

\section{Euler characteristics in the case of $\SL_4(\R)$}
We now consider the particular case $G=\SL_4(\R)$.  Let $K=\SO(4)$, a maximal compact subgroup of $G$ and let $\Ga$ be a discrete, cocompact subgroup of $G$.  Let $A$ be the rank one torus
$$
A=\left\{ \matrixfour{a}{a}{a^{-1}}{a^{-1}}:a>0\right\},\index{$A$}
$$
and let $B$ be the compact group
$$
B=\matrixtwo {\SO(2)}{\SO(2)}.\index{$B$}
$$
$B$ is a compact Cartan subgroup of
\begin{eqnarray*}
M & = & \SS\matrixtwo{\SL_2^{\pm}(\R)}{\SL_2^{\pm}(\R)} \\
  & \cong & \left\{(x,y)\in \Mat_2(\R)\times \Mat_2(\R)|\begin{array}{c}\det(x),\det(y)=\pm 1 \\ \det(x)\det(y)=1 \end{array}\right\}.\index{$M$}
\end{eqnarray*}
Let
$$
N=\matrix{I_2}{\Mat_2(\R)}{0}{I_2}\index{$N$}
$$
and let $P=MAN$ be a parabolic subgroup of $G$ with Levi component $MA$ and unipotent radical $N$. Let
$$
A^-=\left\{ \matrixfour{a}{a}{a^{-1}}{a^{-1}}:0<a<1\right\},\index{$A^-$}
$$
be the negative Weyl chamber in $A$ with respect to the root system given by the choice of parabolic.  Let $\CE_P(\Ga)$ \index{$\CE_P(\Ga)$} be the set of $\Ga$-conjugacy classes of elements $\ga\in\Ga$ which are conjugate in $G$ to an element of $A^-B$.  We shall use this notation for the rest of this chapter.

We say $g\in G$ is \emph{regular} \index{regular} if its centraliser is a torus and \emph{non-regular} \index{non-regular} (or \emph{singular}) \index{singular} otherwise.  Clearly, for $\ga\in\Ga$ regularity is a property of the $\Ga$-conjugacy class $[\ga]$.  Let $[\ga]\in\CE_P(\Ga)$ and write $G_{\ga}$ for the centraliser of $\ga$ in $G$.  The element $\ga$ is conjugate in $G$ to an element $a_{\ga}b_{\ga}\in A^-B$ and we define the length $l_{\ga}$ of $\ga$ to be $l_{\ga}=b(\log a_{\ga},\log a_{\ga})^{1/2}$.  It follows that if $\ga$ is regular then $G_{\ga}\cong AB$.

In the first case $K_{\ga}=B$ is a maximal compact subgroup of $G_{\ga}$, the group $B$ is then clearly a Cartan subgroup of $K_{\ga}$, the product $AB$ is a $\th$-stable Cartan subgroup of $G_{\ga}$ and $A$ is central in $G_{\ga}$.  If we let $\Ga_{\ga}=\Ga\cap G_{\ga}$ denote the centraliser of $\ga$ in $\Ga$ then $\Ga_{\ga}$ is discrete and cocompact in $G_{\ga}$.  Let $\Ga'$ be a torsion free subgroup of finite index in $\Ga$ and let $\Ga_{\ga}'=\Ga'\cap G_{\ga}$.  Then $\Ga_{\ga}'$ is a torsion free subgroup of finite index in $\Ga_{\ga}$.  We define $\Ga_{\ga,A}$, $\Ga'_{\ga,A}$ as in Section \ref{sec:EulerChars}.  The first higher Euler characteristic $\chi_1(\Ga_{\ga})$ of $\Ga_{\ga}$ is then defined as in (\ref{eqn:fhEC}) as:
\begin{equation}
\label{eqn:EulerChar}
\chi_1(\Ga_{\ga})=\frac{\left[\Ga_{\ga,A}:\Ga_{\ga,A}'\right]}{\left[\Ga_{\ga}:\Ga_{\ga}'\right]}\chi_1(\Ga_{\ga}').
\end{equation}

In the second case $K_{\ga}\cong\SS(\O(2)\x\O(2))$ is a maximal compact subgroup of $G_{\ga}$.  Furthermore, $B$ is a Cartan subgroup of $K_{\ga}$, the product $AB$ is a $\th$-stable Cartan subgroup of $G_{\ga}$ and $A$ is central in $G_{\ga}$.  The definition of $\chi_1(\Ga_{\ga})$ then proceeds exactly as above.

\begin{lemma}
\label{lem:ECharEqn}
For $\ga\in\CE_P(\Ga)$ we have that
$$
\chi_1(\Ga_{\ga})=\frac{C_{\ga}\vol(\Ga_{\ga}\bs G_{\ga})}{l_{\ga_0}},
$$
where $C_{\ga}$ is an explicit constant depending only on $\ga$, which is equal to one when $\ga$ is regular, and $\ga_0$ is the primitive geodesic underlying $\ga$.
\end{lemma}

\prf
Since $\Ga_{\ga}'$ is torsion free we may take from the second proposition in section 2.4 of \cite{Deitmar95} the equation
$$
\chi_1(\Ga_{\ga}')=\frac{C_{\ga}\vol(\Ga_{\ga}'\bs G_{\ga})}{\la_{\ga}'},
$$
where $\la_{\ga}'$ denotes the volume of the maximal compact flat in $\Ga'\bs G/K$ containing $\ga$, and $C_{\ga}$ is an explicit constant depending only on $\ga$.  We note that in \cite{Deitmar95} there is an extra factor in the equation which does not show up here.  The reason is that in \cite{Deitmar95} a differential form, not a measure, was constructed.  The value of the constant $C_{\ga}$ is given in \cite{Deitmar95} in terms of the root system of $G_{\ga}$ with respect to a Cartan subgroup.  In the case that $\ga$ is regular, as noted above we have that $G_{\ga}=AB$ and hence $G_{\ga}$ is a Cartan subgroup of itself.  It is then easy to see that value of $C_{\ga}$ in this case is one.

We denote by $\la_{\ga}$ the volume of the maximal compact flat in $\Ga\bs G/K$ containing $\ga$.  The values of $\la_{\ga}$ and $\la_{\ga}'$ are given by the volumes of the images of $\Ga_{\ga,A}\bs A$ and $\Ga_{\ga,A}'\bs A$ respectively under their embeddings into $\Ga\bs G/K$ and $\Ga'\bs G/K$ respectively.  Since $A$ is one-dimensional, $\la_{\ga}=l_{\ga_0}$.  Furthermore, from the definition of the groups $\Ga_{\ga,A}$ and $\Ga_{\ga,A}'$ it follows that
$$
\left[\Ga_{\ga,A}:\Ga_{\ga,A}'\right]=\la_{\ga}'/\la_{\ga}=\la_{\ga}'/l_{\ga_0}.
$$
Since
$$
\vol(\Ga_{\ga}'\bs G_{\ga})=\left[\Ga_{\ga}:\Ga_{\ga}'\right]\vol(\Ga_{\ga}\bs G_{\ga}),
$$
the lemma follows from (\ref{eqn:EulerChar}).
\qed

\begin{theorem}
\label{thm:ECharPos}
For all $\ga\in\CE_P(\Ga)$ the Euler characteristic $\chi_1(\Ga_{\ga})$ is positive.  If $\ga$ is regular then
$$
\chi_1(\Ga_{\ga})=\frac{\left[\Ga_{\ga,A}:\Ga_{\ga,A}'\right]}{\left[\Ga_{\ga}:\Ga_{\ga}'\right]}.
$$
\end{theorem}

\prf
If $\ga$ is regular then $G_{\ga}\cong AB$, $K_{\ga}\cong B$ and $\Ga_{\ga}'$ is a complete lattice in $G_{\ga}$ so $X_{\Ga_{\ga}'}=\Ga_{\ga}'\bs G_{\ga}/K_{\ga}\cong S^1$, the unit circle in $\C$.  The Betti numbers $h^j(S^1)$ are equal to zero for $j\neq 0,1$ and one for $j=0,1$, hence $\chi_1(\Ga_{\ga}')=1$.  The claimed value for $\chi_1(\Ga_{\ga})$ then follows immediately from (\ref{eqn:EulerChar}).

Suppose $\ga$ is not regular.  Since $\ga\in\CE_P(\Ga)$, we have that $\ga$ is conjugate in $G$ to $a_{\ga}b_{\ga}\in A^-B\subset AM$.  Let $(\si,V_{\si})$ be a finite dimensional representation of $M$ and let $K_M$ be the maximal compact subgroup $\SS(\O(2)\x\O(2))$ of $M$,  which contains the compact Cartan $B$ of $M$.  We saw above that there exist Euler-Poincar$\eac$ functions of $\si$ on $M$.  Let $f_{\si}$ be one such.  We denote by $M_{\ga}$ the centraliser of $b_{\ga}$ in $M$ and by $\CO^M_{b_{\ga}}(f_{\si})$ the orbital integral
$$
\int_{M/M_{\ga}} f_{\si}\left(xb_{\ga}x^{-1}\right) dx.
$$

From Proposition 1.4 of \cite{Deitmar00} we get the equation
$$
\CO^M_{b_{\ga}}(f_{\si})=C_{\ga}\tr\si(b_{\ga}),
$$
where $C_{\ga}$ is the constant from Lemma \ref{lem:ECharEqn}.  Together with Lemma \ref{lem:ECharEqn} this gives us
$$
\chi_1(\Ga_{\ga})
  =\frac{\CO^M_{b_{\ga}}(f_{\si})\vol(\Ga_{\ga}\bs G_{\ga})}{l_{\ga_0}\tr\si(b_{\ga})}.
$$
Choosing $\si=1$, the trivial representation of $M$, this simplifies to
\begin{equation}
\label{eqn:EChar}
\chi_1(\Ga_{\ga})=\frac{\CO^M_{b_{\ga}}(f_1)\vol(\Ga_{\ga}\bs G_{\ga})}{l_{\ga_0}}.
\end{equation}

To complete the proof of the theorem we shall show that the orbital integral $\CO^M_{b_{\ga}}(f_1)$ is positive.  In the case we are considering
$$
b_{\ga}=\matrixtwo{\pm 1}{\pm 1},
$$
where $1$ denotes the identity matrix in $\SL_2(\R)$, hence it is central and $M_{\ga}=M$ so we have simply that
\begin{equation}
\label{eqn:OrbInt}
\CO^M_{b_{\ga}}(f_1)=f_1(b_{\ga}).
\end{equation}

The group $\bar{M}\cong\SL_2(\R)\x\SL_2(\R)$ is the connected component of $M\cong\SS(\SL_2^{\pm}(\R)\x\SL_2^{\pm}(\R))$.  $M$ has a maximal compact subgroup $K_M\cong\O(2)\x\O(2)$ and compact Cartan $T_M\cong\SO(2)\x\SO(2)$.  We have that $T_M \bar{M}=\bar{M}$.  Hence by Lemma \ref{lem:EPfnrestrict}, since $f_1$ is an Euler-Poincar$\eac$ function for the trivial representation on $M$ we have that $\bar{f}_1=f_1|_{\bar{M}}$ is an Euler-Poincar$\eac$ function for the trivial representation on $\bar{M}$.

Let $g_1$, $h_1$ be Euler-Poincar$\eac$ functions of the trivial representation on $\SL_2(\R)$.  By Lemma \ref{lem:EPfnproduct} the function $\tilde{f}_1=g_1 h_1$ is an Euler-Poincar$\eac$ function of the trivial representation on $\bar{M}$.

We recall that
$$
b_{\ga}=\matrixtwo{\pm 1}{\pm 1},
$$
which is central in $M$ and deduce from Lemma \ref{lem:EPCentral} that
$$
f_1(b_{\ga})=\bar{f}_1(b_{\ga})=\tilde{f}_1(b_{\ga})=g_1(\pm 1)h_1(\pm 1).
$$
From Lemmas \ref{lem:EPCentral} and \ref{lem:EPfnvalue} it follows that
$$
g_1(1)=g_1(-1)=h_1(1)=h_1(-1)\in\R
$$
and so $f_1(b_{\ga})$ is positive.
\qed

\section{The unitary dual of $K_M$}

Let $K_M=K\cap M\cong\SS(\O(2)\x\O(2))$.\index{$K_M$}  We shall need in the proof of later results to know the unitary dual $\widehat{K_M}$ of $K_M$, which is given in the following proposition.

First we define the following one dimensional representations of $\SO(2)$ and $\SO(2)\times \SO(2)$.
$$
\ep_l R(\th)=e^{il\th},\ \ \textnormal{for all}\ l\in\Z,
$$
$$
\ep_{l,k} \matrixtwo{R(\th)}{R(\eta)} = e^{i(l\theta + k\eta)},\ \ \textnormal{for all}\ l,k\in\Z,
$$
where
$$
R(\th)=\matrix{\cos\th}{-\sin\th}{\sin\th}{\cos\th}\in \SO(2).
$$
Note that $\ep_0$ and $\ep_{0,0}$ are the trivial representation of their respective groups.

\begin{proposition}
\label{pro:KMDual}
\index{unitary dual!of $K_M$}
For $l,k\in\Z$ not both zero there are unique representations $\d_{l,k}$ of $K_M$ on $\C^2$ with
$$
\d_{l,k}|_{\SO(2)\times \SO(2)}=\ep_{l,k}\oplus\ep_{-l,-k},
$$
and
$$
\d_{l,k}\matrixfour{-1}{1}{-1}{1}(z_1,z_2)=(z_2,z_1).
$$
We can also define the representation $\d$ of $K_M$ on $\C$ by
$$
\d(X,Y)(z)=(\det X)z=(\det Y)z.
$$
We have that $\widehat{K_M}=\{triv,\d,\d_{l,k}:l,k\in\Z\ \textrm{not both zero}\}$.
\end{proposition}
\prf
By Lemma \ref{lem:SO2Dual}, we have that $\widehat{\SO(2)}=\{\ep_l:l\in\Z\}$.
In general, for two locally compact groups $H$ and $K$, the map $(\si,\tau)\mapsto\si\otimes\tau$ defines an isomorphism $\hat{H}\times\hat{K}\cong\widehat{H\times K}$.  Thus the map from $\widehat{\SO(2)}\times\widehat{\SO(2)}$ to $\widehat{\SO(2)\times \SO(2)}$ given by $(\tau_1,\tau_2)\mapsto\tau_1\otimes\tau_2$ is an isomorphism.  Hence we have that $\widehat{\SO(2)\times \SO(2)}$ is the set
$$
\{\ep_l\ox\ep_k|\,l,k\in\Z\}=\{\ep_{l,k}|\,l,k\in\Z\}.
$$

Let
$$
T = \matrixfour{-1}{1}{-1}{1}
$$
and
$$
R(\th,\eta)=\matrixtwo{R(\th)}{R(\eta)}.
$$

For $l,k\in\Z$ we note that
$$
\d_{l,k}(T)\d_{l,k}(T)(z_1,z_2)=(z_1,z_2)=\d_{l,k}(T^2)(z_1,z_2)
$$
and
\begin{eqnarray*}
\d_{l,k}(T)\d_{l,k}(R(\th,\eta))\d_{l,k}(T)(z_1,z_2) & = & (e^{-i(l\th+k\eta)}z_1,e^{i(l\th+k\eta)}z_2) \\
    & = & \d_l(R(-\th,-\eta))(z_1,z_2) \\
    & = & \d_l(T R(\th,\eta) T)(z_1,z_2),
\end{eqnarray*}
so $\d_{l,k}$ is indeed a representation.  We shall show that the representations given in the proposition are in fact the only irreducible unitary representations of $K_M$.

Let $\tau\in\widehat{K_M}$.  Then by \cite{Knapp86}, Theorem 1.12(d) the representation $\tau$ restricted to ${\SO(2)\x\SO(2)}$ is a direct sum of irreducible representations, that is
$$
\tau|_{\SO(2)\x\SO(2)}=\ep_{l_1,k_1}\oplus\ep_{l_2,k_2}\oplus\cdots\oplus\ep_{l_n,k_n},
$$
for some $l_1,\dots,l_n,k_1,\dots,k_n\in\Z$.  Let $v\in\tau$ be an $(\SO(2)\x\SO(2))$-eigenvector with eigenvalue $e^{i(l_1\th+k_1\eta)}$.  Then the equation
$$
TR(-\th,-\eta)=R(\th,\eta)T
$$
tells us that $\tau(T)v$ is also an $\SO(2)$-eigenvector, with eigenvalue $e^{-i(l_1\th+k_1\eta)}$.  If $\tau(T)v$ is a scalar multiple of $v$ then $l_1=k_1=0$ and the equation $T^2=I$ tells us that $\tau=triv$ or $\d$.  Otherwise $\tau|_{\SO(2)\x\SO(2)}=\ep_{l_1,k_1}\oplus\ep_{-l_1,-k_1}$ and $\tau=\d_{l_1,k_1}$.

For $l=k=0$ we have $\d_{0,0}\cong triv\oplus\d$.  For all other values of $l$ and $k$ the representation $\d_{l,k}$ is irreducible.  Also, for $l,k\in\Z$ we have that $\d_{l,k}$ is unitarily equivalent to $\d_{-l,-k}$.  There are no other equivalences between the representations $\d_{l,k}$.

The proposition follows.
\qed

\section{Infinitesimal characters}

Let $G$ be a connected reductive group with maximal compact subgroup $K$, real Lie algebra $\g_{\R}$ and complexified Lie algebra $\g$.  Let $\pi$ be a representation of $G$.  By \cite{Knapp86}, Theorem 1.12(d), the restriction of $\pi$ to $K$ splits into an orthogonal sum of irreducible representations of $K$.  The irreducible $K$-representations occuring in this decomposition are called the \emph{$K$ types} \index{$K$ types} occuring in $\pi$.    By \cite{Knapp86}, Theorem 1.12(e), the $K$ types of $\pi$ occur in the above decomposition with well definied multiplicity, which is either a non-negative integer or $+\infty$.  We say that $\pi$ is \emph{admissible} \index{representation!admissible} if $\pi(K)$ acts by unitary operators and if each $K$ type occurs with finite multiplicity in $\pi|_K$.  In particular $\pi$ is admissible for all $\pi\in\hat{G}$ (\cite{Knapp86}, Chapter VIII, Theorem 8.1).  We say a vector in the representation space $\pi$ is \emph{$K$-finite} \index{$K$-finite} if its $K$-translates span a finite dimensional space.  Let $\pi_K$ \index{$\pi_K$} denote the space of $K$-finite vectors in $\pi$.

For $X\in\g_{\R}$ and $v\in\pi_K$ define
$$
\pi(X)v=\lim_{t\ra 0}\frac{\pi(\exp\,tX)v-v}{t}.
$$
By \cite{Knapp86}, Propositions 3.9 and 8.5 the limit exists and defines a linear map $\pi_K\ra\pi_K$.  Thus we get a representation of $\g_{\R}$ on $\pi_K$, which extends to a representation of the universal enveloping algebra $U(\g)$ of $\g$ on $\pi_K$.  This representation we also denote by $\pi$.

Let $\h$ be a Cartan subalgebra of $\g$.  The Weyl group \index{Weyl group} $W=W(\g,\h)$ is a finite group generated by the root reflections in the root system $(\g,\h)$.  The group $W$ acts on $\h$ and hence on the dual space $\h^*$.

Suppose that $\pi$ is an admissible representation of $G$ such that $\pi(Z)$ acts as a scalar on $\pi_K$ for all $Z$ in the centre $\z_G$ of the universal enveloping algebra of $\g$.  In particular this condition is satisfied when $\pi$ is irreducible admissible (\cite{Knapp86}, Chapter VIII, Corollary 8.14).  In this situation the representation $\pi$ gives us a character of $\z_G$.  Via the Harish-Chandra homomorphism (see \cite{Knapp86}, Chapter VIII, \S\S5,6) the set of characters of $\z_G$ can be identified with the set of Weyl group orbits in $\h^*$.  If $\La\in\h^*$ corresponds under this identification to the character of $\z_G$ given by $\pi$ we say that $\pi$ has \emph{infinitesimal character} \index{infinitesimal character} $\La$ and we write $\La_{\pi}=\La$.  The infinitesimal character $\La_{\pi}$ is defined up to the action of the Weyl group.  It follows from the definition of the infinitesimal character that if $\tau$ is a sub- or quotient representation of $\pi$ then $\La_{\tau}=\La_{\pi}$.

Up until now we have assumed that $G$ is connected, however if we merely assume that $G$ is a Lie group and $\pi$ is an irreducible representation of $G$, then by \cite{Knapp86}, Corollary 3.12, we have that $\pi(Z)$ commutes with $\pi(g)$ for all $Z\in\z_G$ and all $g\in G$.  Then by Schur's Lemma (\cite{Knapp86}, Proposition 1.5) we see that $\pi(Z)$ is a scalar operator for all $Z\in\z_G$, so we can define the infinitesimal character of $\pi$ in this case also.

Let $G=\SL_4(\R)$ and the subgroups $K$ and $P=MAN$ be as above.  Let $\h$ be the diagonal subalgebra of $\g$.  In this case the Weyl group $W(\g,\h)$ acts on $\h$ by interchanging elements of the leading diagonal.  We shall see in the next chapter that the analytic properties of the zeta functions considered there are related to the infinitesimal characters of elements of $\hat{G}$.  The following proposition gives us information about these infinitesimal characters which will be required in the following chapter.

Let $\rho_P$ be the half-sum of the positive roots of the system $(\g,\a)$, where $\a$ is the complexified Lie algebra of $A$, so that
$$
\rho_P \matrixfour{a}{a}{-a}{-a} = 4a.
$$
Let $\si$ be a finite dimensional virtual representation of $M$, whose $K_M$ types are all contained in the set $\{triv,\d,\d_{l,k}:l,k\in\{0,2\}\}$.  Let $\hat{M}_{\si}$ be the subset of all $\xi\in\hat{M}$ such that $\tr\xi(f_{\si})=0$ for all Euler-Poincar$\eac$ functions $f_{\si}$ for $\si$ on $M$.  (Note that the value of $\tr\xi(f_{\si})=0$ depends only on $\xi$ and $\si$ and not on the choice of Euler-Poincar$\eac$ function $f_{\si}$.)  Let $\hat{G}_{\si}$ be the set of all elements of $\hat{G}$ except for: the trivial representation; representations induced from parabolic subgroups other than $P=MAN$ and representations induced from $\xi\in\hat{M}_{\si}$.  We define an order on the real dual space of $\a$ by $\la>\mu$ if and only if $(\la-\mu)=t\rho_P$ for some $t>0$.

\begin{proposition}
\label{pro:InfChars}
Let $\si$ be a finite dimensional virtual representation of $M$ and let $\pi\in\hat{G}_{\si}$.  Then the infinitesimal character $\La_{\pi}$ of $\pi$ satisfies
$$
\Re(w\La_{\pi})|_{\a}\geq-\rez{2}\rho_P\ \textrm{ or }\ -\rho_P\geq\Re(w\La_{\pi})|_{\a}
$$
for all $w\in W(\g,\h)$.
\end{proposition}
\prf
We shall prove the proposition by considering different cases in turn.  By Theorem \ref{thm:SL4Dual} we have the following cases to consider.  The case when $\pi$ is a principal series representation induced from $P=MAN$; the case when $\pi$ is a complementary series representation induced from $P=MAN$; the case when $\pi$ is a limit of complementary series representation and the case when $\pi$ is one of the representations $\pi_m$ for $m\in\N$.

First we consider the principal series representations.  Let $\pi=\pi_{\xi,\nu}=\Ind_P^G(\xi\ox\nu)$ be induced from $P$, where $\xi$ is an irreducible, unitary representation of $M \cong \SS(\SL_2^{\pm}(\R)\times \SL_2^{\pm}(\R))$ and $\nu\in\a^*$ is imaginary (ie. $\nu\in i\a_{\R}^*$, where $\a_{\R}^*$ is the real dual space of $\a$).

Let $\xi'$ be an irreducible subspace of $\xi |_{\SL_2(\R)\times \SL_2(\R)}$.  Then $\xi'$ has infinitesimal character $\La_{\xi'}$ and $\La_{\xi}=\La_{\xi'}$.  We have that $\xi'\cong\xi_1 \otimes \xi_2$ where $\xi_1$ and $\xi_2$ are irreducible unitary representations of $\SL_2(\R)$.  

To limit the possibilities for $\xi_1$ and $\xi_2$ that we need to consider we use the double induction formula (see \cite{Knapp86}, (7.4)).
\begin{lemma}(Double induction formula)

Let $M_{\diamond}A_{\diamond}N_{\diamond}$ be a parabolic subgroup of $M$, so that $M_{\diamond}(A_{\diamond}A)(N_{\diamond}N)$ is a parabolic subgroup of $G$.  If $\si_{\diamond}$ is a unitary representation of $M_{\diamond}$ and $\nu_{\diamond}\in\a_{\diamond}^*=(\Lie_{\C}A_{\diamond})^*$, $\nu\in\a^*$, then there is a canonical equivalence of representations
$$
\Ind_{MAN}^G\left(\Ind_{M_{\diamond}A_{\diamond}N_{\diamond}}^M\left(\si_{\diamond}\ox\nu_{\diamond}\right)\ox\nu\right)\cong\Ind_{M_{\diamond}(A_{\diamond}A)(N_{\diamond}N)}^G\left(\si_{\diamond}\ox(\nu_{\diamond}+\nu)\right).
$$
\end{lemma}
\qed

We may assume that neither $\xi_1$ nor $\xi_2$ are induced since then, by the double induction formula, we would have the case that $\pi$ was induced from a parabolic other than $P=MAN$, which was excluded.  By Theorem \ref{thm:SL2Dual}, the remaining possibilities for $\xi_1$ and $\xi_2$ are the trivial representation and the discrete series and limit of discrete series representations.

Let $\La_{\pi}$ be the infinitesimal character of $\pi$ and $\La_{\xi_{i}}$ be the infinitesimal character of $\xi_{i}$, then $\La_{\xi}=\La_{\xi_1}+\La_{\xi_2}$.  Recall that $\h$ is the diagonal subalgebra of $\g$.  We lift $\La_{\xi}$ and $\nu$ to $\h$ by defining $\La_{\xi}$ to be zero on $\a$ and $\nu$ to be zero on the diagonal elements of $\m$, so that $\La_{\pi} = \La_{\xi}+\nu$ (\cite{Knapp86}, Proposition 8.22).  The Weyl group $W=W(\g,\h)$ acts on $\La_{\pi}$.  Let $w\in W$, we will show that either $\tr\xi(f_{\si})=0$ or
$$
\Re(w\La_{\pi} |_{\a})\geq-\frac{1}{2}\rho_P
$$
or
$$
-\rho_P\geq\Re(w\La_{\pi} |_{\a}).
$$

First we take $\xi_1$ and $\xi_2$ to be the trivial representation.  This gives us that
$$
\La_{\xi} \matrixfour {s}{-s}{t}{-t} = s+t.
$$
If
$$
\nu \matrixfour {a}{a}{-a}{-a} = \al a ,\ \ \ \ \al\in\textit{i}\R,
$$
then
\begin{eqnarray}
\label{eqn:InfChar}
\lefteqn{\La_{\pi}\matrixfour{a}{b}{c}{-a-b-c}} \nonumber \\
  & = & \La_{\xi}\matrixfour{\frac{a-b}{2}}{\frac{b-a}{2}}{c+\frac{a+b}{2}}{-c-\frac{a+b}{2}} + \al\frac{a+b}{2} \\
  & = & \frac{a-b}{2} + \left(c+\frac{a+b}{2}\right) + \al\frac{a+b}{2} \nonumber \\
  & = & a+c+\frac{\al}{2}(a+b). \nonumber
\end{eqnarray}

If we let $w=1$ then from above we see that $\Re(\La_{\pi}|_{\a})=0$.  If we take $w$ to be the transposition interchanging $b$ and $c$ we get
$$
w\La_{\pi} \matrixfour {a}{b}{c}{-a-b-c} = a+b+\frac{\al}{2}(a+c),
$$
the real part of which when restricted to $\a$ gives $\frac{1}{2} \rho_P$.  All other Weyl group elements are dealt with similarly and we see that $-\frac{1}{2}\rho_P\leq\Re(\La_{\pi}|_{\a})\leq\frac{1}{2} \rho_P$ for all $w\in W$.

It remains to consider the case when either or both of $\xi_1$ and $\xi_2$ are a discrete series representation or limit of discrete series representation with parameter $m_i$.  We set $\CD_0^+=\CD_0^- = triv$ and let $\xi_i = \CD_{m_i}^+$ or $\xi_i = \CD_{m_i}^-$, where $m_i \geq 0$, and $m_1$ and $m_2$ are not both zero.  From (\ref{eqn:EPFnId}) we get
\begin{equation}
\label{eqn:TraceXif1}
\tr\xi(f_{\si})=\sum_{p=0}^{\dim \p_M}\ (-1)^p\ \dim \left(V_{\xi}\ox \bwedge^p \p_M \ox V_{\si}\right)^{K_M}.
\end{equation}
We shall use Proposition \ref{pro:KMDual} to examine the $K_M$ types to limit the possibilities for $m_1$, $m_2$ for which $\tr\xi(f_{\si})\neq 0$.

\begin{lemma}
\label{lem:XiKMTypes}
For $m_1$, $m_2$ both non-zero we have the following isomorphism of $K_M$-modules:
$$
V_{\xi}\cong \bigoplus_{{j\geq m_1, j\equiv m_1\mod 2} \atop {k\geq m_2, k\equiv m_2\mod 2}}\d_{j,k}\oplus\d_{-j,k}.
$$
If $m_1=0$ we have
$$
V_{\xi}\cong \bigoplus_{k\geq m_2, k\equiv m_2\mod 2}\d_{0,k},
$$
and for $m_2=0$ we have
$$
V_{\xi}\cong \bigoplus_{j\geq m_1, j\equiv m_1\mod 2}\d_{j,0}.
$$
\end{lemma}
\prf
This follows from Lemma \ref{lem:DiscSer}.
\qed

\begin{lemma}
\label{lem:pMKMTypes1}
We have the following isomorphisms of $K_M$-modules:
\begin{eqnarray*}
  \bwedge^0 \p_M & \cong & triv \\
  \bwedge^1 \p_M & \cong & \d_{2,0}\oplus\d_{0,2} \\
  \bwedge^2 \p_M & \cong & \d\oplus\d\oplus\d_{2,2}\oplus\d_{2,-2} \\
  \bwedge^3 \p_M & \cong & \d_{2,0}\oplus\d_{0,2} \\
  \bwedge^4 \p_M & \cong & triv.
\end{eqnarray*}
\end{lemma}
\prf
$K_M$ acts on $\p_M$ by the adjoint representation and we can compute
$$
\p_M \cong \d_{2,0}\oplus\d_{0,2}.
$$
The other isomorphisms follow straightforwardly from this.
\qed

Let $v=v_1\ox v_2\ox v_3\in V_{\xi}\ox\bwedge^*\p_M\ox V_{\si}$, where the $v_i$'s all lie in a single $K_M$ type of their respective representation spaces.  Lemma~\ref{lem:pMKMTypes1} tells us that $v_2$ is in one of the following $K_M$ types: $triv$, $\d$, $\d_{2,0}$, $\d_{0,2}$, $\d_{2,2}$, $\d_{2,-2}$ by our assumption on $\si$, the possibilities for the $K_M$ type containing $v_3$ are also the same.  It follows that $\tr\xi(f_{\si})$ is non-zero only if $m_1,m_2\in\{0,2,4\}$.

It follows from the exact sequences in the proof of Lemma \ref{lem:DiscSer} that $\CD^+_m\oplus\CD^-_m\subset\CP^{\pm,(m-1)\rho_1}$, where the index on $\CP$ is $+$ if $m$ is even and $-$ if $m$ is odd.  The infinitesimal character of $\CP^{\pm,(m-1)\rho_1}$ is simply $(m-1)\rho_1$ (see \cite{Knapp86}, Proposition 8.22), hence it follows that the infinitesimal characters of $\CD^+_m$ and $\CD^-_m$ are both equal to $(m-1)\rho_1$.  This gives us in the cases where $\tr\xi(f_{\si})$ is non-zero:
\begin{equation}
\label{eqn:infchar}
\La_{\xi} \matrixfour {s}{-s}{t}{-t} = (m_1 -1)s+(m_2 -1)t,\ \ m_1,m_2\in\{0,2,4\}.
\end{equation}
By computing the action of the different Weyl group elements we can see that in all cases either $w\La_{\xi}|_{\a}\geq-\rez{2}\rho_P$ or $w\La_{\xi}|_{\a}\leq-\rho_P$.  Since $\La_{\pi}=\La_{\xi}+\nu$ and $\nu$ is imaginary we see that $\Re(\La_{\pi}|_{\a})=\La_{\xi}|_{\a}$, so the claim follows.

The complementary series representations are dealt with similarly.  We recall from Theorem \ref{thm:SL4Dual} that the complementary series induced from $P$ are $\pi=\Ind^G_P\left(\bar{\pi}_m\ox t\rho_P\right)$, for $m\in\N$ and $0<t<\rez{2}$.  We may argue as above to find the possibilities for $m$ such that $\tr\bar{\pi}_m(f_{\si})\neq 0$.  In this way we see that there are two possibilities for $m$ for which $\tr\bar{\pi}_m(f_{\si})\neq 0$, namely $m=2$ and $m=4$.  In the first case, $w\La_{\pi}|_{\a}\geq -\rez{2}\rho_P$, for all $w\in W(\g,\h)$.  In the second case we have
$$
\La_{\pi}\matrixfour{a}{b}{c}{-a-b-c}=(3+2t)a+2tb+3c.
$$
If $w\in W(\g,\h)$ is the element which swaps the first and fourth diagonal entries then $w\La_{\pi}|_{\a}=-\frac{3}{2}\rho_P$.  In all other cases we have $w\La_{\pi}|_{\a}\geq -\rez{2}\rho_P$.

The limit of complementary series representations are closely related to the family of representations $\pi_m$, $m\in\N$, so we shall deal with them together.

For $m\in\N$, we denote by $\bar{\pi}_m$ the representation of $M$ induced from the representation $\CD_m^+\otimes\CD_m^+$ of $\SL_2(\R)\x\SL_2(\R)$.  For $m\in\N$ we have the limit of complementary series representation given as an irreducible unitary subrepresentation of $I_m=\Ind_P^G(\bar{\pi}_m\ox\rez{2}\rho)$, which we will denote by $\pi_m^c$.  The representations $\pi_m$ are the Langlands quotients of the $I_m$.

Let $\La_{\pi_m}$ be the infinitesimal character of $\pi_m$ and $\La_{\pi_m^c}$ the infinitesimal character of $\pi_m^c$.  Clearly $\La_{\pi_m}=\La_{\pi_m^c}=\La_{\bar{\pi}_m}+\rez{2}\rho_P$, so we need only consider $\La_{\pi_m}$.  The value of
$$
\La_{\pi_m}\matrixfour{a}{b}{c}{-a-b-c}
$$
is equal to
$$
\La_{\bar{\pi}_m}\matrixfour{\frac{a-b}{2}}{\frac{b-a}{2}}{\frac{a+b+2c}{2}}{\frac{a+b-2c}{2}} + \rez{2}\rho_P \matrixfour {\frac{a+b}{2}}{\frac{a+b}{2}}{-\frac{a+b}{2}}{-\frac{a+b}{2}}
$$
\begin{eqnarray*}
  & = & (m-1)\frac{a-b}{2} + (m-1)(c+\frac{a+b}{2}) + (a+b) \\
  & = & ma+b+(m-1)c.
\end{eqnarray*}

Restricting to $\a$ we get
$$
\La_{\pi_m}\matrixfour{a}{a}{-a}{-a} = 2a = \frac{1}{2}\rho_P.
$$
The Weyl group $W$ acts on $\La_{\pi_m}$.  Now for all elements $w\in W$ except one we get $w\La_{\pi_m}|_{\a}\geq-\frac{1}{2}\rho_P$ for all $m\geq 1$.  The exception is $w_1$ which swaps the first and fourth diagonal entries.  For this we have
$$
w_1\La_{\pi_m}\matrixfour{a}{a}{-a}{-a} = -2(m-1)a.
$$
Thus $w_1\La_{\pi_m}|_{\a}\geq -\frac{1}{2}\rho_P$ if $m=1,2$ and if $m\geq 3$ then $w_1\La_{\pi_m}|_{\a}\leq-\rho_P$.

This completes the proof of the proposition.
\qed

  \chapter{Analysis of the Ruelle Zeta Function}
    \label{ch:Ruelle}
    \markright{\textnormal{\thechapter{. Analysis of the Ruelle Zeta Function}}}
    Let $G=\SL_4(\R)$ \index{$G$} and $\Ga\subset G$ \index{$\Ga$} be discrete and cocompact.  Let $\g_{\R}=\sl_4(\R)$ and $\g=\sl_4(\C)$ \index{$\g$} be respectively the Lie algebra and complexified Lie algebra of $G$.  As in the previous chapter, all Haar measures will be normalised as in \cite{HarishChandra75}, \S 7.  Recall that this normalisation depends on the choice of an invariant bilinear form $b$ on $\g$.  Let $b$ be the following multiple of the trace form on $\g$:
\begin{equation}
\label{eqn:norm}
\index{$b(X,Y)$}
b(X,Y)=16\tr(XY).
\end{equation}
We choose this normalisation in order to get the first zero of the Ruelle zeta function at the point $s=1$ in Theorem \ref{thm:RuelleMain} below.  Let $K \subset G$ \index{$K$} be the maximal compact subgroup $\SO(4)$.  Let $\k_{\R} \subset \g_{\R}$ be its Lie algebra and let $\p_{\R} \subset \g_{\R}$ be the orthogonal space of $\k_{\R}$ with respect to the form $b$.  Then $b$ is positive definite and Ad($K$)-invariant on $\p_{\R}$ and thus defines a $G$-invariant metric on $X=G/K$, the symmetric space attached to $G$.

Let
$$
A=\left\{ \matrixfour{a}{a}{a^{-1}}{a^{-1}}:a>0\right\},\index{$A$}
$$
$$
B=\matrixtwo {\SO(2)}{\SO(2)}.\index{$B$}
$$
$B$ is a compact Cartan subgroup of
\begin{eqnarray*}
M & = & \SS\matrixtwo{\SL_2^{\pm}(\R)}{\SL_2^{\pm}(\R)} \\
  & \cong & \left\{(x,y)\in \Mat_2(\R)\times \Mat_2(\R)|\begin{array}{c}\det(x),\det(y)=\pm 1 \\ \det(x)\det(y)=1 \end{array}\right\}.
\index{$M$}
\end{eqnarray*}
We also define
$$
N=\matrix{I_2}{\Mat_2(\R)}{0}{I_2}\textrm{ and }\ \bar{N}=\matrix{I_2}{0}{\Mat_2(\R)}{I_2},\index{$N$}\index{$\bar{N}$}
$$
and set $K_M=K\cap M$.\index{$K_M$}

Let $\m$ denote the complexified Lie algebra of $M$ and let $\m=\k_M\oplus\p_M$ \index{$\p_M$} be its polar decomposition, where $\k_M$ is the complexified Lie algebra of $K_M=K\cap M$.  Let $\h$ \index{$\h$} be the Cartan subalgebra of $\g$ consisting of all diagonal matrices, and let $\a$\index{$\a$}, $\n$ \index{$\n$} and $\bar{\n}$ \index{$\bar{\n}$} be the  complexified Lie algebras of $A$, $N$ and $\bar{N}$ respectively.  Let $P$ denote the parabolic with Langlands decomposition $P=MAN$ \index{$P$} and $\bar{P}$ the opposite parabolic with Langlands decomposition $\bar{P}=MA\bar{N}$.  Let $\rho_P$ be the half-sum of the positive roots of the system $(\g,\a)$, so that
$$
\rho_P \matrixfour{a}{a}{-a}{-a} = 4a.\index{$\rho_P$}
$$

Let
$$
H_1=\frac{1}{8}\matrixfour{-1}{-1}{1}{1} \in \a_{\R}.\index{$H_1$}
$$
Then it follows that $b(H_1)=b(H_1,H_1)=1$ and $\rho_P(H_1)=-\rez{2}$.  Let $A^- = \{\exp(tH_1):t>0\}$ \index{$A^-$} be the negative Weyl chamber in $A$.  Let $\CE_P(\Ga)$ \index{$\CE_P(\Ga)$} be the set of all conjugacy classes $[\ga]$ in $\Ga$ such that $\ga$ is conjugate in G to an element $a_{\ga}b_{\ga}$ of $A^-B$.

An element $\ga\in\Ga$ is called \emph{primitive} if for $\d\in\Ga$ and $n\in\N$ the equation $\d^n =\ga$ implies that $n=1$.  Clearly primitivity is invariant under conjugacy, so we may view it as a property of conjugacy classes.  Let $\CE_P^p(\Ga)\subset\CE_P(\Ga)$ \index{$\CE_P^p(\Ga)$} be the subset of primitive classes.

For $[\ga]\in\CE_P(\Ga)$ we define the length of $\ga$ to be $l_{\ga}=b(\log a_{\ga},\log a_{\ga})^{1/2}$\index{$l_{\ga}$}.  Let $G_{\ga}$ \index{$G_{\ga}$} be the centraliser of $\ga$ in $G$, let $\Ga_{\ga}=\Ga\cap G_{\ga}$ \index{$\Ga_{\ga}$} be the centraliser of $\ga$ in $\Ga$ and let $\chi_1(\Ga_{\ga})$ be the first higher Euler characteristic as in the previous chapter.

\section{The Selberg trace formula}
\label{sec:STF}

Let $\hat{G}$ be the unitary dual of $G$.  The group $G$ acts on the Hilbert space $L^2 (\Ga\bs G)$ by translations from the right.  Since $\Ga\bs G$ is compact this representation decomposes discretely:
$$
L^2 (\Ga\bs G)=\bigoplus_{\pi\in\hat{G}} N_{\Ga}(\pi)\pi
$$
with finite multiplicities $N_{\Ga}(\pi)$\index{$N_{\Ga}(\pi)$}, (see \cite{Gelfand69}).

Recall that for $\ga\in\Ga$ we denote by $G_{\ga}$ and $\Ga_{\ga}$ the centraliser of $\ga$ in $G$ and $\Ga$ respectively.  For a function $f$ on $G$, denote by $\CO_{\ga}(f)$\index{$\CO_{\ga}(f)$} the orbital integral\index{orbital integral}
$$
\CO_{\ga}(f)=\int_{G/G_{\ga}}f(x\ga x^{-1})\,dx.
$$

The \emph{Selberg trace formula} \index{Selberg trace formula}(see \cite{Wallach76}, Theorem 2.1) says that for suitable functions $f$ on $G$ the following identity holds:
\begin{equation}
\label{eqn:STF}
\sum_{\pi\in\hat{G}}N_{\Ga}(\pi)\tr\pi(f)=\sum_{[\ga]}\vol(\Ga_{\ga}\bs G_{\ga})\CO_{\ga}(f),
\end{equation}
where the sum on the right is over all conjugacy classes in $\Ga$.  The set of suitable functions includes, but is not limited to, all $\dim G+1$ times continuously differentiable functions of compact support on $G$.  In fact, we shall need to extend the set of test functions for which the trace formula is valid.

\begin{lemma}
Let $d\in\N$, $d\geq 16$, that is $d=2d'$ for some $d'>\dim G/2$.  Let $f$ be a $d$-times differentiable function on $G$ such that $Df\in L^1(G)$ for all left invariant differentiable operators $D$ on $G$ with complex coefficients and of degree $\leq d$.  Then the trace formula is valid for $f$.
\end{lemma}
\prf
This is Lemma 1.3 of \cite{Deitmar04} in the case $G=\SL_4(\R)$.
\qed

For $g\in G$ and $V$ any complex vector space on which $g$ acts linearly let $E(g|V)\subset\R_+$ be defined by
$$
E(g|V)=\{|\mu|:\mu\textrm{ is an eigenvalue of }g\textrm{ on }V\}.
$$
Let $\la_{\min}(g|V)=\min(E(g|V))$ and $\la_{\max}(g|V)=\max(E(g|V))$.

For $am\in AM$ define
$$
\la(am)=\frac{\la_{\min}(a|\n)}{\la_{\max}(m|\g)^2},
$$
where we are considering the adjoint action of $G$ on $\n$ and $\g$ resepectively.  Define the set
$$
(AM)^{\sim}=\{am\in AM:\la(am)<1\}.
$$
An element of $G$ is said to be \emph{elliptic} \index{elliptic element} if it is conjugate to an element of a compact torus.  Let $M_{\ell}$ denote the set of elliptic elements in $M$.

\begin{lemma}
The set $(AM)^{\sim}$ has the following properties:

(1) $A^-M_{\ell}\subset(AM)^{\sim}$;

(2) $am\in(AM)^{\sim}\Ra a\in A^-$;

(3) $am,a'm'\in(AM)^{\sim},g\in G$ with $a'm'=gamg^{-1}\Ra a=a',g\in AM$.
\end{lemma}
\prf
See \cite{Deitmar00}, Lemma 2.4.
\qed

For the construction of our test function we shall need, for given $s\in\C$ and $j\in\N$, a smooth, conjugation invariant, $j$-times continuously differentiable function on $AM$, with support in $(AM)^{\sim}$.  Further, we require that at each point $ab\in A^-B$ the function takes the value $l_a^{j+1}e^{-sl_a}$.  In \cite{Deitmar00} a function $g_s^j$ \index{$g_s^j$} is constructed with these properties, with the one difference that there the positive Weyl chamber $A^+$ is used instead of the negative Weyl chamber $A^-$.  With only very minor modification the construction in \cite{Deitmar00} will yield a function with our required properties, which we shall also call $g_s^j$. 

Let $\eta:N\ra\R$ be a smooth, non-negative function of compact support, which is invariant under conjugation by elements of $K_M$ and satisfies
$$
\int_N\eta(n)\ dn=1.
$$
Let $f:M\ra\C$ be a smooth, compactly supported function, invariant under conjugation by $K_M$.  Suppose further that the orbital integrals of $f$ on $M$ satisfy
$$
\CO_m^M(f)=\int_{M/M_m} f(xmx^{-1})\ dx=0
$$
whenever $m$ is not conjugate to an element of $B$, where $M_m$ denotes the centraliser of $m$ in $M$.  The group $AM$ acts on $\n$ according to the adjoint representation.

Given these data we define $\Phi=\Phi_{f,j,s}:G\ra\C$ by
\begin{equation}
\label{eqn:TestFunction}
\index{$\Phi$}
\Phi(knam(kn)^{-1})=\eta(n)f(m)\frac{g_s^j(am)}{\det(1-(am)^{-1}|\bar{\n})},
\end{equation}
for $k\in K,n\in N,am\in AM$.

To see that $\Phi$ is well defined we recall that, by the decomposition $G=KNAM$, every $g\in G$ which is conjugate to an element of $(AM)^{\sim}$ can be written in the form $knam(kn)^{-1}$.  By the properties of $(AM)^{\sim}$ we see that two of these representations can only differ by an element of $K_M$.  The components of the function $\Phi$ are all invariant under conjugation by $K_M$, and we note that $\det(1-(am)^{-1}|\n)\neq 0$ for all $am\in(AM)^{\sim}$, so we can see that $\Phi$ is well-defined.

\begin{proposition}
\label{pro:STF2}
The function $\Phi$ is $(j-14)$-times continuously differentiable.  For $j$ and $\Re(s)$ large enough it goes into the trace formula and we have:
\begin{equation}
\label{eqn:STF2}
\sum_{\pi\in\hat{G}} N_{\Ga}(\pi)\tr\pi(\Phi)=\sum_{[\ga]\in\CE_P(\Ga)}\vol(\Ga_{\ga}\bs G_{\ga})\CO^M_{b_{\ga}}(f)\frac{l_{a_{\ga}}^{j+1}e^{-sl_{a_{\ga}}}}{\det(1-(a_{\ga}b_{\ga})^{-1}|\bar{\n})}.
\end{equation}
\end{proposition}
\prf
Noting that $2\dim\n+\dim\k=14$ we see that this follows from Proposition 2.5 of \cite{Deitmar00}.  This proposition was proved in the case that the function $f$ is an Euler-Poincar$\eac$ function for some finite dimensional representation of $M$.  However the only properties of Euler-Poincar$\eac$ functions used in the proof are those given above for $f$, namely that it is smooth, of compact support, invariant under conjugation by $K$, and the orbital integrals satisfy the given condition.
\qed

\section{The Selberg zeta function}
\label{sec:SelbergZetaFn}

An element $g\in G$ is said to be \emph{weakly neat} \index{weakly neat} if the adjoint $\Ad(g)$ has no non-trivial roots of unity as eigenvalues.  A subgroup $H$ of $G$ is said to be weakly neat if every element of $H$ is weakly neat.  Let $[\ga]\in\CE^p_P(\Ga)$ so that $\ga$ is conjugate in $G$ to $a_{\ga}b_{\ga}\in A^-B$.  We want to know which roots of unity can occur as eigenvalues of $\Ad(\ga)$.  These are the same as the roots of unity which occur as eigenvalues of $\Ad(a_{\ga}b_{\ga})$.  If
$$
b_{\ga}=\matrixtwo{R(\th)}{R(\phi)},
$$
where
$$
R(\th)=\matrix{\cos\th}{-\sin\th}{\sin\th}{\cos\th},\index{$R(\th)$}
$$
then the eigenvalues of $\Ad(a_{\ga}b_{\ga})$ are $e^{\pm2i\th}$ and $e^{\pm2i\phi}$.  Define the sets
$$
R_{\th}=\{n\in\N:n\th\in\pi\Z\}\ \textrm{ and }\ R_{\phi}=\{n\in\N:n\phi\in\pi\Z\}
$$
and
$$
R_{\ga}=\{\min R_{\th},\min R_{\phi}\}.\index{$R_{\ga}$}
$$
Then $R_{\ga}$ contains either $0,1$ or $2$ elements.  We can see that $\ga$ weakly neat is equivalent to $R_{\ga}=\varnothing$.  For $\ga\in\CE_P(\Ga)$, where $\ga=\ga_0^n$ for $\ga_0$ primitive, the value of $\chi_1(\Ga_{\ga})$ depends on whether or not $n\in R_{\ga_0}$.

For $I\subset R_{\ga}$ with $I\neq\varnothing$ define $n_I$ \index{$n_I$} to be the least common multiple of the elements of $I$ and set $n_{\varnothing}=1$.  Further, define
$$
\chi_I(\ga)=\frac{(-1)^{|I|}}{n_I}\sum_{J\subset I}(-1)^{|J|}\chi_1(\Ga_{\ga^{n_J}}).\index{$\chi_I(\ga)$}
$$

Let $z\in\C\smallsetminus\{0\}$ and $q\in\Q$.  We define $z^q$ to be equal to $e^{q\log z}$, where we take the branch of the logarithm with imaginary part in the interval $(-\pi,\pi]$.

For any finite-dimensional virtual representation $\si$ of $M$ we define, for $\Re(s)\gg 0$, the \emph{generalised Selberg zeta function}\index{generalised Selberg zeta function}\index{Selberg zeta function, generalised}
\begin{equation}
\label{eqn:GenSelb}
\index{$Z_{P,\si}(s)$}
Z_{P,\si}(s)= \prod_{[\ga]\in\CE_P^p(\Ga)}\prod_{n\geq 0}\prod_{I\subset R_{\ga}}\det\left(1-e^{-sn_Il_{\ga}}\ga^{n_I}\,|V_{\si}\ox S^n(\n)\right)^{\chi_I(\ga)},
\end{equation}
where $S^n (\n)$ denotes the $n^{{\rm th}}$ symmetric power of the space $\n$ and $\ga$ acts on $V_{\si}\ox S^n(\n)$ via $\si(b_{\ga})\ox \Ad^n(a_{\ga}b_{\ga})$.  In the case that $\Ga$ is weakly neat this simplifies to
$$
Z_{P,\si}(s)= \prod_{[\ga]\in\CE_P^p(\Ga)}\prod_{n\geq 0}\det\left(1-e^{-sl_{\ga}}\ga\,|V_{\si}\ox S^n(\n)\right)^{\chi_1(\Ga_{\ga})},
$$
(see \cite{Deitmar00}).

Let $\pi\in\hat{G}$.  We recall that a vector in the representation space $\pi$ is said to be $K$-finite if its $K$-translates span a finite dimensional space and let $\pi_K$ denote the space of $K$-finite vectors in $\pi$.  We say that a real or complex vector space $V$ is a \emph{($\g$,$K$)-module} \index{($\g$,$K$)-module} if it is a $\g$-module and a locally finite and semi-simple $K$-module.  Further, the actions of $\g$ and $K$ must satisfy,
$$
k\cdot X\cdot v = \textnormal{Ad}k(X)\cdot k\cdot v,
$$
for $v\in V, k\in K$ and $X\in U(\g)$, where $U(\g)$ is the universal enveloping algebra of $\g$, and we must have that the action of $K$ is differentiable on every $K$-stable finite dimensional subspace of $V$ and has $\pi|_{\k}$ as its differential.

The space $\pi_K$ is a $(\g,K)$-module (see \cite{BorelWallach80}, I, 2.2).  The Lie algebra $\n$ acts on $\pi_K$ and we denote by $H^{\bullet}(\n,\pi_K)$ \index{$H^{\bullet}(\n,\pi_K)$} (resp. $H_{\bullet}(\n,\pi_K)$) the corresponding Lie algebra cohomology (resp. homology) (see \cite{BorelWallach80},\cite{CartanEilenberg56}).  We have the following isomorphism of $AM$-modules (see \cite{HechtSchmid83}, p57):
\begin{equation}
\label{eqn:AM-iso}
H_p(\n,\pi_K)\cong H^{4-p}(\n,\pi_K)\otimes \bwedge^4 \n.
\end{equation}

For $\la\in\a^*$ let
\begin{equation}
\label{eqn:MLamdaPi}
\index{$m_{\la}(\pi)$}
m_{\la}(\pi)=\sum_{p,q\geq 0} (-1)^{p+q} \dim\left( H^q(\n, \pi_K)^{\la}\otimes\bwedge^p \p_M \otimes V_{\si}\right)^{K_M},
\end{equation}
where for an $\a$-module V and $\la\in\a^*$, $V^{\la}$ \index{$V^{\la}$} denotes the generalised $\la$-eigenspace $\{v\in V|\ \exists n\in\N\ {\rm such\ that}\ (a-\la(a))^n v = 0\ \forall a\in\a\}$ and the superscript $K_M$ denotes the subspace of $K_M$ invariants.

We say that an admissible representation $\pi$ of a linear connected reductive group $G'$ has a \emph{global character} \index{global character} $\Th=\Th^{G'}_{\pi}$ \index{$\Th^G_{\pi}$} if for all smooth, compactly supported functions $f$ on $G'$ the operator $\pi(f)$ is of trace class and $\Th^{G'}_{\pi}$ is a locally $L^2$ function on $G'$ satisfying
$$
\tr\pi(f)=\int_{G'}\Th^{G'}_{\pi}(g)f(g)\,dg
$$
for each such $f$.  By \cite{Knapp86}, Theorem 10.2, every irreducible admissible, and in particular every irreducible unitary, representation of $G'$ has a global character.

\begin{theorem}
\label{thm:SelbergZeta}
Let $\si$ be a finite dimensional virtual representation of $M$.  The function $Z_{P,\si}(s)$ extends to a meromorphic function on the whole of $\C$.  For $\la\in\a^*$, the vanishing order of $Z_{P,\si}(s)$ at the point $s=\la(H_1)$ is
\begin{equation}
\label{eqn:v-order}
\sum_{\pi\in\hat{G}}N_{\Ga}(\pi)m_{\la}(\pi).
\end{equation}
Further, all the poles and zeros of $Z_{P,\si}(s)$ lie in $\R\cup(\frac{1}{2}+\textit{i}\R)$.
\end{theorem}
\prf
The analogue of this theorem in the case that $G$ has trivial centre and $\Ga$ is weakly neat is proved in \cite{Deitmar00}, Theorem 2.1.  If $G$ has trivial centre then $\Ga$ weakly neat implies $\Ga$ torsion free, since non weakly neat torsion elements must be central.  With a few modfications the proof of \cite{Deitmar00}, Theorem 2.1 becomes valid in our case also.  In fact the assumption that $\Ga$ is weakly neat is used in two places.  We sketch the proof here, pointing out the necessary modifications for it to be valid in our case.

The theorem is proved by setting $f=f_{\si}$ in the test function $\Phi=\Phi_{j,s}=\Phi_{\si,j,s}$ defined in (\ref{eqn:TestFunction}), where $f_{\si}$ is an Euler-Poincar$\eac$ function for $\si$ on $M$.  We then get that for $j\in\N$ and $\Re(s)$ sufficiently large, the right hand side of (\ref{eqn:STF2}) is equal to
\begin{equation}
\label{eqn:STFRHS}
\sum_{[\ga]\in\CE_P(\Ga)}\vol(\Ga_{\ga}\bs G_{\ga})\CO^M_{b_{\ga}}(f_{\si})\frac{l_{\ga}^{j+1}e^{-sl_{\ga}}}{\det(1-(a_{\ga}b_{\ga})^{-1}|\bar{\n})}.
\end{equation}

In \cite{Deitmar00} the following equation from \cite{Deitmar95}, proved under the assumption that $\Ga$ is torsion free, is used:
\begin{equation}
\label{eqn:ECharVol}
\chi_1(\Ga_{\ga})=\frac{C_{\ga}\vol(\Ga_{\ga}\bs G_{\ga})}{l_{\ga_0}},
\end{equation}
where $C_{\ga}$ is an explicit constant depending only on $\ga$, and $\ga_0$ denotes the primitive geodesic underlying $\ga$.  We have shown in Lemma \ref{lem:ECharEqn} that (\ref{eqn:ECharVol}) holds for all $\ga\in\CE_P(\Ga)$ in our case also.  From \cite{Deitmar00}, Proposition 1.4 we take the equation
$$
\CO^M_{b_{\ga}}(f_{\si})=C_{\ga}\tr\si(b_{\ga}),
$$
which, together with (\ref{eqn:ECharVol}) gives us
$$
\CO^M_{b_{\ga}}(f_{\si})\vol(\Ga_{\ga}\bs G_{\ga})=l_{\ga_0}\tr\si(b_{\ga})\chi_1(\Ga_{\ga}).
$$
Substituting this into (\ref{eqn:STFRHS}) we get
$$
\sum_{[\ga]\in\CE_P(\Ga)}l_{\ga_0}\tr\si(b_{\ga})\chi_1(\Ga_{\ga})\frac{l_{\ga}^{j+1}e^{-sl_{\ga}}}{\det(1-(a_{\ga}b_{\ga})^{-1}|\bar{\n})}.
$$
We claim that this is equal to
$$
(-1)^{j+1}\left(\frac{\partial\ }{\partial s}\right)^{j+2}\log Z_{P,\si}(s).
$$
Indeed
$$
\log Z_{P,\si}(s)=-\sum_{[\ga]\in\CE_P^p(\Ga)}\sum_{m=1}^\infty\sum_{I\subset R_{\ga}}\chi_I(\ga)\frac{e^{-sml_{\ga}}}{m}\tr\si(b_{\ga}^m)\sum_{n\geq 0}\tr\left((a_{\ga}b_{\ga})^m|S^n(\n)\right).
$$
We note that
\begin{eqnarray*}
\sum_{n\geq 0}\tr\left((a_{\ga}b_{\ga})^m|S^n(\n)\right) & = & \det\left(1-(a_{\ga}b_{\ga})^m|\n\right)^{-1} \\
  & = & \det\left(1-(a_{\ga}b_{\ga})^{-m}|\bar{\n}\right)^{-1}
\end{eqnarray*}
and the claimed equality follows.

In \cite{Deitmar00}, since $\Ga$ is weakly neat it follows that $G_{\ga^n}=G_{\ga}$ for all $\ga\in\Ga$ and $n\in\N$.  In \cite{Duistermaat79} it is shown that $X_{\ga}\cong\Ga_{\ga}\bs G_{\ga}/K_{\ga}$ and so it follows that for all $n\in\N$ we have $X_{\ga^n}\cong X_{\ga}$ and hence $\chi_1(\Ga_{\ga^n})=\chi_1(\Ga_{\ga})$.  Thus the above equality involving the logarithmic derivative of $Z_{P,\si}(s)$ is derived with the simpler Euler product expansion for $Z_{P,\si}(s)$ given above.

In our case we may have an element $[\ga]\in\CE_P^p(\Ga)$ with a non-trivial root of unity as an eigenvalue.  For such a $\ga$ we have $\chi_1(\Ga_{\ga^n})\neq\chi_1(\Ga_{\ga})$ for $n\in R_{\ga}$.  For this reason we have had to introduce the more complicated Euler product expansion for $Z_{P,\si}(s)$ so that the above equality involving the logarithmic derivative of $Z_{P,\si}(s)$ still holds.

On page 909 of \cite{Deitmar00} it is shown that
\begin{equation}
\label{eqn:TracePiPhi1}
\tr\pi(\Phi_s) = \int_{MA^-}f_{\si}(m)\Th^{MA}_{H^{\bullet}(\n,\pi_K)}(ma)dm\,g_s^j(a)da.
\end{equation}
Using the property (\ref{eqn:EPFnId}) of $f_{\si}$ we get
\begin{eqnarray}
\label{eqn:TracePiPhi}
\tr\pi(\Phi_s) & = & \int_{A^-}\tr\left(a\left|\left(H^{\bullet}(\n,\pi_K)\ox\bwedge^*\p_M\ox V_{\si}\right)^{K_M}\right.\right)l_a^{j+1}e^{-sl_a}\,da \nonumber \\
  & = & \int_0^{\infty}\sum_{\la\in\a^*}m_{\la}(\pi)\,e^{(\la(H_1)-s)t}\,t^{j+1}\,dt \nonumber \\
  & = & (-1)^{j+1}\left(\frac{\partial\ }{\partial s}\right)^{j+1}\sum_{\la\in\a^*}m_{\la}(\pi)\rez{s-\la(H_1)}.
\end{eqnarray}
Proposition \ref{pro:STF2} then gives us that
\begin{equation}
\label{eqn:HighLogDer}
\left(\frac{\partial\ }{\partial s}\right)^{j+2}\log Z_{P,\si}(s)=\sum_{\pi\in\hat{G}}N_{\Ga}(\pi)\left(\frac{\partial\ }{\partial s}\right)^{j+1}\sum_{\la\in\a^*}m_{\la}(\pi)\rez{s-\la(H_1)},
\end{equation}
from which it follows that the vanishing-order of $Z_{P,\si}(s)$ at $s=\la(H_1)$ is
$$
\sum_{\pi\in\hat{G}}N_\Ga(\pi) m_{\la}(\pi).
$$

Two further comments are in order.  First, in \cite{Deitmar00} the positive Weyl chamber $A^+=\{\exp(tH_1):t<0\}$ in $A$ is considered, where we have instead considered the negative chamber $A^-$.  For this reason we have interchanged the positions of the Lie algebras $\n$ and $\bar{\n}$ from the way they are used in \cite{Deitmar00}.  This is easily seen to give a precisely equivalent result.

Secondly, the results of \cite{Deitmar00} are stated for a finite dimensional representation $\si$ of $M$.  By linearity the results extend in a straightforward way to the case where $\si$ is a virtual representation, which we use here.
\qed

\section{A functional equation for $Z_{P,\si}(s)$}

\begin{proposition}
\label{pro:Bounds}
For $\la\in\a^*$ let $\|\la\|$ be the norm given by the form $b$ in \textnormal{(\ref{eqn:norm})}.  There are $m_1\in\N$ and $C>0$ such that for every $\pi\in\hat{G}$ and every $\la\in\a^*$ we have
$$
\left|\sum_{q=0}^4(-1)^q\dim(H^q(\n,\pi_K)^{\la})\right|\leq C(1+\|\la\|)^{m_1}.
$$
Further, let $S$ denote the setof all pairs $(\pi,\la)\in\hat{G}\x\a^*$ such that
$$
\sum_{q=0}^4(-1)^q\dim\left(H^q(\n,\pi_K)^{\la}\right)\neq 0.
$$
Then there is $m_2\in\N$ such that
$$
\sum_{(\pi,\la)\in S}\frac{N_{\Ga}(\pi)}{(1+\|\la\|)^{m_2}}<\infty.
$$
\end{proposition}
\prf
These results follow from \cite{Deitmar04a}, Proposition 2.4 and Lemma 2.7.
\qed

An entire function $f$ is said to be of \emph{finite order} \index{finite order, function of} if there is a positive constant $a$ and an $r>0$ such that $|f(z)|<\exp\left(|z|^a\right)$ for $|z|>r$.  If $f$ is of finite order then the \emph{order of $f$} \index{order!of a function} is defined to be the infimum of such $a$'s.

It is well known that the classical Selberg zeta function is an entire function of order two (see \cite{Selberg56}).  Our next lemma gives an analogous result for the generalised Selberg zeta function we are considering here.

Let $\si$ be a finite dimensional virtual representation of $M$.  Theorem \ref{thm:SelbergZeta} tells us that $Z_{P,\si}(s)$ is meromorphic and hence it may be written as the quotient of two entire functions:
$$
Z_{P,\si}(s)=\frac{Z_1(s)}{Z_2(s)},
$$
where the zeros of $Z_1(s)$ correspond to the zeros of $Z_{P,\si}(S)$ and the zeros of $Z_2(s)$ correspond to the poles of $Z_{P,\si}(s)$.  The orders of the zeros of $Z_1(s)$ (resp. $Z_2(s)$) equal the orders of the corresponding zeros (resp. poles) of $Z_{P,\si}(s)$.  The functions $Z_1(s)$ and $Z_2(s)$ are not uniquely determined, but clearly their respective sets of zeros, together with the orders of the zeros, are.  For $i=1,2$ let $R_i$ denote the set of zeros of $Z_i(s)$ counted with multiplicity.

\begin{lemma}
\label{lem:FiniteOrder}
There exist $Z_i(s)$, for $i=1,2$, with the above properties, which are, in addition, of finite order.
\end{lemma}

\prf
Let $Z_1(s)$, $Z_2(s)$ be as above.  We shall show that we may take them to be of finite order.

For $\pi\in\hat{G}$ let $\La(\pi)$ be the set of all $\la\in\a^*$ with $\la\neq 0$ and $m_{\la}(\pi)\neq 0$.  Let $\hat{G}(\Ga)$ be the set of $\pi\in\hat{G}$ such that $N_{\Ga}(\pi)\neq 0$ and let $S$ denote the set of all pairs $(\pi,\la)$ such that $\pi\in\hat{G}(\Ga)$ and $\la\in\La(\pi)$.  For $\la\in\a^*$ let $\|\la\|$ be the norm given by the form $b$ in (\ref{eqn:norm}).

The expression (\ref{eqn:v-order}) tells us that $s\neq 0$ is a zero or pole of $Z_P(s)$ if and only if $s=\la(H_1)$ for some $\la\in\a^*$, for which there exists $\pi\in\hat{G}(\Ga)$ such that $(\pi,\la)\in S$.

Since the function $Z_{P,\si}(s)$ is meromorphic and non-zero, it follows that there exists $c>0$ such that
\begin{equation}
\label{eqn:Bound1}
|\la(H_1)|\geq c\left(1+\|\la\|\right)
\end{equation}
for all $\la$ such that $(\pi,\la)\in S$ for some $\pi$.

By the definition of $m_{\la}(\pi)$ we deduce immediately from Proposition \ref{pro:Bounds} that there exist $m_1\in\N$ and $C>0$ such that for every $\pi\in\hat{G}$ and every $\la\in\a^*$ we have
\begin{equation}
\label{eqn:Bound2}
|m_{\la}(\pi)|\leq C\left(1+\|\la\|\right)^{m_1},
\end{equation}
and that there exists $m_2\in\N$ such that
\begin{equation}
\label{eqn:Bound3}
\sum_{(\pi,\la)\in S} \frac{N_{\Ga}(\pi)}{\left(1+\|\la\|\right)^{m_2}}<\infty.
\end{equation}
Let $k = m_1 + m_2$.  Then, for $(\pi,\la)\in S$, by (\ref{eqn:Bound1}) and (\ref{eqn:Bound2}) we have,
\begin{eqnarray*}
\left|\frac{m_{\la}(\pi)}{\la(H_1)^k}\right| & \leq & \rez{c^k}\cdot\frac{|m_{\la}(\pi)|}{\left(1+\|\la\|\right)^k} \\
  & \leq & \frac{C}{c^k}\cdot\rez{\left(1+\|\la\|\right)^{m_2}}.
\end{eqnarray*}
It then follows from (\ref{eqn:v-order}) and (\ref{eqn:Bound3}) that, for $i=1,2$,
\begin{equation}
\label{eqn:FiniteRank}
\sum_{\rho\in R_i\smallsetminus\{0\}}\rez{|\rho|^k}\leq\sum_{(\pi,\la)\in S}\frac{N_{\Ga}(\pi)|m_{\la}(\pi)|}{|\la(H_1)|^k}<\infty.
\end{equation}
We say that $Z_1(s)$ and $Z_2(s)$ are of \emph{finite rank}\index{finite rank, function of}.

By the Weierstrass Factorisation Theorem \index{Weierstrass Factorisation Theorem} (\cite{Conway78}, VII.5.14) and (\ref{eqn:FiniteRank}), there exist entire functions $g_1(s)$, $g_2(s)$ such that
\begin{equation}
\label{eqn:Weierstrass}
Z_i(s)=s^{n_i}e^{g_i(s)}\prod_{\rho\in R_i\smallsetminus\{0\}}\left(1-\frac{s}{\rho}\right)\exp\left(\frac{s}{\rho}+\frac{s^2}{2\rho^2}\cdots+\frac{s^k}{k\rho^k}\right),
\end{equation}
where $n_i$ is the order of the zero of $Z_i(s)$ at $s=0$.

From (\ref{eqn:HighLogDer}) we recall the equation
\begin{equation}
\label{eqn:LogDerEqn}
\left(\frac{d\ }{ds}\right)^j\log Z_{P,\si}(s)=\sum_{\pi\in\hat{G}(\Ga)}N_{\Ga}(\pi)\left(\frac{d\ }{ds}\right)^{j-1}\sum_{\la\in\La(\pi)\cup\{0\}}\frac{m_{\la}(\pi)}{s-\la(H_1)},
\end{equation}
where $j\in\N$ is sufficiently large.  Let $J=\max(j,k)$.  It follows from (\ref{eqn:Weierstrass}) and (\ref{eqn:LogDerEqn}) that
$$
\left(\frac{d\ }{ds}\right)^{J-1}\left(g'_1(s)-g'_2(s)\right) + \sum_{\rho\in R_1}\frac{(-1)^{J-1}(J-1)!}{(s-\rho)^J} - \sum_{\rho\in R_2}\frac{(-1)^{J-1}(J-1)!}{(s-\rho)^J}
$$
is equal to
$$
\sum_{\pi\in\hat{G}(\Ga)}\sum_{\la\in\La(\pi)\cup\{0\}}N_{\Ga}(\pi)m_{\la}(\pi)\frac{(-1)^{J-1}(J-1)!}{(s-\la(H_1))^J}.
$$

Remembering that the zeros of $Z_i(s)$ are included with multiplicity in $R_i(s)$ and bearing in mind the expression (\ref{eqn:v-order}) for the vanishing order of $Z_{P,\si}(s)$, we see that this implies that
$$
\left(\frac{d\ }{ds}\right)^{J-1}\left(g'_1(s)-g'_2(s)\right)=0.
$$
It follows that $g_1(s)-g_2(s)$ is a polynomial of degree at most $J$.  Without loss of generality we may take $g_2(s)$ to be zero, so that $g_1(s)$ is itself a polynomial of degree at most $J$.  Finally, Theorem XI.2.6 of \cite{Conway78} tells us that since $Z_1(s)$ and $Z_2(s)$ are of finite rank and $g_1(s)$ and $g_2(s)$ are both polynomials of degree at most $J$, it follows that $Z_1(s)$ and $Z_2(s)$ are both of order at most $J$.
\qed

Before we give the next lemma we make a couple of definitions.  Let $G'$ be a linear connected reductive group with maximal compact subgroup $K'$ and let $\xi$ be a representation of $G'$.  Let $\g'$ be the complexified Lie algebra of $G'$, let $\h'\subset\g'$ be the diagonal subalgebra and let $\rho'\in(\h')^*$ be the half sum of the positive roots of the system $(\g',\h')$.  We say that $\xi$ is \emph{tempered} \index{representation!tempered} if for all $K'$-finite vectors $u,v\in \xi_{K'}$ there exist constants $c_{u,v}$ such that for all $g\in G'$
$$
|\langle\xi(g)u,v\rangle|\leq c_{u,v}\int_{K'}e^{-\rho'(H(g^{-1}k))}dk,
$$
where $H(g^{-1}k)$ denotes the logarithm of the split part of $g^{-1}k$ under the Iwasawa decomposition and $\langle\cdot,\cdot\rangle$ is the inner product on $V_{\xi}$.  An admissible representation of $G'$ is called \emph{standard} \index{representation!standard} if it is induced from an irreducible tempered representation of $M'\subset G'$, where $P'=M'A'N'$ is a parabolic subgroup of $G'$.

The Weyl group $W(G,A)$ has two elements, let $w$ be the nontrivial element therein.  Then $w$ acts on $M$ by conjugation and for a representation $\si$ of $M$ we can define the representation $^w\si$ by $^w\si(m)=\si(wmw^{-1})$.

\begin{lemma}
\label{lem:Zeros}
Suppose that $\si\cong\,^w\!\si$ as $K_M$-modules.  Then $Z_{P,\si}(s)$ and $Z_{P,\si}(1-s)$ have the same poles and zeros with multiplicity.
\end{lemma}
\prf
By Theorem \ref{thm:SelbergZeta} $\la\in\a^*$, the vanishing order of $Z_{P,\si}(s)$ at the point $s=\la(H_1)$ is equal to $\sum_{\pi\in\hat{G}}N_{\Ga}(\pi)m_{\la}(\pi)$.  We shall show that
\begin{equation}
\label{eqn:MLaPiEq}
m_{\la}(\pi)=m_{-\la-2\rho_P}(\pi)
\end{equation}
for all $\pi\in\hat{G}$.  Since $-2\rho_P(H_1)=1$ the lemma will follow.

Let $\pi\in\hat{G}$.  Recall from (\ref{eqn:MLamdaPi}) that
$$
m_{\la}(\pi)=\sum_{p,q\geq 0} (-1)^{p+q} \dim\left( H^q(\n, \pi_K)^{\la}\otimes\bwedge^p \p_M \otimes V_{\si}\right)^{K_M}.
$$
The global character $\Th^G_{\pi}$ of $\pi$ on $G$ can be written as a linear combination with integer coefficients of characters of standard representations (\cite{Knapp86}, Chapter X, \S10.2).  We are interested in the values taken by $\Th^G_{\pi}$ on $MA^-$.  By \cite{Knapp86}, Proposition 10.19, the only characters which are non-zero on $MA$ are those characters of representations induced from $P=MAN$.  So there exist $n\in\N$ and for all $1\leq i\leq n$ integers $c_i$, tempered representations $\xi_i$ of $M$ and $\nu_i\in\a^*$ such that
\begin{equation}
\label{eqn:TempRepSum}
\Th^G_{\pi}=\sum_{i=1}^n c_i\Th^G_{\pi^i},
\end{equation}
where $\pi^i=\Ind_P^G(\xi_i\ox\nu_i)$.  From \cite{HechtSchmid83}, Theorem 3.6 and (\ref{eqn:AM-iso}) it follows that for all regular $ma\in MA^-$ we have
$$
\Th^G_{\pi}(ma)=\Th^{MA}_{H^{\bullet}(\n,\pi_K)}(ma)\frac{\det\left(a|\bwedge^4\n\right)}{\det(1-ma|\n)},
$$
and the same holds for all $\pi^i$.  Together with (\ref{eqn:TempRepSum}) this implies that
$$
\Th^{MA}_{H^{\bullet}(\n,\pi_K)}(ma)=\sum_{i=1}^n\Th^{MA}_{H^{\bullet}(\n,\pi^i_K)}(ma)
$$
for all regular $ma\in MA^-$.  Substituting this into (\ref{eqn:TracePiPhi1}) and proceeding as in (\ref{eqn:TracePiPhi}) we see that it suffices to show (\ref{eqn:MLaPiEq}) for the representations $\pi^i$.

The Weyl group element $w$ acts on the group $MA$ by conjugation, which has the effect of swapping the two components.  For $\nu\in\a^*$ we therefore have $w\nu=-\nu$.  Recall that for a representation $\xi$ of $M$ we let $^w\xi$ denote the representation $^w\xi(m)=\xi(wmw^{-1})$.  By \cite{Knapp86}, Theorem 8.38, the representations $\pi^i$ and $^w\pi^i=\Ind_P^G(^w\xi_i\ox(-\nu_i))$ are equivalent.  From the definition of the induced group action we can see that
$$
H^q(\n,\pi^i_K)^{\la}=H^q(\n,^w\!\pi^i_K)^{-\la-2\rho_P}.
$$
Thus we see that
$$
\dim\left( H^q(\n, \pi^i_K)^{\la}\otimes\bwedge^p \p_M \otimes V_{\si}\right)^{K_M}
$$
is equal to
$$
\dim\left( H^q(\n, ^w\!\pi^i_K)^{-\la-2\rho_P}\otimes\bwedge^p \p_M \otimes V_{\si}\right)^{K_M}.
$$

Using the notation of Proposition \ref{pro:KMDual} we note that the $K_M$-types satisfy $^wtriv_{K_M}=triv_{K_M}$, $^wdet_{K_M}=det_{K_M}$ and $^w\d_{l,m}=\d_{m,l}\cong\d_{-m,-l}$.  Lemma \ref{lem:pMKMTypes1} and the isomorphism of $K_M$-modules
$$
\n\cong\d_{2,2}\oplus\d_{2,-2}
$$
tell us that the $K_M$-modules $\bwedge^p\p_M$ and $\n$ are invariant under the action of $w$ and we have assumed that $^w\si\cong\si$ as $K_M$-modules.  Since $w^2$ is the identity element of $W(G,A)$ we can conclude that
$$
\dim\left( H^q(\n, \pi^i_K)^{\la}\otimes\bwedge^p \p_M \otimes V_{\si}\right)^{K_M}
$$
is equal to
$$
\dim\left( H^q(\n,\pi^i_K)^{-\la-2\rho_P}\otimes\bwedge^p \p_M \otimes V_{\si}\right)^{K_M}.
$$
Hence $m_{\la}(\pi^i)=m_{-\la-2\rho_P}(\pi^i)$ and the lemma follows.
\qed

\begin{theorem}
\label{thm:FunctEqn}
Suppose that $\si\cong\,^w\!\si$ as $K_M$-modules.  Then there exists a polynomial $G(s)$ such that $Z_{P,\si}(s)$ satisfies the functional equation
$$
Z_{P,\si}(s)=e^{G(s)}Z_{P,\si}(1-s).
$$
\end{theorem}
\prf
Let $R_i$, $n_i$ be as above for $i=1,2$.  We saw in the proof of Lemma \ref{lem:FiniteOrder} that there exist entire functions $Z_1(s)$ and $Z_2(s)$ of finite order such that
$$
Z_{P,\si}(s)=\frac{Z_1(s)}{Z_2(s)}
$$
and polynomials $g_1(s)$ and $g_2(s)$ such that
$$
Z_i(s)=s^{n_i}e^{g_i(s)}\prod_{\rho\in R_i\smallsetminus\{0\}}\left(1-\frac{s}{\rho}\right)\exp\left(\frac{s}{\rho}+\frac{s^2}{2\rho^2}\cdots+\frac{s^k}{k\rho^k}\right).
$$
Lemma \ref{lem:Zeros} tells us that $Z_{P,\si}(s)$ and $Z_{P,\si}(1-s)$ have the same poles and zeros.  Hence we can in the same way conclude that there exist entire functions $W_1(s)$ and $W_2(s)$ of finite order such that
$$
Z_{P,\si}(1-s)=\frac{W_1(s)}{W_2(s)}
$$
and polynomials $h_1(s)$ and $h_2(s)$ such that
$$
W_i(s)=s^{n_i}e^{h_i(s)}\prod_{\rho\in R_i\smallsetminus\{0\}}\left(1-\frac{s}{\rho}\right)\exp\left(\frac{s}{\rho}+\frac{s^2}{2\rho^2}\cdots+\frac{s^k}{k\rho^k}\right).
$$
Setting
$$
G(s)=g_1(s)+h_2(s)-g_2(s)-h_1(s)
$$
we can see that the claimed functional equation does indeed hold.
\qed

\section{The Ruelle zeta function}
\label{sec:RuelleZetaFn}
For any finite-dimensional virtual representation $\si$ of $M$ we define, for $\Re(s)\gg 0$, the \emph{generalised Ruelle zeta function}\index{generalised Ruelle zeta function}\index{Ruelle zeta function, generalised}
$$
R_{\Ga,\si}(s) = \prod_{[\ga]\in\CE_P^p(\Ga)}\prod_{I\subset R_{\ga}}\det\left(1-e^{-sn_Il_{\ga}}\ga^{n_I}\,|V_{\si}\right)^{\chi_I(\ga)}.\index{$R_{\Ga,\si}(s)$}
$$
We have the following theorem giving a relationship between the generalised Selberg zeta function and the generalised Ruelle zeta function.

\begin{theorem}
\label{thm:RuelleZeta}
Let $\si$ be a finite dimensional virtual representation of $M$.  The function $R_{\Ga,\si}(s)$ extends to a meromorphic function on $\C$.  More precisely
$$
R_{\Ga,\si}(s) = \prod_{q=0}^4 Z_{P,(\bwedge^q\bar{\n}\ox\si)}\left( s+\frac{q}{4}\right)^{(-1)^q}.
$$
\end{theorem}
\prf
Let $\ga\in\CE_P(\Ga)$ let $\mu(\ga)\in\N$ be the least such that $\ga=\ga_0^{\mu(\ga)}$ for some $\ga_0\in\CE_P^p(\Ga)$.  We compute at first
\begin{eqnarray*}
\log R_{\Ga,\si}(s) & = & \sum_{[\ga]\in\CE_P^p(\Ga)}\sum_{I\subset R_{\ga}}\chi_I(\ga)\tr\log\left(1-e^{-sn_Il_{\ga}}\ga^{n_I}\,|V_{\si}\right) \\
  & = & \sum_{[\ga]\in\CE_P^p(\Ga)}\sum_{I\subset R_{\ga}}\chi_I(\ga)\sum_{m=1}^{\infty}\rez{m}e^{-smn_Il_{\ga}}\tr\si(b_{\ga}^{n_I}) \\
  & = & \sum_{[\ga]\in\CE_P(\Ga)}\sum_{I\subset R_{\ga}}\chi_I(\ga)\frac{e^{-sn_Il_{\ga}}}{\mu(\ga)}\tr\si(b_{\ga}^{n_I}).
\end{eqnarray*}
Similarly, for $\tau_q=\bwedge^q\bar{\n}\ox \si$ we have
\begin{eqnarray*}
\log Z_{P,\tau_q}(s) & = & \sum_{[\ga]\in\CE^p_P(\Ga)}\sum_{I\subset R_{\ga}}\chi_I(\ga)\sum_{n=0}^{\infty}\tr\log\left((1-e^{-sn_Il_{\ga}}\ga^{n_I}|V_{\tau_q}\ox S^n(\n)\right) \\
  & = & \sum_{[\ga]\in\CE^p_P(\Ga)}\sum_{I\subset R_{\ga}}\chi_I(\ga)\sum_{n=0}^{\infty}\sum_{m=1}^{\infty}\frac{-1}{m}e^{-smn_Il_{\ga}}\tr\left(\ga^{mn_I}|V_{\tau_q}\ox S^n(\n)\right) \\
  & = & -\sum_{[\ga]\in\CE_P(\Ga)}\sum_{I\subset R_{\ga}}\chi_I(\ga)\frac{e^{-sn_Il_{\ga}}}{\mu(\ga)}\tr\si(b_{\ga}^{n_I})\frac{\tr(b_{\ga}^{n_I}|\bwedge^q\bar{\n})}{\det\left(1-(a_{\ga}b_{\ga})^{-n_I}|\bar{\n}\right)} \\\
  & = & -\sum_{[\ga]\in\CE_P(\Ga)}\sum_{I\subset R_{\ga}}\chi_I(\ga)\frac{e^{-sn_Il_{\ga}}}{\mu(\ga)}\tr\si(b_{\ga}^{n_I})\frac{\tr(b_{\ga}^{n_I}|\bwedge^q\bar{\n})}{\tr\left((a_{\ga}b_{\ga})^{-n_I}|\bwedge^*\bar{\n}\right)}.
\end{eqnarray*}
Since $\n$ is an $M$-module defined over the reals we conclude that the trace $\tr\left((a_{\ga}b_{\ga})^{-1}|\bwedge^*\bar{\n}\right)$ is a real number.  Therefore it equals its complex conjugate, which is $\tr\left(a_{\ga}^{-1}b_{\ga}|\bwedge^*\bar{\n}\right)$.  Summing over $q$ we get
$$
\log R_{\Ga,\si}(s)=\sum_{q=0}^4(-1)^q\log Z_{P,\tau_q}\left(s+\frac{q}{4}\right).
$$
The theorem follows.
\qed

We shall be interested in the zeta function $R_{P,\si}(s)$ in the case where $\si=triv=1$ is the trivial representation of $M$ and in the case where $\si$ is the following virtual representation of $M$.  Define
$$
\tilde{\si}^+=15\bwedge^0\m+6\bwedge^2\m+\bwedge^4\m
$$
and
$$
\tilde{\si}^-=10\bwedge^1\m+3\bwedge^3\m
$$
and let $\tilde{\si}=\tilde{\si}^+ -\tilde{\si}^-$\index{$\tilde{\si}$}, where $M$ acts on $\bwedge^n\m$ according to $\Ad^n$.  The reason for this choice of $\tilde{\si}$ is that it allows us to separate the contribution of the non-regular elements of $\CE_P(\Ga)$ from the regular ones, as made clear in the following lemma.

\begin{lemma}
\label{lem:SigmaTilde}
Let $\ga\in\CE_P(\Ga)$ with
$$
b_{\ga}=\matrixtwo{R(\th)}{R(\phi)},
$$
where
$$
R(\th)=\matrix{\cos\th}{-\sin\th}{\sin\th}{\cos\th}.
$$
Then
$$
\tr\tilde{\si}(b_{\ga})=\det(1-b_{\ga}:\m/\b)=4(1-\cos 2\th)(1-\cos 2\phi),
$$
where $\b$ is the complexified Lie algebra of the group $B$.  In particular we have $\tr\tilde{\si}(b_{\ga})\geq 0$ for all $\ga\in\CE_P(\Ga)$, and $\tr\tilde{\si}(b_{\ga})=0$ if and only if $\ga$ is non-regular.
\end{lemma}
\prf
As $B$-modules, we have the isomorphism $\m\cong\b\oplus(\m/\b)$.  The group $B$ is abelian and so acts trivially on the 2-dimensional Lie algebra $\b$.  Hence, for $0\leq n\leq 6$, we have the following isomorphism of $B$-modules
\begin{equation}
\label{eqn:BModIso}
\bwedge^n\m\cong\bigoplus_{p+q=n}\left({2}\atop{p}\right)\bwedge^q\m/\b.
\end{equation}
We set $a_0=1$, $a_1=-3$, $a_2=6$, $a_3=-10$, $a_4=15$ and note that, for $k=0,...,4$, these satisfy
\begin{equation}
\label{eqn:ASum}
\sum_{m=0}^{k}a_{k-m}\left({2}\atop{m}\right)=(-1)^k.
\end{equation}
The $B$-module isomorphism
$$
\bwedge^*(\m/\b)\cong\sum_{n=0}^4 a_{4-n}\bwedge^n\m=\tilde{\si}
$$
then follows from (\ref{eqn:BModIso}) and (\ref{eqn:ASum}).  We can then see that
$$
\tr\tilde{\si}(b_{\ga})=\tr\left(b_{\ga}:\bwedge^*\m/\b\right)=\det(1-b_{\ga}:\m/\b).
$$
The adjoint action of $B$ on $\m/\b$ can easily be computed to give the claimed value.  Finally we note that $\tr\tilde{\si}(b_{\ga})=0$ if and only if $\th,\phi\in\Z\pi$ and recall from the proof of Lemma \ref{lem:centralisers} that this occurs precisely when $\ga$ is non-regular.
\qed

The main result of this chapter is the following theorem.

\begin{theorem}
\label{thm:RuelleMain}
The function $R_{\Ga,1}(s)$ has a double zero at $s=1$.  The function $R_{\Ga,\tilde{\si}}(s)$ has a zero of order eight at $s=1$.  Apart from that, for $\si\in\{1,\tilde{\si}\}$, all poles and zeros of $R_{\Ga,\si}(s)$ are contained in the union of the interval $\left[-1,\frac{3}{4}\right]$ with the five vertical lines given by $\frac{k}{4} + i\R$ for $k=-2,-1,0,1,2$.
\end{theorem}

The theorem will follow from the following proposition, which we prove in the remainder of the chapter.

\begin{proposition}
\label{pro:SelbergPolesZeros}
The function $Z_{P,1}(s)$ has a double zero at $s=1$.  The function $Z_{P,\tilde{\si}}(s)$ has a zero of order eight at $s=1$.  Apart from that, for $\si\in\{1,\tilde{\si}\}$, all poles and zeros of $Z_{P,\si}(s)$ lie in the half-plane $\{\Re(s)\leq \frac{3}{4}\}$.  Further, for $\tau_q$ being the representation of $M$ on $\bwedge^q\bar{\n}\ (q\in\{1,\ldots,4\})$ obtained from the adjoint representation and $\si\in\{\tau_q,\tau_q\ox\tilde{\si}\}$, all poles and zeros of $Z_{P,\si}(s)$ lie in the half-plane $\{\Re(s)\leq 1\}$.
\end{proposition}

Let us assume for a moment that the proposition has been proved.  The representations $\si$ of $M$ considered in the proposition all satisfy the isomorphism of $K_M$-modules $\si\cong\,^w\!\si$, where $w$ is the non-trivial element of the Weyl group $W(G,A)$.  Hence we can apply the functional equation given in Theorem \ref{thm:FunctEqn} to see that the poles and zeros of the functions $Z_{P,\si}(s)$ all lie in the region $0\leq\Re(s)\leq 1$.  We can then apply Theorem \ref{thm:SelbergZeta} to see that the poles and zeros in fact lie in $[0,1]\cup(\rez{2}+i\R)$.  Finally, an application of Theorem \ref{thm:RuelleZeta} completes the proof of Theorem \ref{thm:RuelleMain}.

The proof of the proposition will take up the rest of the chapter.  We see from (\ref{eqn:v-order}) that a representation $\pi\in\hat{G}$ makes a contribution to the vanishing order of $Z_{P,\si}(s)$ only if $m_{\la}(\pi)\neq 0$ for some $\la\in\a^*$.  From (\ref{eqn:TracePiPhi}) it follows that, if $\tr\pi(\Phi_{s,\si})=0$ for some $\pi\in\hat{G}$ and finite dimensional representation $\si$ of $M$, then $m_{\la}(\pi)=0$ for all $\la\in\a^*$.  Thus, if we can show that $\tr\pi(\Phi_{s,\si})=0$ for some $\pi$ and $\si$, then we will know that $\pi$ makes no contribution to the vanishing order of $Z_{P,\si}$.

We define a partial order on the real dual space $\a_{\R}^*$ by setting $\mu >\nu$ if and only if $\mu-\nu=t\rho_P$ for some $t>0$.  Since $\rho_P(H_1)=-\frac{1}{2}$, we shall be able to prove Proposition \ref{pro:SelbergPolesZeros} with the following steps.  Firstly, irrespective of the choice of $\si$, all eigenvalues $\la$ of $\a$ on $H^{\bullet}(\n, \pi_K)$ satisfy
\begin{equation}
\label{eqn:CharIneq2}
\Re(\la)\geq -2\rho_P,
\end{equation}
with equality only if $\pi=triv$.  Secondly, in the case $\si\in\{1,\tilde{\si}\}$, if $\pi$ is the trivial representation on $G$, then $H^4(\n, \pi_K)^{-2\rho_P}\neq 0$.  If $\si=1$ then this gives a double zero at the point $s=1$ and if $\si=\tilde{\si}$ this gives a zero of order eight at $s=1$.  Apart from this, for all $\pi\in\hat{G}$ and all $\la\in\a^*$, we have either $\tr\pi(\Phi)=0$, or that $H^{\bullet}(\n, \pi_K)^{\la}\neq 0$ implies \begin{equation}
\label{eqn:CharIneq1}
\Re(\la)\geq -\frac{3}{2}\rho_P.
\end{equation}

The first step will be proven in the following two sections using the Hochschild-Serre spectral sequence.  Then the second step will be proven using a result of Hecht and Schmid which relates the action of $\a$ on $H^{\bullet}(\n, \pi_K)$ to the infinitesimal character of $\pi$.

\section{Spectral sequences}

In what follows we shall need to use the notion of a spectral sequence.  We give here a brief definition, for more details see \cite{Lang02}, Chapter XX,\S 9.

By a \emph{complex} we mean a sequence $(K^i,d^i)$, for $i\geq 0$, of objects and homomorphisms in a given abelian category satisfying
\begin{displaymath}
\xymatrix{K^0 \ar[r]^{d^0} & K^1 \ar[r]^{d^1} & K^2 \ar[r]^{d^2} & }
\end{displaymath}
and such that $d^{i+1}\circ d^i=0$.  We say that a complex $L=(L^i,e^i)$ is a \emph{subcomplex} of $K=(K^i,d^i)$, and we write $L\subset K$, if $L^i\subset K^i$ and $e^i=d^i|_{L^i}$ for all $i\geq 1$.  For a complex $K=(K^i,d^i)$, a \emph{filtration} $FK$ of $K$ is a decreasing sequence of subcomplexes
$$
K=F^0 K\supset F^1 K\supset F^2 K\supset\cdots\supset F^n K\supset F^{n+1} K=\{0\}.
$$
The quotient of two complexes can be formed in the obvious way.

For a complex $K=(K^i,d^i)$ we define the \emph{$i$-th cohomology group} of $K$ as $H^i(K)=\ker d^i/\Im d^{i-1}$.  The associated graded object
$$
H(K)=\bigoplus_{i\geq 0}H^i(K)
$$
is called the \emph{cohomology} of $K$.

To a filtered complex $FK$ is associated the \emph{graded complex}
$$
\Gr FK=\Gr K=\bigoplus_{i\geq 0}\Gr^i K=\bigoplus_{i\geq 0}F^i K/F^{i+1}K.
$$
A filtration $F^iK$ on $K$ induces a filtration $F^i H(K)$ on the cohomology by
$$
F^iH^j(K)=\frac{\ker (d^j|_{F^i K})}{\Im (d^{j-1}|_{F^i K})}.
$$
The associated graded object is
$$
\Gr H(K)=\bigoplus_{i,j}\Gr^i H^j(K)=\bigoplus_{i,j}\frac{F^iH^j(K)}{F^{i+1}H^j(K)}.
$$

A \emph{spectral sequence} \index{spectral sequence} is a sequence $(E_r,d_r)$, for $r\geq 0$, of bigraded objects
$$
E_r=\bigoplus_{p,q\geq 0} E^{p,q}_r
$$
together with homomorphisms (known as \emph{differentials})
$$
d_r:E^{p,q}_r\ra E^{p+r,q-r+1}_r \textrm{ satisfying } d_r^2=0,
$$
and such that the cohomology of $E_r$ is $E_{r+1}$, that is
\begin{equation}
\label{eqn:SSCohomCondition}
H(E_r)=\ker d_r/\Im d_r=E_{r+1}.
\end{equation}
If we have $E_{r_0}=E_{r_0+1}=\cdots$ for some $r_0$, then this limit object is called $E_{\infty}$ and we say that the spectral sequence \emph{abuts} to $E_{\infty}$.  If there exist $m,n\in\N$ such that for all $r$ we have $E_r^{p,q}\neq 0$ only if $p\leq m$ and $q\leq n$, then it follows from (\ref{eqn:SSCohomCondition}) that $E_2^{m,n}=E_3^{m,n}=\cdots=E_{\infty}^{m,n}$.  In \cite{Lang02}, Chapter XX, Proposition 9.3, it is shown that given a filtered complex $FK$ one can derive a spectral sequence $(E_r)$, which in particular has
$$
E_{\infty}^{p,q}=\Gr^p(H^{p+q}(K)).
$$

\section{The Hochschild-Serre spectral sequence}

We shall need the following proposition.

\begin{proposition}
\label{pro:Vanishing}
Let $P'=M'A'N'$ be a parabolic subgroup of $G$ and let $\a'$,$\n'$ be the complexified Lie algebras of $A'$,$N'$ respectively.  Define a partial order $>_{\a'}$ on the dual space $\a'^*$ of $\a'$ by $\mu>\nu$ if and only if $\mu-\nu$ is a non-zero linear combination with positive integral coefficients of positive roots of $(\g,\a)$.

Let $\la\in\a'^*$ and $0\leq p<\dim\n$ be such that $H^p(\n,\pi_K)^{\la}\neq 0$.  Then there exists $\mu\in\a'^*$ with $\la >_{\a'}\mu$ such that $H^{\dim\n}(\n,\pi_K)^{\mu}\neq 0$.
\end{proposition}
\prf
This follows from \cite{HechtSchmid83}, Proposition 2.32 and the isomorphism (\ref{eqn:AM-iso}).
\qed

Let $M_0\subset G$ be the subgroup of diagonal matrices with each diagonal entry equal to $\pm 1$, let $A_0\subset G$ be the subgroup of diagonal matrices with positive entries and let $N_0\subset G$ be the subgroup of upper triangular matrices with ones on the diagonal.  Then $P_0=M_0A_0N_0$ is the minimal parabolic subgroup of $G$.  Let $\a_0$ and $\n_0$ be the complexified Lie algebras of $A_0$ and $N_0$ respectively, and let $\a_0^*$ be the dual of $\a_0$.  Let $\rho_0\in\a_0^*$ be defined as follows:
$$
\rho_0\matrixfour{a}{b}{c}{d}=3a+2b+c.
$$
Then $\rho_0$ is the half-sum of the positive roots of the system $(\g,\a_0)$.  Let $\a_{0,\R}$ be the real Lie algebra of $A_0$, which we may consider as a subalgebra of $\a_0$.  Let $\a_{0,\R}^-$ be the set of all $X\in\a_{0,\R}$ such that $\al(X)<0$ for all positive roots $\al$ of the system $(\g,\a_0)$.  Let $\a_{0,\R}^{*,+}$ be the set of all $\nu\in\a_{0,\R}^*$ such that $\nu(X)<0$ for all $X\in\a_{0,\R}^-$.  Then for all $\nu\in\a_{0,\R}^{*,+}$ we have $\nu=\sum_{\al}\la_{\al}\al$, where the sum is over the positive roots of $(\g,\a_0)$ and every $\la_{\al}>0$.

For $\pi\in\hat{G}$ a \emph{matrix-coefficient} \index{matrix-coefficient} of $\pi$ is any function $G\ra\C$ of the form
$$
f_{u,v}(g)=\langle\pi(g)u,v\rangle
$$
for some $u,v\in\pi$.

\begin{lemma}
\label{lem:MinParAAction}
Let $\pi\in\hat{G}$ and let $\la\in\a_0^*$ be such that $H^{\bullet}(\n_0,\pi_K)^{\la}\neq 0$.  Then $\la\in-2\rho_0+\a_{0,\R}^{*,+}$, except in the case $\pi=triv$, when $H^{6}(\n_0,\pi_K)^{-2\rho_0}\neq 0$ and other than this $H^{\bullet}(\n_0,triv)^{\la}\neq 0$ implies $\la\in-2\rho_0+\a_{0,\R}^{*,+}$.
\end{lemma}
\prf
According to \cite{HechtSchmid83}, Theorem~4.16, there exists a countable set $\CEE(\pi)\subset\a_1^*$ and a collection of polynomial functions $p^{\nu}_{u,v}$ indexed by $u,v\in\pi_K$ and $\nu\in\CEE(\pi)$ such that if $f_{u,v}$, for $u,v\in\pi_K$, is a matrix coefficient of $\pi$ then:
\begin{equation}
\label{eqn:mtxcoeff}
f_{u,v}(\exp\ X)=\sum_{\nu\in\CEE(\pi)} p^{\nu}_{u,v}(X)e^{(\nu +\rho_0)(X)},
\end{equation}
for all $X\in\a_{0,\R}^-$.  We assume that $\CEE(\pi)$ is minimal and hence that each polynomial $p^{\nu}_{u,v}$ is non-zero.

If $\pi=triv$ then every matrix coefficient of $\pi$ is a constant function so in (\ref{eqn:mtxcoeff}) we have $\CEE(\pi)=\{-\rho_0\}$ and $p^{-\rho_0}_{u,v}=C_{u,v}$, where $C_{u,v}$ is a constant for all $u,v\in\pi$.

If $\pi\in\hat{G}\smallsetminus\{triv\}$ then $\pi$ is infinite dimensional so by \cite{HoweMoore79}, Theorem 5.1 all the matrix coefficients of $\pi$ vanish at infinity.  Hence, for all $\nu\in\CEE(\pi)$ and for all $X\in\a_{0,\R}^-$, we have $(\nu+\rho_0)(X)<0$.  In other words, for all $\nu\in\CEE(\pi)$ we have $\nu\in-\rho_0+\a_{0,\R}^{*,+}$.

We have assumed that $H^{\bullet}(\n_0,\pi_K)^{\la}\neq 0$.  Let $>_{\a_0}$ be the partial order defined on $\a_0^*$ by $\mu>_{\a_0}\nu$ if and only if $\mu-\nu$ is a linear combination with positive integral coefficients of positive roots of $(\g,\a_0)$.  By Proposition \ref{pro:Vanishing} there exists $\mu\in\a_0^*$ such that $H^6(\n_0,\pi_K)^{\mu}\neq 0$ and $\la\geq_{\a_0}\mu$.

Let $\La=\{\nu+\rho_0\in\a_0^*:H^6(\n_0,\pi_K)^{\nu}\neq 0\}$ and let $\La^{\min}$ be the set of elements of $\La$ which are minimal with respect to $>_{\a_0}$.  Let $\CEE^{\min}(\pi)$ be the set of $\nu\in\CEE(\pi)$ which are minimal with respect to $>_{\a_0}$.  Theorem~4.25 of \cite{HechtSchmid83} says that $\La^{\min}=\CEE(\pi)^{\min}$.

If $\pi=triv$ then $\La^{\min}=\{-\rho_0\}$ and the claims of the proposition follow from the definition of $\La$ and \cite{HechtSchmid83}, Proposition~2.32 as quoted above.

If $\pi\in\hat{G}\smallsetminus\{triv\}$ then $\La^{\min}=\CEE(\pi)^{\min}$ implies that there exists $\nu\in\CEE(\pi)$ such that $\la\geq_{\a_0}\nu-\rho_0$.  We saw above that $\nu\in-\rho_0+\a_{0,\R}^{*,+}$ so the claim follows.
\qed

\begin{proposition}
\label{pro:AAction}
Let $\pi\in\hat{G}$ and let $\la\in\a^*$ be such that $H^{\bullet}(\n,\pi_K)^{\la}\neq 0$.  Then $\Re(\la)\geq-2\rho_P$, with equality only if $\pi=triv$.
\end{proposition}
\prf
Let $\n_M=\n_0\cap\m$ and $\a_M=\a_0\cap\m$.  Then $\n_0=\n_M\oplus\n$ and $\a_0=\a\oplus\a_M$.  In \cite{HochschildSerre53b, HochschildSerre53a} a filtered complex is constructed so that the spectral sequence derived from it has
$$
E_2^{p,q} \cong H^p(\n_M,H^q(\n,\pi_K))
$$
and
$$
E_{\infty}^{p,q} \cong \Gr^p H^{p+q}(\n_0,\pi_K),
$$
where $H^{p+q}(\n_0,\pi_K)$ is appropriately filtered.  This is the so-called \emph{Hochschild-Serre spectral sequence}\index{spectral sequence!Hochschild-Serre}.

By Proposition \ref{pro:Vanishing} it suffices to prove the proposition in the case
$$
H^4(\n,\pi_K)^{\la}\neq 0.
$$
In this case it then follows that
$$
H^2(\n_M,H^4(\n,\pi_K)^{\la})\neq 0
$$
and so there exists $\la_M\in\a_M^*$ such that
$$
H^2(\n_M,H^4(\n,\pi_K)^{\la})^{\la_M}\neq 0.
$$
Since $A$ acts trivially on $\n_M$, this equals
$$
H^2(\n_M,H^4(\n,\pi_K))^{\la+\la_M}=(E_2^{2,4})^{\la+\la_M},
$$
where we consider $\la+\la_M$ as an element of $\a_0^*=\a^*\oplus\a_M^*$.  For all $r$ we have that $E_r^{p,q}\neq 0$ only if $0\leq p\leq 2$ and $0\leq q\leq 4$, so it follows that $E_2^{2,4}=E_{\infty}^{2,4}$.  Since the action of $A_0$ commutes with the differentials of the spectral sequence it follows that
$$
(E_{\infty}^{2,4})^{\la+\la_M}\neq 0
$$
and hence that
$$
H^6(\n_0,\pi_K)^{\la+\la_M}\neq 0.
$$
The proposition then follows from Lemma \ref{lem:MinParAAction} by projection of $\la+\la_M$ onto the $\a^*$ component.
\qed

\section{Contribution of the trivial representation}

For a $K_M$-module $\eta$ let $2\eta$ denote the module $\eta\oplus\eta$.  We shall need the following:
\begin{lemma}
\label{lem:pMKMTypes}
We have the following isomorphisms of $K_M$-modules:
\begin{eqnarray*}
  \bwedge^0 \p_M & \cong & triv \\
  \bwedge^1 \p_M & \cong & \d_{2,0}\oplus\d_{0,2} \\
  \bwedge^2 \p_M & \cong & 2\d\oplus\d_{2,2}\oplus\d_{2,-2} \\
  \bwedge^3 \p_M & \cong & \d_{2,0}\oplus\d_{0,2} \\
  \bwedge^4 \p_M & \cong & triv.
\end{eqnarray*}
and
\begin{eqnarray*}
\bwedge^0 \m & \cong & triv \\
\bwedge^1 \m & \cong & 2\d\oplus\d_{2,0}\oplus\d_{0,2} \\
\bwedge^2 \m & \cong & triv\oplus2\d\oplus2(\d_{2,0}\oplus\d_{0,2})\oplus\d_{2,2}\oplus\d_{2,-2} \\
\bwedge^3 \m & \cong & 2triv\oplus2triv\oplus2(\d_{2,0}\oplus\d_{0,2}\oplus\d_{2,2}\oplus\d_{2,-2}) \\
\bwedge^4 \m & \cong & triv\oplus2\d\oplus2(\d_{2,0}\oplus\d_{0,2})\oplus\d_{2,2}\oplus\d_{2,-2}
\end{eqnarray*}
\end{lemma}

\prf
The isomorphisms for $\p_M$ were given in Lemma \ref{lem:pMKMTypes1}.  For $\m$ note that $K_M$ acts on $\m$ by the adjoint representation and we can compute
$$
\m\cong 2\d\oplus\d_{2,0}\oplus\d_{0,2}.
$$
The other isomorphisms follow straightforwardly from this.
\qed

\begin{proposition}
\label{pro:TrivAAction}
Let $\pi$ be the trivial representation on $G$.

For $\si$ the trivial representation on $M$, the representation $\pi$ contributes a double zero of $Z_{P,\si}(s)$ at the point $s=1$.  For $\si=\tilde{\si}$ the representation $\pi$ contributes a zero of $Z_{P,\si}(s)$ of order eight at the point $s=1$.  In both cases, all other poles and zeros contributed by $\pi$ are in $\left\{\Re(s)\leq\frac{3}{4}\right\}$.
\end{proposition}
\prf
The space $H_0(\n,\pi_K)=\pi_K /\n\pi_K$ is one dimensional with trivial $\a$-action.  The action of $\a$ on $\n$ is given by $\rez{2}\rho_P$, so the isomorphism (\ref{eqn:AM-iso}) tells us that $\a$ acts on the one dimensional space $H^4(\n,\pi_K)$ according to $-2\rho_P$.  Proposition \ref{pro:Vanishing} tells us that for $q=0,1,2,3$ and $\la\in\a^*$, $H^q(\n,\pi_K)^{\la}\neq 0$ implies $\la\geq -\frac{3}{2}\rho_P$.  Since $-2\rho_P(H_1)=1$ this gives a pole or zero at $s=1$ and evaluation of the relevant characters at $H_1$ also gives the other poles and zeros contributed by $\pi$ in $\{\Re(s)\leq\frac{3}{4}\}$.

It remains to compute the vanishing order of $Z_{P,\si}(s)$ for $\si\in\{1,\tilde{\si}\}$ at the point $s=1$.  Since $\dim H^4(\n,\pi_K)^{-2\rho_P}=1$, we get from (\ref{eqn:v-order}) the following expression for the vanishing order:
$$
N_{\Ga}(\pi) \sum_{p\geq 0} (-1)^{p} \dim\left(\bwedge^p \p_M\ox V_{\si}\right)^{K_M}.
$$
For $\si=1$ this is equal to
$$
N_{\Ga}(\pi) \sum_{p\geq 0} (-1)^{p} \dim\left(\bwedge^p \p_M\right)^{K_M}
$$
and for $\si=\tilde{\si}$ to
$$
N_{\Ga}(\pi) \sum_{p,q\geq 0} (-1)^{p+q} \dim\left(\bwedge^p \p_M\ox a_q\bwedge^q\m\right)^{K_M},
$$
where $a_0=15$, $a_1=10$, $a_2=6$, $a_3=3$, $a_4=1$ and $a_q=0$ for $q\geq 5$.  The only functions on $L^2(\Ga\bs G)$ invariant under the action of $G$ are the constant functions, hence $N_{\Ga}(\pi)=1$.  We can then use Lemma \ref{lem:pMKMTypes} to see that the above expressions take the claimed values.
\qed

\section{Contribution of the other representations}

The following proposition gives a relationship between the action of $\a$ on $H^{\bullet}(\n,\pi_K)$ and the infinitesimal character of $\pi$.

\begin{proposition}
\label{pro:aCohomologyAction}
Let $\pi\in\hat{G}$.

Suppose $H^{\bullet}(\n, \pi_K)^{\la}\neq 0$ for some $\la\in\a^*$.  Then $\la=w\La_{\pi}|_{\a}-\rho_P$, where $w\in W(\g,\h)$ and $\La_{\pi}$ is a representative of the infinitesimal character of $\pi$.  Moreover, for $p\in\Z$ we have
$$
H^p(\n,\pi_K)=\bigoplus_{w\in W(\g,\h)} H^p(\n,\pi_K)^{w\La_{\pi}-\rho_P}.
$$
\end{proposition}
\prf
This follows from \cite{HechtSchmid83}, Corollary~3.32 and the isomorphism (\ref{eqn:AM-iso}).
\qed

In light of this proposition, in order to complete the proof of Proposition \ref{pro:SelbergPolesZeros}, and hence of Theorem \ref{thm:RuelleMain}, it will be sufficient to show that for all $\pi\in\hat{G}\smallsetminus\{triv\}$ and for $\si\in\{1,\tilde{\si}\}$ either $\tr\pi(\Phi_{\si})=0$ or the infinitesimal character $\La_{\pi}$ of $\pi$ satisfies
$$
\Re(w\La_{\pi})|_{\a}\geq-\rez{2}\rho_P\ \textrm{ or }\ -\rho_P\geq\Re(w\La_{\pi})|_{\a}
$$
for all $w\in W(\g,\h)$.

We first consider the elements of $\hat{G}$ induced from parabolic subgroups other that $P=MAN$.

We say that two parabolics $P'=M'A'N'$ and $P''=M''A''N''$, where $A', A''\subset A_0 =\left\{\diag(a,b,c,(abc)^{-1})|a,b,c>0\right\}$ are \emph{associate} if there exists a member $w$ of the Weyl group $W(\g,\h)$ such that $w^{-1}MAw=M'A'$.  Representations of $G$ induced from associate parabolics are equivalent.

Up to association G has four parabolic subgroups: $P$ defined above; the minimal parabolic $P_0=M_0 A_0 N_0$ where $M_0 \cong \{\pm 1\}$ and $N_0$ is the group of real, upper triangular matrices with ones on the diagonal; the parabolic $P'$ with Langlands decomposition $P'=M'A'N'$ where $M'\cong \SL_2^{\pm}(\R)\times\{\pm 1\}$, $A'=\left\{\diag(a,a,b,a^{-2}b^{-1})|a,b>0\right\}$ and $N'$ is the group of real, upper triangular matrices with ones on the diagonal and otherwise with zeros in the second column; and the parabolic $P''$ with Langlands decomposition $P''=M''A''N''$ where $M''\cong \SL_3^{\pm}(\R)$, $A''=\left\{\diag(a,a,a,a^{-3})|a>0\right\}$ and $N''$ is the group of real, upper triangular matrices with ones on the diagonal whose only non-zero entries above the diagonal are in the rightmost column.  

\begin{proposition}
\label{pro:parabolics}
Let $\pi$ be a principal series or complementary series representation induced from the parabolic $\bar{P}=\bar{M}\bar{A}\bar{N}$, where $\bar{P} = P_0$, $P'$ or $P''$ and let $\si$ be a finite dimensional representation on $M$.  Then $\tr\pi(\Phi_{\si})=0$ so $\pi$ contributes no zeros or poles to $Z_{P,\si}(s)$.
\end{proposition}
\prf
Let $\Th_{\pi}$ be the global character of the irreducible unitary representation $\pi$, so that
$$
\tr\pi(\Phi_{\si}) = \int_G \Phi_{\si}(x)\Th_{\pi}(x)\ dx.
$$
The Weyl integration formula can be applied (see \cite{Deitmar00}, p908) to give us
$$
\tr\pi(\Phi_{\si}) = \sum_L \rez{|W(L)|} \int_{A^+ L} \int_{M/L} f_{\si}\left( mlm^{-1}\right)dm\ \Th_{\pi}(al)d(al)\ dl,
$$
where the sum is over conjugacy classes of Cartan subgroups $L$ of $M$, and we denote by $W(L)=W(L,M)$ the Weyl group of $L$ in $M$, by $f_{\si}$ an Euler-Poincar$\eac$ function for $\si$ and $d(al)$ is an explicitly given function on $A^+ L$.  Proposition~1.4 of \cite{Deitmar00} tells us that for $L\neq B$
$$
\int_{M/L} f_{\si}\left( mlm^{-1}\right)dm=0,
$$
hence
$$
\tr\pi(\Phi_{\si}) = \rez{|W(B)|} \int_{A^+ B} \int_{M/B} f_{\si}\left( mbm^{-1}\right)dm\ \Th_{\pi}(ab)d(ab)\ db.
$$

The character $\Th_{\pi}$ of $\pi$ is non-zero only on Cartan subgroups of $G$ that are $G$-conjugate to Cartan subgroups of $\bar{M}\bar{A}$ (see \cite{Knapp86}, Proposition~10.19).  The subgroup $BA$ is not G-conjugate to any Cartan subgroup of $\bar{M}\bar{A}$, so it follows that $\tr\pi(\Phi_{\si})=0$.
\qed

Now let $\pi=\Ind_P^G(\xi\ox\nu)$ for some $\xi\in\hat{M}$ and $\nu\in\a^*$.  For $\tau \in \widehat{K_M}$ let $P_{\tau}:V_{\xi}\rightarrow V_{\xi}(\tau)$ be the projection onto the $\tau$-isotype.  For any function $f$ on $G$ which is sufficiently smooth and of sufficient decay the operator $\pi(f)$ is of trace class.  Its trace is
$$
\sum_{\tau \in \widehat{K_M}}\int_K \int_{MAN} a^{\nu + \rho_P} f(k^{-1}mank)\ \tr P_{\tau}\xi(m)P_{\tau}\ dman\ dk.
$$
Plugging in the test function $f=\Phi_{\si}$, where $\si\in\{triv,\tilde{\si}\}$, this gives us, as in \cite{Deitmar02}:
$$
\tr\pi(\Phi_{\si})=\int_{A^+} C(a)\ \tr\xi(f_{\si})\ da,
$$
where $C(a)$ depends only on $a$ and $f_{\si}$ is an Euler-Poincar\'{e} function on M attached to the representation $\si$.  We can see that $\tr\pi(\Phi_{\si})$ is non-zero only if $\tr\xi(f_{\si})$ is also.

Theorem \ref{thm:RuelleMain} now follows from Proposition \ref{pro:InfChars}.

  \chapter{A Prime Geodesic Theorem for $\SL_4(\R)$}
    \label{ch:PGT}
    \markright{\textnormal{\thechapter{. A Prime Geodesic Theorem for $\SL_4(\R)$}}}
    We continue using the notation defined in the previous chapters, in particular we take $G=\SL_4(\R)$ and take $\Ga\subset G$ to be a discrete, cocompact subgroup.  For $\ga\in\Ga$ let $N(\ga)=e^{l_{\ga}}$ \index{$N(\ga)$} and define for $x>0$
$$
\pi(x)=\sum_{{[\ga]\in\CE^{p}_P(\Ga)}\atop{N(\ga)\leq x}}\chi_1(\Ga_{\ga}),\index{$\pi(x)$}\ \textrm{ and }\ \ 
\tilde{\pi}(x)=\sum_{{[\ga]\in\CE^{p,\reg}_P(\Ga)}\atop{N(\ga)\leq x}}\chi_1(\Ga_{\ga})\tr\tilde{\si}(b_{\ga}),\index{$\tilde{\pi}(x)$}
$$
where $\Ga_{\ga}$ is the centraliser of $\ga$ in $\Ga$, we denote by $\chi_1(\Ga_{\ga})$ the first higher Euler characteristic and $\tilde{\si}$ is the virtual representation defined in Section \ref{sec:RuelleZetaFn}.  We recall from Theorem \ref{thm:ECharPos} that the first higher Euler characteristics are all positive and from Lemma \ref{lem:SigmaTilde} that the traces $\tr\tilde{\si}(b_{\ga})$ are also all positive, so both $\pi(x)$ and $\tilde{\pi}(x)$ are monotonically increasing functions.

\begin{theorem}\textnormal{(Prime Geodesic Theorem)}
\label{thm:PGT}
\index{Prime Geodesic Theorem}
For $x\rightarrow\infty$ we have
$$
\pi(x)\sim\frac{2x}{\log x}\ \textrm{ and }\ \tilde{\pi}(x)\sim\frac{8x}{\log x}.
$$
More sharply,
$$
\pi(x)=2\,\li(x)+O\left(\frac{x^{3/4}}{\log x}\right)\ \textrm{ and }\ \tilde{\pi}(x)=8\,\li(x)+O\left(\frac{x^{3/4}}{\log x}\right)
$$
as $x\rightarrow\infty$, where $\li(x)=\int_2^x \rez{\log t} dt$ \index{$\li(x)$} is the integral logarithm.
\end{theorem}

We also prove the following Prime Geodesic Theorem, which will be needed for our application to class numbers.

Let $B^0$ \index{$B^0$} be a closed subset of the compact group $B$ with the following properties: it is of measure zero; it is invariant under the map $b\mapsto b^{-1}$ and contains all fixed points of this map; and its complement $B^1=B\smallsetminus B^0$ \index{$B^1$} in $B$ is homeomorphic to an open subset of Euclidean space each of whose connected components is contractible.  The assumption that $B^0$ contains all fixed points of the map $b\mapsto b^{-1}$ is equivalent to the assumption that $B^0$ contains all non-regular elements of $B$.  Let $\CE_P^{p,1}(\Ga)$ be the subset of all $[\ga]\in\CE_P^p(\Ga)$ such that $\ga$ is conjugate in $G$ to an element of $A^-B^1$.  We define for $x>0$
$$
\pi^1(x)=\sum_{{[\ga]\in\CE^{p,1}_P(\Ga)}\atop{N(\ga)\leq x}}\chi_1(\Ga_{\ga}).\index{$\pi^1(x)$}
$$
Then as in the case of $\pi(x)$ above, $\pi^1(x)$ is a monotonically increasing function.

\begin{theorem}
\label{thm:PGT2}
For $x\rightarrow\infty$ we have
$$
\pi^1(x)\sim\frac{2x}{\log x}.
$$
\end{theorem}

The proof of these two theorems will occupy the rest of the chapter.

\section{Analytic properties of $R_{\Ga,\si}(s)$}

This section and the following proceed according to the methods of \cite{Randol77} and \cite{Randol78}, in which the Selberg zeta function for quotients of the hyperbolic plane are considered.  In this section the analogs of a series of lemmas are proved in our context.  The next section translates the main theorem of \cite{Randol77} into our context.

Recall from Theorem \ref{thm:FunctEqn} that for a finite dimensional virtual representation $\si$ of $M$ we have the functional equation
\begin{equation}
\label{eqn:FunctionalEqn}
Z_{P,\si}(1-s)=e^{-G(s)}Z_{P,\si}(s),
\end{equation}
where $G(s)$ is a polynomial.  Let $D$ be the degree of the polynomial $G(x)$.

\begin{lemma}
\label{lem:RuelleLogDerEst1}
Let $H$ be a half-plane of the form $\left\{\Re(s)<-(1+\ep)\right\}$ for some $\ep>0$ and let $\si$ be a finite dimensional virtual representation of $M$.  Then there exists a constant $C>0$ such that for $s\in H$ we have
$$
|R_{\Ga,\si}'(s)/R_{\Ga,\si}(s)|\leq C|s|^{D-1}.
$$
\end{lemma}

\prf
From Theorem \ref{thm:RuelleZeta} we get the identity
$$
R_{\Ga,\si}(s) = \prod_{q=0}^4 Z_{P,(\bwedge^q\bar{\n}\ox V_{\si})}\left( s+\frac{q}{4}\right)^{(-1)^q},
$$
which implies
\begin{equation}
\label{eqn:RuelleSelbergLogDer}
\frac{R'_{\Ga,\si}(s)}{R_{\Ga,\si}(s)} = \sum_{q=0}^4 (-1)^q \frac{Z'_{P,(\bwedge^q\bar{\n}\ox V_{\si})}\left( s+\frac{q}{4}\right)}{Z_{P,(\bwedge^q\bar{\n}\ox V_{\si})}\left( s+\frac{q}{4}\right)}.
\end{equation}
Considering this identity, it will suffice to prove that when $K$ is a half-plane of the form $\left\{\Re(s)<-\ep\right\}$ for some $\ep>0$, there exists a constant $C>0$ such that for $s\in K$ we have
$$
|Z'_{P,(\bwedge^q\n\ox V_{\si})}(s)/Z_{P,(\bwedge^q\bar{\n}\ox V_{\si})}(s)|\leq C|s|^{D-1}
$$
for all $q=0,\ldots,4$.  Since the proof does not depend on the value of $q$ or on $\si$ we shall abbreviate our notation for the zeta function to $Z_P(s)$, which notation we shall use for the rest of the section.

It follows that
$$
-\frac{Z'_P(1-s)}{Z_P(1-s)}=\frac{Z'_P(s)}{Z_P(s)} - G'(s).
$$
From the definition (\ref{eqn:GenSelb}) of $Z_P(s)$ and Proposition \ref{pro:SelbergPolesZeros}, we can see that $Z_P(s)$ is both bounded above and bounded away from zero on the half plane $K'=\left\{\Re(s)>1+\ep\right\}$.  It follows that $Z'_P(s)/Z_P(s)$ is bounded on $K'$.  This implies the lemma.
\qed

For $t>0$, let $N(t)$\index{$N(t)$}, denote the number of poles and zeros of $Z_P(s)$ at points $s=\rez{2}+x$, where $0<x<t$.

\begin{lemma}
\label{lem:NEstimate}
$N(t)=O(t^D)$
\end{lemma}

\prf
Define $\xi(s)=\left(Z_P(s)\right)^2 e^{-G(s)}$, where $G(s)$ is the polynomial in the functional equation (\ref{eqn:FunctionalEqn}) for $Z_P(s)$.  We note that, in light of its role in the functional equation, the polynomial $G(s)$ must satisfy $G(s)=-G(1-s)$.  It then follows that
$$
\xi(1-s)=\xi(s).
$$
Note that $\xi(s)$ is real on the real axis and so $\xi(\bar{s})=\overline{\xi(s)}$.

Fix a real number $1<a<5/4$ and let $t>0$.  Let $R$ be the rectangle defined by the inequalities $1-a\leq\Re s\leq a$ and $-t\leq\Im a\leq t$.  We assume that $t$ has been chosen so that no zero or pole occurs on the boundary of $R$.  Then
$$
N(t)=\rez{4}\cdot\rez{2\pi i}\int_{\partial R}\frac{\xi'(s)}{\xi(s)}ds - N_0=\rez{4}\cdot\rez{2\pi}\Im\left(\int_{\partial R}\frac{\xi'(s)}{\xi(s)}ds\right) - N_0,
$$
where $N_0$ is the number of poles and zeros of $Z_P(s)$ on the real line.  It then follows from the functional equation for $\xi(s)$ and from the fact that $\xi(\bar{s})=\overline{\xi(s)}$, that we have
$$
N(t)=\rez{2\pi}\Im\left(\int_{C}\frac{\xi'(s)}{\xi(s)}ds\right) - N_0,
$$
where $C$ is the portion of $\partial R$ consisting of the vertical segment from $a$ to $a+it$ plus the horizontal segment from $a+it$ to $\rez{2}+it$.

Now the definition of $\xi(s)$ gives us
$$
\frac{\xi'(s)}{\xi(s)}=2\frac{Z'_P(s)}{Z_P(s)} - G'(s),
$$
so that
\begin{eqnarray*}
\Im\left(\int_{C}\frac{\xi'(s)}{\xi(s)}ds\right) & = & 2.\Im\left(\int_{C}\frac{Z'_P(s)}{Z_P(s)}ds\right) - \Im\left(\int_{C}G'(s)ds\right) \\
  \\
  & = & 2.\Im\left(\int_{C}\frac{Z'_P(s)}{Z_P(s)}ds\right) - \Im\left(G\left(\rez{2}+it\right)+G(a)\right) \\
  \\
  & = & 2.\Im\left(\int_{C}\frac{Z'_P(s)}{Z_P(s)}ds\right) + O(t^D).
\end{eqnarray*}

It thus remains to show that
$$
S(t)=\Im\left(\int_{C}\frac{Z'_P(s)}{Z_P(s)}ds\right)=O(t^D).
$$
Note that $S(t)$ is the variation of the argument of $Z_P(s)$ along $C$.  We may extend the definition of $S(t)$ to those values of $t$ for which $\rez{2}+it$ is a pole or zero of $Z_P(s)$ by defining it to be $\lim_{\ep\ra 0}\rez{2}(S(t+\ep)+S(t-\ep))$.

From the definition of $Z_P(s)$ we can see that $S(t)=h(t)+O(1)$, where $h(t)$ is the variation of the argument of $Z_P(s)$ along the segment from $a+it$ to $\rez{2}+it$.  The value of $h(t)$ is bounded by a multiple of the number of zeros of $\Re (Z_P(s))$ on this segment, since the point $Z_P(s)$ cannot move between the right and left half-planes without crossing the imaginary axis.  On the segment, the real part of $Z_P(s)$ coincides with
$$
f_P(w)=\rez{2}(Z_P(w+it)+Z_P(w-it)),
$$
where $w$ runs from $\rez{2}$ to $a$ on the real axis.  Since we have assumed that $\rez{2}+it$ is not a zero or pole of $Z_P(s)$, the function $f_P(w)$ is holomorphic in a neighbourhood of the closed disc $S$, centred at $a$, of radius $a-\rez{2}$.  As this disc contains the interval from $\rez{2}$ to $a$, we may apply Jensen's Formula (\cite{Conway78}, XI.1.2) to conclude that
\begin{eqnarray*}
h(t) & = & O\left(\int_{\partial S} \log|f_P(w)|dw\right) \\
  & = & O\left(\int_{\partial S}\log |Z_P(w+it)+Z_P(w-it)|dw\right) \\
  & = & O\left(\int_{\partial S}\log |Z_P(w+it)|+\log|Z_P(w-it)|dw\right) \\
  & = & O\left(\sum_{i=1,2}\int_{\partial S}\log |Z_i(w+it)|dw+\int_{\partial S}\log |Z_i(w-it)|dw\right).
\end{eqnarray*}

There remains one final step in the proof of the lemma:
\begin{equation}
\label{eqn:Estimate}
|Z_i(x+it)|=e^{O(|t|^D)}
\end{equation}
uniformly in the strip $1-a\leq x\leq a$ for $i=1,2$.  We know that the function $Z_P(s)$ is bounded and bounded away from zero on the line $\Re s=a$.  Hence, by Lemma \ref{lem:FiniteOrder} and the functional equation (\ref{eqn:FunctionalEqn}), we can apply the Phragm\'en Lindel\"of Theorem \index{Phragm\'en Lindel\"of Theorem} (\cite{Conway78}, Theorem VI.4.1) to give (\ref{eqn:Estimate}).
\qed

\begin{lemma}
\label{lem:RuelleLogDerEst2}
Given $a<b\in\R$, there exists a sequence $(y_n)$ tending to infinity such that
$$
\left|\frac{R_{\Ga}'(x+iy_n)}{R_{\Ga}(x+iy_n)}\right|=O(y_n^{2D})
$$
for $a<x<b$.
\end{lemma}

\prf
As in Lemma \ref{lem:RuelleLogDerEst1} it will suffice to prove the result for $Z_P(s)$.

Using the notation of Lemma \ref{lem:FiniteOrder} we have that
$$
\frac{Z'_P(s)}{Z_P(s)}=\rez{s}(n_1-n_2)+g'_1(s)-g'_2(s)+\sum_{i=1,2}(-1)^{i-1}\sum_{\rho\in R_i\smallsetminus\{0\}}s^k\rho^{-k}(s-\rho)^{-1}.
$$
Let $t_0>2$ be fixed and consider the segment of the line $\Re s=\rez{2}$ given by $\rez{2}+it$ for $t_0-1<t\leq t_0+1$.  Let $N(t)$ be as above, then by Lemma \ref{lem:NEstimate} we know that $N(t)=O(t^D)$.  It follows immediately that the number of roots on the segment is $O(t_0^D)$.

By the Dirichlet principle, there exists a $\rez{2}+iy$ in the segment whose distance from any pole or zero is greater that $C/T^D$, for some fixed $C>0$.  We conclude that the portion of the sum $\sum_{i=1,2}(-1)^{i-1}\sum_{\rho}s^k\rho^{-k}(s-\rho)^{-1}$ corresponding to poles and zeros in the segment for $s_x=x+iy$ is $O(y^{2D})$, since $|s_x^k\rho^{-k}|=O(1)$ for these $\rho$, when $a<x<b$.

To deal with the segments $\rez{2}+it$ for $0<t\leq t_0-1$ and $t_0+1<t<\infty$, we proceed as follows.  The portions of the sum $\sum_{i=1,2}(-1)^{i-1}\sum_{\rho}s^k\rho^{-k}(s-\rho)^{-1}$ corresponding to the first and second segments respectively, can be written
$$
s_x^k\int_0^{t_p-1}\left(\rez{2}+it\right)^{-k}\left(s_x-\rez{2}-it\right)^{-1} dN(t)
$$
and
$$
s_x^k\int_{t_p+1}^{\infty}\left(\rez{2}+it\right)^{-k}\left(s_x-\rez{2}-it\right)^{-1} dN(t).
$$
Recalling that $N(t)=O(t^D)$, we easily conclude that both of these expressions are $O(T^{2D})$.
\qed

\section{Estimating $\psi(x)$ and $\tilde{\psi}(x)$}

To simplify notation, in what follows we write $\ga$ for an element of $\CE_P(\Ga)$ and $\ga_0$ for a primitive element and recall that $\ga\in\CE_P(\Ga)$ implies that $\ga$ is conjugate in $\Ga$ to an element $a_{\ga}b_{\ga}\in A^-B$.  Unless otherwise specified, sums involving $\ga$ or $\ga_0$ will be taken over conjugacy classes in $\CEE_P(\Ga)$ and $\CEE_P^p(\Ga)$ respectively.  If $\ga$ and $\ga_0$ occur in the same formula it is understood that $\ga_0$ wil be the primitive element underlying $\ga$.  We denote by $\CE_P^{\reg}(\Ga)$ \index{$\CE_P^{\reg}(\Ga)$} the subset of regular elements in $\CE_P(\Ga)$.  For $x>0$ let
$$
\psi(x)=\sum_{{[\ga]\in\CE_P(\Ga)}\atop{N(\ga)\leq x}} \chi_1(\Ga_{\ga})l_{\ga_0}
$$\index{$\psi(x)$}
and
$$
\tilde{\psi}(x)=\sum_{{[\ga]\in\CE_P^{\reg}(\Ga)}\atop{N(\ga)\leq x}}\chi_1(\Ga_{\ga_0})\tr\tilde{\si}(b_{\ga})l_{\ga_0}.\index{$\tilde{\psi}(x)$}
$$

Let $\si$ be a finite dimensional virtual representation of $M$.  We have for $\Re(s)>1$:
\begin{equation}
\label{eqn:Ruellelogder}
\frac{R_{\Ga,\si}'(s)}{R_{\Ga,\si}(s)}=\sum_{\ga}\chi_1(\Ga_{\ga})\tr\si(b_{\ga})l_{\ga_0}e^{-sl_{\ga}}.
\end{equation}

The following propositions are analogs of Theorem 2 of \cite{Randol77}, from which, in the next section, we prove the Prime Geodesic Theorem using standard techniques of analytic number theory.

\begin{proposition}
\label{pro:PsiHatEst}
$\psi(x)=2x+O\left(x^{3/4}\right)$
\end{proposition}
\prf
Let $D$ be the degree of the polynomial $G(s)$, as in the previous section, and suppose $k\geq 2D$ is an integer and $x,c>1$.  Then, by (\ref{eqn:Ruellelogder}) and \cite{HardyRiesz64}, Theorem~40,
\begin{eqnarray}
\label{eqn:Ruellelogderint}
\lefteqn{\rez{2\pi i}\int_{c-i\infty}^{c+i\infty} \frac{R_{\Ga,1}'(s)}{R_{\Ga,1}(s)} s^{-1}(s+1)^{-1}\cdots(s+k)^{-1}x^s ds} \nonumber \\
  \nonumber \\
  & = & \rez{2\pi i}\int_{c-i\infty}^{c+i\infty} \left(\sum_{\ga} \chi_1(\Ga_{\ga})l_{\ga_0}e^{-sl_{\ga}}\right) s^{-1}(s+1)^{-1}\cdots(s+k)^{-1}x^s ds \nonumber \\
  \nonumber \\
  & = & \rez{k!} \sum_{N(\ga)\leq x} \chi_1(\Ga_{\ga})l_{\ga_0}\left(1-\frac{N(\ga)}{x}\right)^k.
\end{eqnarray}

Theorem~\ref{thm:RuelleMain} tells us that all poles of $R_{\Ga}'(s)/R_{\Ga}(s)$ lie in the strip $-1\leq\Re(s)\leq 1$.  By virtue of Lemmas~\ref{lem:RuelleLogDerEst1} and \ref{lem:RuelleLogDerEst2} it is permissible to shift the line of integration in (\ref{eqn:Ruellelogderint}) into the half plane $\Re(s)<-1$, taking into account the residues of the poles of $R_{\Ga}'(s)/R_{\Ga}(s)$.  Hence, for $c'<-1$
\begin{eqnarray}
\label{eqn:RuellePolesSum}
\lefteqn{\rez{k!} \sum_{N(\ga)\leq x} \chi_1(\Ga_{\ga})l_{\ga_0}\left(1-\frac{N(\ga)}{x}\right)^k} \\
 & = & \sum_{{\al\in S_k}\atop{\Re(\al)>c'}} c_k(\al)x^{\al} + \rez{2\pi i}\int_{c'-i\infty}^{c'+i\infty} \frac{R_{\Ga,1}'(s)}{R_{\Ga,1}(s)} s^{-1}(s+1)^{-1}\cdots(s+k)^{-1}x^s ds, \nonumber
\end{eqnarray}
where $S_k$ denotes the set of poles of $(R_{\Ga,1}'(s)/R_{\Ga,1}(s))s^{-1}(s+1)^{-1}\cdots(s+k)^{-1}$ and $c_k(\al)$ denotes the residue at $\al$.

For $x>1$, Lemma~\ref{lem:RuelleLogDerEst1} implies the integral in (\ref{eqn:RuellePolesSum}) tends to zero as $c'\rightarrow -\infty$.  If we define
$$
\psi_0(x) =\psi(x),\ \ \psi_j(x)=\int_0^x \psi_{j-1}(t)dt,\ j\in\N,
$$
it is well known that
$$
\psi_j(x)=\rez{j!}\sum_{N(\ga)\leq x} \chi_1(\Ga_{\ga})l_{\ga_0} (x-N(\ga))^j
$$
and we deduce from (\ref{eqn:RuellePolesSum}) that
\begin{equation}
\label{eqn:Psik}
\psi_k(x)=\sum_{\al\in S_k} c_k(\al) x^{k+\al}.
\end{equation}

Let $d>0$.  For a function $f:\R\rightarrow\R$ define the operator $\D$ by setting
$$
\D f(x)=\sum_{i=0}^{2D}(-1)^i\left({2D}\atop{i}\right)f(x+(2D-i)d).
$$
It follows from (\ref{eqn:Psik}) and Theorem~\ref{thm:RuelleMain}, setting $k=2D$, that
\begin{equation}
 \psi_{2D}(x)= \frac{2}{(2D+1)!}x^{2D+1} + \sum_{\al\in S_{2D}^{\R}} c_{2D}(\al)x^{2D+\al} + \sum_{p=-2}^2 \sum_{\al\in S_{2D}^{p/4}} c_{2D}(\al)x^{2D+\al},
\end{equation}
where $S_{2D}^{\R}=S_{2D}\cap (\R\smallsetminus\{1\})$, the real elements of $S_{2D}$ not including $\al=1$, and $S_{2D}^q=S_{2D}\cap (q+i(\R\smallsetminus \{0\}))$, the non-real elements of $S_{2D}$ on the line $\Re(s)=q$.  The coefficient $((2D+1)!)^{-1}$ on the leading term comes from the fact that $(R_{\Ga,1}'(s)/R_{\Ga,1}(s))$ has a double pole at $s=1$ and from the factors $s^{-1}...(s+2D)^{-1}$.

In general if $f$ is at least $2D$ times differentiable,
\begin{equation}
 \label{eqn:DeltaId}
 \D f(x)=\int_x^{x+d} \int_{t_{2D}}^{t_{2D}+d}\cdots\int_{t_2}^{t_2+d} f^{(2D)}(t_1)\ dt_1...dt_{2D}.
\end{equation}
By applying the Mean Value Theorem we get
\begin{equation}
 \label{eqn:DeltaMVT}
 \D x^r = d^{2D} r(r-1)...(r-(2D-1))\tilde{x}^{r-2D},
\end{equation}
where $\tilde{x}\in [x,x+2Dd]$.  In particular we notice that $\D(x^{2D+1})=O(x)$, hence
$$
\D\left(\frac{2}{(2D+1)!}x^{2D+1}\right)=ax+b
$$
for some $a,b\in\R$.  Computations show that
\begin{equation}
\label{eqn:Deltax3}
a=2\sum_{j=0}^{2D}(-1)^j\rez{j!(2D-j)!}((2D-j)d)^{2D}=2d^{2D}.
\end{equation}

By definition $\psi_0(x) =\psi_{2D}^{(2D)}(s)$ so from (\ref{eqn:DeltaId}) we have
\begin{equation}
 \label{eqn:DeltaIdPsi2}
 \D\psi_{2D}(x) = \int_x^{x+d} \int_{t_{2D}}^{t_{2D}+d}\cdots\int_{t_2}^{t_2+d} \psi_0(t_1)\ dt_1...dt_{2D}.
\end{equation}
By Theorem \ref{thm:ECharPos}, the Euler characteristics $\chi_1(\Ga_{\ga})$ are all positive.  Hence $\psi_0(x)$ is non-decreasing and it follows from (\ref{eqn:DeltaIdPsi2}) that
\begin{equation}
\label{eqn:PsiIneq}
\psi_0(x)\leq d^{-2D}\D\psi_{2D}(x)\leq\psi_0(x+2Dd).
\end{equation}
It also follows from (\ref{eqn:DeltaMVT}) and (\ref{eqn:Deltax3}) that
$$
d^{-2D}\D\left(\frac{2}{(2D+1)!}x^{2D+1} + \sum_{\al\in S_{2D}^{\R}} c_{2D}(\al)x^{2D+\al}\right) = 2x + O(x^{3/4}).
$$
Thus it remains to show that for $p\in\{ -2,-1,0,1,2\}$
$$
d^{-2D}\D\left(\sum_{\al\in S_{2D}^{p/4}} c_{2D}(\al)x^{2D+\al}\right) = O(x^{3/4})
$$
and the proposition follows by (\ref{eqn:PsiIneq}).

Let $p\in\{ -2,-1,0,1,2\}$.  In order to estimate
$$
d^{-2D}\D\left(\sum_{\al\in S_{2D}^{p/4}} c_{2D}(\al)x^{2D+\al}\right)
$$
we need two estimates for $\D\left(c_{2D}(\al)x^{2D+\al}\right)$, where $\al\in S_{2D}^{p/4}$.  The residues at the poles of $R_{\Ga,1}'(s)/R_{\Ga,1}(s)$ are $O(1)$, so for $\al\in S_{2D}^{p/4}$
$$
d^{-2D}\D\left(c_{2D}(\al)x^{2D+\al}\right) = O\left(d^{-2D}|\al|^{-(2D+1)}x^{2D+p/4}\right).
$$
On the other hand, it follows from (\ref{eqn:DeltaMVT}) that
$$
d^{-2D}\D\left(c_{2D}(\al)x^{2D+\al}\right) = O\left(|\al|^{-1}x^{p/4}\right).
$$
Define $N_p(t)$ to be the number of poles of $R_{\Ga,1}'(s)/R_{\Ga,1}(s)$ on the interval $\frac{p}{4} + i(0,t]$.  From Lemma \ref{lem:NEstimate} we have that $N_p(t)=O\left( t^D\right)$.  Thus
\begin{eqnarray}
\label{eqn:NonRealPoles}
\lefteqn{d^{-2D}\D\left(\sum_{\al\in S_{2D}^{p/4}} c_{2D}(\al)x^{2D+\al}\right)} \nonumber \\
  & = & O\left(x^{p/4}\int_1^{K^D} t^{-1}dN(t) + d^{-2D}x^{2D+p/4}\int_{K^D}^{\infty}t^{-(2D+1)}dN(t)\right) \nonumber \\
  & = & O\left(K^{D-1}x^{p/4} + K^{-(D+2)}d^{-2D}x^{2+p/4}\right)
\end{eqnarray}
If we choose $K$ and $d$ appropriately, then we can conclude from (\ref{eqn:NonRealPoles}) that $d^{-2}\D\left(\sum_{\al\in S_2^{p/4}} c_2(\al)x^{2+\al}\right) = O(x^{3/4})$, as required.
\qed

\begin{proposition}
\label{pro:PsiTildeEst}
$\tilde{\psi}(x)=8x+O\left(x^{3/4}\right)$
\end{proposition}
\prf
The proof follows exactly as for the previous proposition, replacing $R_{\Ga,1}(s)$ with $R_{\Ga,\tilde{\si}}(s)$ throughout.  It follows in particular from the fact (Theorem \ref{thm:RuelleMain}) that $R_{\Ga,\tilde{\si}}(s)$ has a zero of order eight at the point $s=1$ and from the fact (Lemma \ref{lem:SigmaTilde}) that $\tr\tilde{\si}(b_{\ga})\geq 0$ for all $\ga\in\CE_P^p(\Ga)$ and $\tr\tilde{\si}(b_{\ga})=0$ if and only if $\ga$ is non-regular.
\qed

\section{The Wiener-Ikehara Theorem}

We shall use the following version of the Wiener-Ikehara theorem (see also \cite{Chandrasekharan68}, Chapter XI, Theorem 2, and \cite{Deitmar04}, Theorem 3.2).

Let $R_k(s)$, $k\in\N$ be a sequence of rational functions on $\C$.  for an open set $U\subset\C$ let $\N(U)$ be the set of natural numbers $k$ such that the pole divisor of $R_k$ does not intersect $U$.  We say that the series
$$
\sum_{k\in\N}R_k(s)
$$
\emph{converges weakly locally uniformly on $\C$} if for every open $U\subset\C$ the series
$$
\sum_{k\in\N(U)}R_k(s)
$$
converges locally uniformly on $U$.

\begin{theorem}\textnormal{(Wiener-Ikehara)}
\label{thm:WienerIkehara}
\index{Wiener-Ikehara Theorem}
Let $A(x)\geq 0$ be a monotonic measurable function on $\R_+$.  Suppose that the integral
$$
L(s)=\int_0^{\infty}A(x)e^{-sx}\ dx
$$
converges for $s\in\C$ with $\Re(s)>1$.  Suppose further that there are $j\in\N$, $r\in\R$ and a countable set $I$, and for each $i\in I$ there is $\th_i\in\C$ with $\Re(\th_i)<1$ and $c_i\in\Z$, such that the function
$$
L(s)-\left(\frac{\partial\ }{\partial s}\right)^{j+1}\frac{r}{s-1}-\sum_{i\in I}c_i\left(\frac{\partial\ }{\partial s}\right)^{j+1}\rez{s-\th_i}
$$
extends to a holomorphic function on the half-plane $\Re(s)\geq 1$.  Here we assume the sum converges weakly locally uniformly absolutely on $\C$.  Then
$$
\lim_{x\ra\infty}A(x)x^{-(j+1)}e^{-x}=r.
$$
\end{theorem}
\prf
Let $B(x)=A(x)x^{-(j+1)}e^{-x}$.  Since
$$
\frac{r}{s-1}=r\int_0^{\infty}e^{-(s-1)x}\ dx
$$
we get
$$
\left(\frac{\partial\ }{\partial s}\right)^{j+1}\frac{r}{s-1}=r\int_0^{\infty}x^{j+1}e^{-(s-1)x}\ dx,
$$
and similarly
$$
\left(\frac{\partial\ }{\partial s}\right)^{j+1}\frac{1}{s-\th_i}=\int_0^{\infty}x^{j+1}e^{-(s-\th_i)x}\ dx.
$$
Let $f$ be a smooth function of compact support on $\R$ which is real-valued and even.  Then its Fourier transform $\hat{f}$ will also be real valued.  We further assume $f$ to be of the form $f=f_1*f_1$ for some $f_1$.  Then $\hat{f}=(\hat{f}_1)^2$ is positive on the reals.  Let $I(f)$ be the set of $i\in I$ such that $\Im(\th_i)\in\supp f$.  Then the function
$$
g(s)=L(s)-\left(\frac{\partial\ }{\partial s}\right)^{j+1}\frac{r}{s-1}-\sum_{i\in I(f)}c_i\left(\frac{\partial\ }{\partial s}\right)^{j+1}\rez{s-\th_i}.
$$
extends to a holomorphic function on the set
$$
\{s\in\C:\Re(s)\geq 1, \Im(s)\in\supp f\}.
$$
We have that
$$
g(s)=\int_0^{\infty}(B(x)-r)x^{j+1}e^{-(s-1)x}\ dx-\sum_{i\in I(f)}c_i\int_0^{\infty}x^{j+1}e^{-(s-\th_i)x}\ dx.
$$

For $\ep>0$ let $g_{\ep}(t)=g(1+\ep+it)$.  We then have for $y\in\R$
\begin{eqnarray*}
\int_{\R}f(t)e^{iyt}g_{\ep}(t)dt & = & \int_{\R}f(t)e^{iyt}\left(\int_0^{\infty}(B(x)-r)x^{j+1}e^{-(\ep+it)x}dx\right)dt \\
 & & -\sum_{i\in I(f)}c_i\int_{\R}f(t)e^{iyt}\int_0^{\infty}x^{j+1}e^{-(1+\ep+it-\th_i)x}\ dx\,dt.
\end{eqnarray*}
We want to change the order of integration.  The only potential problem is with the summand involving $B(x)$.  Since $A(x)$ is non-negative and non-decreasing, we have for real $s$ and $x>0$
$$
L(s)=\int_0^{\infty}A(u)e^{-us}du\geq A(x)\int_x^{\infty}e^{-us}du=\frac{A(x)e^{-xs}}{s},
$$
that is $A(x)\leq sL(s)e^{xs}$.  Since $L(s)$ is holomorphic for $\Re(s)>1$, it follows that $A(x)=O(e^{xs})$ for every $s>1$, which implies that $A(x)=o(e^{xs})$ for every $s>1$.  Hence $B(x)e^{-\d x}=A(x)e^{-(1+\d)x}=o(1)$ for every $\d>0$.  This implies that the integral
$$
\int_0^{\infty}(B(x)-r)x^{j+1}e^{-(\ep+it)x}dx
$$
converges uniformly in the interval $-2\la\leq t\leq 2\la$.  Hence we can interchange the order of integration to obtain
\begin{eqnarray*}
\int_{\R}f(t)e^{iyt}g_{\ep}(t)dt & = & \int_0^{\infty}(B(x)-r)x^{j+1}e^{-\ep x}\left(\int_{\R}e^{i(y-x)t}f(t)dt\right)dx \\
  & & -\sum_{i\in I(f)}c_i \int_0^{\infty}x^{j+1}e^{-(1+\ep-\th_i)x}\int_{\R}e^{i(y-x)t}f(t)dt\,dx. \\
\\
  & = & \int_0^{\infty}(B(x)-r)x^{j+1}e^{-\ep x}\hat{f}(y-x)\ dx \\
  & & -\sum_{i\in I(f)}c_i \int_0^{\infty}x^{j+1}e^{-(1+\ep-\th_i)x}\hat{f}(y-x)\ dx
\end{eqnarray*}
We want to take the limit of both sides as $\ep\ra 0$ and show that the limit passes under the integration signs.

Since $g(s)$ is holomorphic for $\Re(s)\geq 0$ it follows that $g_{\ep}(t)\ra g(1+it)$ uniformly on the support of $f$, as $\ep\ra 0$.  Thus the limit of the left hand side exists and
$$
\lim_{\ep\ra 0}\int_{\R}f(t)e^{iyt}g_{\ep}(t)\ dt=\int_{\R}f(t)e^{iyt}g(1+it)\ dt.
$$
We note that since $\hat{f}$ is rapidly decreasing
$$
\lim_{\ep\ra 0}\int_0^{\infty}x^{j+1}e^{-\ep x}\hat{f}(y-x)\ dx=\int_0^{\infty}x^{j+1}\hat{f}(y-x)\ dx,
$$
and since $B(x)$ is non-negative and non-decreasing it follows that
$$
\lim_{\ep\ra 0}\int_0^{\infty}(B(x)-r)x^{j+1}e^{-\ep x}\hat{f}(y-x)\ dx=\int_0^{\infty}(B(x)-r)x^{j+1}\hat{f}(y-x)\ dx.
$$
Considering the sum over $I(f)$, we note that the imaginary parts of the $\th_i$ for $i\in I_{\la}$ lie in a compact set and hence the convergence of the sum is uniform in $\ep$ and the limit may be drawn under the summation.  Finally, we note once more that the integrand is non-negative and monotonically increasing as $\ep\ra 0$, so the limit may once more be taken under the integral.  We conclude that
\begin{eqnarray*}
\int_{\R}f(t)e^{iyt}g(1+it)\ dt & = & \int_0^{\infty}(B(x)-r)x^{j+1}\hat{f}(y-x)\ dx \\
  & & -\sum_{i\in I(f)}c_i \int_0^{\infty}x^{j+1}e^{-(1-\th_i)x}\hat{f}(y-x)\ dx.
\end{eqnarray*}
By the Riemann-Lebesgue lemma the left hand side tends to zero as $y\ra\infty$.  For $y\geq 0$ and $\th\in\C$ with $\Re(\th)<1$ we estimate
\begin{eqnarray*}
\int_0^{\infty}x^{j+1}e^{(\th-1)x}\hat{f}(y-x)\ dx & \leq & \int_{-\infty}^{\infty}x^{j+1}e^{(\th-1)x}\hat{f}(y-x)\ dx \\
  & = & \int_{-\infty}^{\infty}(x+y)^{j+1}e^{(\th-1)(x+y)}\hat{f}(-x)\ dx \\
  & \leq & Cy^{j+1}e^{(\th-1)y},
\end{eqnarray*}
for some constant $C$.  It follows that the sum over $I(f)$ tends to zero as $y\ra\infty$.  We therefore have
\begin{equation}
\label{eqn:WienerIkehara}
\lim_{y\ra\infty}\int_0^{\infty}B(x)x^{j+1}\hat{f}(y-x)\ dx=r\lim_{y\ra\infty}\int_0^{\infty}x^{j+1}\hat{f}(y-x)\ dx.
\end{equation}

\begin{lemma}
For every non-negative integer $k$
$$
\lim_{y\ra\infty}\rez{y^k}\int_0^{\infty}x^k\hat{f}(y-x)dx=2\pi f(0).
$$
\end{lemma}
\prf
We prove the lemma using induction.  First, let $k=0$.  Then
\begin{eqnarray*}
\lim_{y\ra\infty}\int_0^{\infty}\hat{f}(y-x)dx & = & \lim_{y\ra\infty}\int_{-\infty}^{\infty}\hat{f}(y-x)dx \\
  & = & \int_{-\infty}^{\infty}\hat{f}(-x)dx \\
  & = & 2\pi f(0).
\end{eqnarray*}
We now assume the lemma has been proved for a fixed value of $k$.
\begin{eqnarray*}
\lim_{y\ra\infty}\rez{y^{k+1}}\int_0^{\infty}x^{k+1}\hat{f}(y-x)dx & = & \lim_{y\ra\infty}\rez{y^{k+1}}\int_{-\infty}^{\infty}x^{k+1}\hat{f}(y-x)dx \\
  & = & \lim_{y\ra\infty}\rez{y^{k+1}}\int_{-\infty}^{\infty}(x+y)^{k+1}\hat{f}(-x)dx \\
  & = & \lim_{y\ra\infty}\rez{y^{k+1}}\int_{-\infty}^{\infty}x^{k+1}\hat{f}(-x)dx \\
  &   & \ \ \ + \lim_{y\ra\infty}\rez{y^k}\int_{-\infty}^{\infty}(x+y)^k\hat{f}(-x)dx \\
  & = & \lim_{y\ra\infty}\rez{y^k}\int_0^{\infty}x^k\hat{f}(y-x)dx \\
  & = & 2\pi f(0).
\end{eqnarray*}
\qed

The lemma together with (\ref{eqn:WienerIkehara}) implies that
$$
\lim_{y\ra\infty}\int_0^{\infty}B(x)\left(\frac{x}{y}\right)^{j+1}\hat{f}(y-x)dx=r2\pi f(0).
$$
Let $S>0$.  Since $A(x)$ is monotonic we have $A(y-S)\leq A(x)\leq A(y+S)$ whenever $y-S\leq x\leq y+S$.  For $x$ in that range we then have by the definition of $B(x)$
$$
B(y-s)(y-S)^{j+1}e^{y-S}\ \ \leq\ \ B(x)x^{j+1}e^x\ \ \leq\ \ B(y+S)(y+S)^{j+1}e^{y+S}.
$$
The first inequality implies
\begin{eqnarray*}
B(x)x^{j+1} & \geq & B(y-S)(y-S)^{j+1}e^{y-x-S} \\
  & \geq & B(y-S)(y-S)^{j+1}e^{-2S}.
\end{eqnarray*}
So for $y\geq S$,
\begin{eqnarray*}
e^{-2S}B(y-S)\frac{(y-S)^{j+1}}{y^{j+1}}\int_{y-S}^{y+S}\hat{f}(y-x)dx & \leq & \int_{y-S}^{y+S}B(x)\frac{x^{j+1}}{y^{j+1}}\hat{f}(y-x)dx \\
  & \leq & \int_0^{\infty}B(x)\frac{x^{j+1}}{y^{j+1}}\hat{f}(y-x)dx,
\end{eqnarray*}
which implies
$$
\limsup_{y\ra\infty}B(y)\leq re^{2S}\frac{2\pi f(0)}{\int_{-S}^{S}\hat{f}(x)dx}.
$$
We vary $f$ so that $\hat{f}$ is small outside $[-S,S]$.  In this way we get
$$
\limsup_{y\ra\infty}B(y)\leq re^{2S}.
$$
Since this holds for any $S>0$ it follows that
$$
\limsup_{y\ra\infty}B(y)\leq r.
$$
The inequality
$$
\liminf_{y\ra\infty}B(y)\geq r.
$$
is obtained in a similar fashion.  The theorem is proven.
\qed

\section{The Dirichlet series}
\label{sec:DirSeries}
Let $B^0$ \index{$B^0$} be a closed subset of the compact group $B$ with the following properties: it is of measure zero; it is invariant under the map $b\mapsto b^{-1}$ and contains all fixed points of this map; and its complement $B^1=B\smallsetminus B^0$ \index{$B^1$} in $B$ is homeomorphic to an open subset of Euclidean space each of whose connected components is contractible.

The Weyl group $W=W(M,B)$ contains two elements and the non-trivial element acts on $B$ by $b\mapsto b^{-1}$, so the invariance condition above says that both $B^0$ and $B^1$ are invariant under the action of $W$.  The fixed points under the action of the Weyl group are precisely the non-regular elements of $B$, so the assumption that $B^0$ contains all these fixed points is equivalent to $B^{\nreg}\subset B^0$, where $B^{\nreg}$ \index{$B^{\nreg}$} denotes the subset of nonregular elements of $B$.  The action of $W$ on $B^1$ permutes the connected components and the assumption $B^{\nreg}\subset B^0$ implies that the quotient space $B^1/W$ is also homeomorphic to an open subset of Euclidean space each of whose connected components is simply connected.

For subsets $S$ and $T$ of $G$ we denote by $S.T$ the subset
$$
S.T=\{sts^{-1}:s\in S,t\in T\}
$$
of $G$.  Let $M_{\ell}^1=M.B^1\subset M$.  Let $\CE_P^0(\Ga)$ \index{$\CE_P^0(\Ga)$} and $\CE_P^1(\Ga)$ \index{$\CE_P^1(\Ga)$} be the subsets of $\CE_P(\Ga)$ consisting of all conjugacy classes $[\ga]$ such that $b_{\ga}\in B^0$ or $b_{\ga}\in B^1$ respectively.  The assumption that $B^0$ is of measure zero is not required for the following lemma, but will be needed later on.

\begin{lemma}
There exist a set $M_c\subset M_{\ell}^1$ with compact closure in $M$ and a monotonically increasing sequence $(g_n)$ \index{$g_n$} of smooth functions on $M$, supported on $M_c$, which are invariant under conjugation by elements of $K_M$, and whose orbital integrals satisfy
$$
\CO_{b_{\ga}}^M(g_n)=\int_{M/B}g_n(xb_{\ga}x^{-1})\ dx\ra 1\ \ \ \textrm{ as }n\ra\infty
$$
for all $\ga\in\CE_P^1(\Ga)$.
\end{lemma}
\prf
We view $M_{\ell}^1$ as a fibre bundle with base space $B^1/W$ and fibres homeomorphic to $M/B$.

Let $d(\cdot,\cdot)$ denote the metric on $B$ given by the form $b$ in (\ref{eqn:norm}).  For $n\in\N$ let $\tilde{B}_n\subset B^1$ be the set $\tilde{B}_n=\{b\in B:d(b,B^0)\geq 1/n\}$ and let $B_n=\tilde{B}_n/W$.  Then, by \cite{Warner83}, Corollary 1.11, for each $n\in\N$ there exists a function $h_n:B/W\ra\R$ such that $0\leq h_n(b)\leq 1$ for all $b\in B/W$, for all $b\in B_n$ we have $h_n(b)=1$ and $h_n$ is supported on $B^1/W$.  We may assume that the $h_n$'s form an increasing series.

Let $U$ be a compact neighbourhood of a point in $M/B$ homeomorphic to a subset of Euclidean space.  Let $k:M/B\ra\R$ be a smooth positive function supported on $U$ and satisfying $\int_{M/B}k(m)dm=1$.

We now define, for $n\in\N$, functions $\tilde{g}_n:M\ra\R$ as follows.  On each connected component $V$ of $B^1/W$ we fix a trivialisation of the restriction to $V$ of the bundle $M/B\ra M_{\ell}^1\ra B^1/W$.  Then for $v\in V$ and $m\in M/B$ define $\tilde{g}_n(mvm^{-1})=h_n(v)k(m)$.  Since $V$ is simply connected there are no problems with global agreement of this definition.  For $m\in M\smallsetminus M_{\ell}^1$ define $\tilde{g}_n(m)=0$.  Finally we define the functions $g_n:M\ra\R$ by
$$
g_n(m)=\rez{2}\int_{K_M}\tilde{g}_n(kmk^{-1})\,dk
$$
for all $m\in M$.  Then the functions $g_n$ form an increasing sequence of smooth $K_M$-invariant functions supported on the compact set $M_c$, which we define to be the closure in $M$ of $K_M.(U.B^1)$.  We will show that their orbital integrals have the required properties.  Let $b\in B^1$.
\begin{eqnarray*}
\CO_{b}^M(g_n) & = & \int_{M/B}g_n(mbm^{-1})\ dm \\
  & = & \rez{2}\int_{M/B}\int_{K_M}\tilde{g}_n(kmbm^{-1}k^{-1})\ dk\,dm \\
  & = & \rez{2}\int_{M/B}\int_{K_M}\tilde{g}_n(mkbk^{-1}m^{-1})\ dk\,dm \\
  & = & \rez{2}\int_{M/B}\int_{B}\tilde{g}_n(mbm^{-1}) + \tilde{g}_n(mb^{-1}m^{-1})\ dk\,dm \\
  & = & \rez{2}\left(h_n(b)+h_n(b^{-1})\right)\int_{M/B}k(m)\ dm \\
  & = & h_n(b).
\end{eqnarray*}
The last equality holds since $bW=b^{-1}W$ in $B^1/W$.  There exists $N\in\N$ such that $b\in B_n$ for all $n\geq N$, hence $\CO_{b}^M(g_n)\ra 1$ as $n\ra\infty$, by the definition of the functions $h_n$.
\qed

For $n\in\N$ let
$$
L^j_n(s)=\sum_{[\ga]\in\CE_P^1(\Ga)}\chi_1(\Ga_{\ga_0})\CO_{b_{\ga}}(g_n)l_{\ga_0}\frac{l_{\ga}^{j+1}e^{-sl_{\ga}}}{\det(1-(a_{\ga}b_{\ga})^{-1}|\n)},\index{$L^j_n(s)$}
$$
let
$$
r_n=\int_M g_n(x)dx,\index{$r_n$}
$$
and let
$$
m_{n,\la}=\sum_{q=0}^4 (-1)^q r_n\dim H^q(\n,\pi_K)^{\la}.\index{$m_{n,\la}$}
$$

\begin{proposition}
\label{pro:DirSeries}
For all $n\in\N$ and for $j\in\N$ large enough the series $L^j_n(s)$ converges locally uniformly in the set $\{s\in\C:\Re(s)>1\}$.

The function $L^j_n(s)$ can be written as a Mittag-Leffler series
$$
L_n^j(s) = \left(\frac{\partial\ }{\partial s}\right)^{j+1}\frac{r_n}{s-1} + \sum_{\pi\in\hat{G}}N_{\Ga}(\pi)\sum_{\la\in\a^*}m_{n,\la}\left(\frac{\partial\ }{\partial s}\right)^{j+1}\rez{s-\la(H_1)},
$$
where the summand of the double series corresponding to $\pi=triv$, $\la=-2\rho_P$ is excluded.  The double series converges weakly locally uniformly on $\C$.  For $\pi\in\hat{G}$, $\la\in\a^*$ such that $(\pi,\la)\neq (triv,-2\rho_P)$ we have $m_{n,\la}\neq 0$ only if $\Re(\la)>-2\rho_P$.  Thus, in particular, the double series converges locally uniformly on $\{s\in\C:\Re(s)>1\}$.
\end{proposition}
\prf
We can put $f=g_n$ in (\ref{eqn:TestFunction}) to define the function $\Phi_{g_n,j,s}$, which by Proposition \ref{pro:STF2}, for $j$ and $\Re(s)$ large enough, goes into the Selberg trace formula to give the equation
$$
\sum_{\pi\in\hat{G}} N_{\Ga}(\pi)\tr\pi(\Phi_{g_n})=\sum_{[\ga]\in\CE_P(\Ga)}\vol(\Ga_{\ga}\bs G_{\ga})\CO_{b_{\ga}}(g_n)\frac{l_{\ga}^{j+1}e^{-sl_{\ga}}}{\det(1-(a_{\ga}b_{\ga})^{-1}|\n)}.
$$
By the definition of the functions $g_n$, the orbital integrals $\CO_{b_{\ga}}(g_n)$ are equal to zero when $\ga\in\CE_P^0(\Ga)$.  Also, for $\ga\in\CE_P^1(\Ga)$, since we have assumed that $B^{\nreg}\subset B^0$ we have that $\ga$ is regular and we get from Lemma \ref{lem:ECharEqn} that
$$
\vol(\Ga_{\ga}\bs G_{\ga})=\chi_1(\Ga_{\ga})l_{\ga_0}.
$$
Thus we can see that the geometric side of the trace formula is equal to $L_n^j(s)$.  It follows that for $j$ and $\Re(s)$ large enough the Dirichlet series $L_n^j(s)$ converges absolutely.

\begin{lemma}
$$
\tr\pi(\Phi_{g_n})=(-1)^{j+1}\left(\frac{\partial\ }{\partial s}\right)^{j+1}r_n\sum_{\la\in\a^*}\dim H^{\bullet}(\n,\pi_K)^{\la}\rez{s-\la(H_1)}.
$$
\end{lemma}
\prf
Replacing $f_{\si}$ with $g_n$ in (\ref{eqn:TracePiPhi1}) we get
\begin{eqnarray*}
\tr\pi(\Phi_{g_n}) & = & \int_{MA^-}g_n(m)\Th^{MA}_{H^{\bullet}(\n,\pi_K)}(ma)dm\,g_s^j(a)da \\
  & = & \int_{A^-}\int_M g_n(m)dm\ \tr\!\left(a|H^{\bullet}(\n,\pi_K)\right)l_a^{j+1}e^{-sl_a}da. \\
  & = & r_n\int_0^{\infty}\sum_{\la\in\a^*}\dim H^{\bullet}(\n,\pi_K)^{\la}e^{(\la(H_1)-s)t}\,t^{j+1}\,dt \\
  & = & (-1)^{j+1}\left(\frac{\partial\ }{\partial s}\right)^{j+1}r_n\sum_{\la\in\a^*}\dim H^{\bullet}(\n,\pi_K)^{\la}\rez{s-\la(H_1)}.
\end{eqnarray*}
\qed

By the above lemma and Proposition \ref{pro:AAction}, and using similar arguments to those used in the proof of Proposition \ref{pro:TrivAAction}, we can see that the spectral side of the trace formula is equal to the Mittag-Leffler series given in the propostion.

Since $L_n^j(s)$ is a Dirichlet series with positive coefficients it will converge locally uniformly for $s$ in some open set.  By proving the convergence of the Mittag-Leffler series we shall show that the Dirichlet series extends to a holomorphic function on $\{\Re(s)>1\}$ and hence, since it has positive coefficients, converges locally uniformly there.

We recall that for $\la\in\a^*$ we write $\|\la\|$ for the norm given by the form $b$ in (\ref{eqn:norm}).  It then follows from Proposition \ref{pro:Bounds} that there exist $m_1\in\N$ and $C>0$ such that for every $\pi\in\hat{G}$ and every $\la\in\a^*$ we have
$$
|m_{n,\la}|\leq C(1+\|\la\|)^{m_1}.
$$
If $S$ denotes the set of all pairs $(\pi,\la)\in\hat{G}\x\a^*$ such that $m_{n,\la}\neq 0$, then there exists $m_2\in\N$ such that
$$
\sum_{(\pi,\la)\in S}\frac{N_{\Ga}(\pi)}{(1+\|\la\|)^{m_2}}<\infty.
$$

Now, let $U\subset\C$ be open and let $S(U)$ be the set of all pairs $(\pi,\la)\in\hat{G}\x\a^*$ such that $m_{n,\la}\neq 0$ and $\la(H_1)\notin U$. Let $V$ be a compact subset of $U$.  We have to show that for some $j\in\N$, which does not depend on $U$ or $V$,
$$
\sup_{s\in V}\sum_{(\pi,\la)\in S(U)}\left|\frac{N_{\Ga}(\pi)m_{n,\la}}{(s-\la(H_1))^{j+2}}\right|<\infty.
$$
Let $m_1$, $m_2$ be as above and let $j\geq m_1+m_2-2$.  Since $V\subset U$ and $V$ is compact there is $\ep>0$ such that $s\in V$ and $(\pi,\la)\in S(U)$ implies $|s-\la(H_1)|\geq\ep$.  Hence there is $c>0$ such that for every $s\in V$ and every $(\pi,\la)\in S(U)$,
$$
|(s-\la(H_1)|\geq c(1+\|\la\|).
$$
Putting this all together we get
\begin{eqnarray*}
\sup_{s\in V}\sum_{(\pi,\la)\in S(U)}\left|\frac{N_{\Ga}(\pi)m_{n,\la}}{(s-\la(H_1))^{j+2}}\right| & \leq & \sup_{s\in V}\sum_{(\pi,\la)\in S(U)} \rez{c^{j+2}}\frac{N_{\Ga}(\pi)|m_{n,\la}|}{(1+\|\la\|)^{j+2}} \\
  & \leq & \sup_{s\in V}\sum_{(\pi,\la)\in S(U)}\frac{C}{c^{j+2}}\frac{N_{\Ga}(\pi)}{(1+\|\la\|)^{m_2}} \\
  & < & \infty,
\end{eqnarray*}
which proves the proposition.
\qed

\section{Estimating $\psi_n(x)$}

Let
$$
\phi^j_n(x)=\sum_{{[\ga]\in\CE_P^1(\Ga)}\atop{N(\ga)\leq x}} \chi_1(\Ga_{\ga_0})\CO_{b_{\ga}}(g_n)l_{\ga_0}\frac{l_{\ga}^{j+1}}{\det(1-(a_{\ga}b_{\ga})^{-1}|\n)}.\index{$\phi^j_n(x)$}
$$

\begin{lemma}
\label{lem:Abel}
$$
\int_0^{\infty}\phi^j_n(e^x)e^{-sx}\ dx=L^j_n(s)
$$
\end{lemma}
\prf
This follows from Abel's summation formula (\cite{Chandrasekharan68}, Chapter VII, Theorem 6).
\qed

Let
$$
\phi_n(x)=\sum_{{[\ga]\in\CE_P^1(\Ga)}\atop{N(\ga)\leq x}} \chi_1(\Ga_{\ga_0})\CO_{b_{\ga}}(g_n)l_{\ga_0}\rez{\det(1-(a_{\ga}b_{\ga})^{-1}|\n)}.\index{$\phi_n(x)$}
$$

\begin{lemma}
\label{lem:Phin}
For each $n\in\N$ we have
$$
\phi_n(x)\sim r_n x.
$$
\end{lemma}
\prf
By Propostion \ref{pro:DirSeries} and Lemma \ref{lem:Abel} we can apply Theorem \ref{thm:WienerIkehara} to the series $L_n^j(s)$ and the function $\phi_n^j(x)$ to deduce that
\begin{equation}
\label{eqn:PhijCgence}
\lim_{x\ra\infty}\frac{\phi^j_n(x)}{x(\log x)^{j+1}}=r_n.
\end{equation}
Also it is clear that
$$
\phi_n(x)\geq\frac{\phi^j_n(x)}{(\log x)^{j+1}},
$$
so it follows that
$$
\liminf_{x\ra\infty}\frac{\phi_n(x)}{x}\geq r_n.
$$
Let $0<\mu<1$.  Then
\begin{eqnarray*}
\phi^j_n(x) & \geq & \sum_{{[\ga]\in\CE_P^1(\Ga)}\atop{x^{\mu}<N(\ga)\leq x}} \chi_1(\Ga_{\ga_0})\CO_{b_{\ga}}(g_n)l_{\ga_0}\frac{l_{\ga}^{j+1}}{\det(1-(a_{\ga}b_{\ga})^{-1}|\n)} \\
  & \geq & \mu^{j+1}(\log x)^{j+1}\sum_{{[\ga]\in\CE_P^1(\Ga)}\atop{x^{\mu}<N(\ga)\leq x}} \chi_1(\Ga_{\ga_0})\CO_{b_{\ga}}(g_n)l_{\ga_0}\rez{\det(1-(a_{\ga}b_{\ga})^{-1}|\n)} \\
  & = & \mu^{j+1}(\log x)^{j+1}(\phi_n(x)-\phi_n(x^{\mu})).
\end{eqnarray*}
From this we get
\begin{equation}
\label{eqn:PhinBound}
\frac{\phi_n(x)}{x}\leq\mu^{-(j+1)}\frac{\phi^j_n(x)}{x(\log x)^{j+1}}+\frac{\phi_n(x^{\mu})}{x^{\mu}x^{1-\mu}}.
\end{equation}

Assume first that $\phi_n(x)/x$ tends to infinity as $x\ra\infty$.  Then there is a sequence $x_m$ of positive real numbers, tending to infinity such that $\phi(x_m)/x_m$ tends to infinity as $m\ra\infty$ and
$$
\frac{\phi_n(x_m)}{x_m}\geq\frac{\phi_n(x_m^{\mu})}{x_m^{\mu}}
$$
for all $m$.  Then
$$
\frac{\phi_n(x_m)}{x_m}\leq\mu^{-(j+1)}\frac{\phi_n^j(x_m)}{x_m(\log x_m)^{j+1}}+\frac{\phi_n(x_m)}{x_m}.\rez{x_m^{1-\mu}},
$$
so that
$$
\frac{\phi_n(x_m)}{x_m}\leq\frac{\mu^{-(j+1)}\frac{\phi_n^j(x_m)}{x_m(\log x_m)^{j+1}}}{1-\rez{(x_m)^{1-\mu}}}.
$$
By (\ref{eqn:PhijCgence}) the right hand side converges, so we have a contradiction.  This implies that
$$
\lim_{x\ra\infty}\frac{\phi_n(x)}{x}
$$
exists and is finite and we set
$$
L=\limsup_{x\ra\infty}\frac{\phi_n(x)}{x}=\limsup_{x\ra\infty}\frac{\phi_n(x^{\mu})}{x^{\mu}}.
$$
From (\ref{eqn:PhijCgence}) and (\ref{eqn:PhinBound}) we get
\begin{eqnarray*}
L & \leq & r_n\mu^{-(j+1)}+L\limsup_{x\ra\infty}\rez{x^{1-\mu}} \\
  & = & r_n\mu^{-(j+1)}.
\end{eqnarray*}
Since this holds for any value of $\mu$ in the interval $0<\mu<1$ we get that $L\leq r$ and the lemma follows.
\qed

For $n\in\N$ and $x>0$ let
$$
\psi_n(x)=\sum_{{[\ga]\in\CE_P^1(\Ga)}\atop{N(\ga)\leq x}} \chi_1(\Ga_{\ga_0})\CO_{b_{\ga}}(g_n)l_{\ga_0}.\index{$\psi_n(x)$}
$$

\begin{proposition}
\label{pro:PsiNEst}
For each $n\in\N$ we have
$$
\psi_n(x)\sim r_n x.
$$
\end{proposition}
\prf
We note that $0<\det(1-a_{\ga}b_{\ga}|\n)<1$ for all $\ga\in\CE_P^{\reg}(\Ga)$ and that the value of the determinant tends to $1$ as $l_{\ga}$ tends to infinity, hence for $0<\ep<1$ there are only finitely many $\ga\in\CE_P^{\reg}(\Ga)$ such that $\det(1-a_{\ga}b_{\ga}|\n)\leq 1-\ep$.  Since $\CE_P^1(\Ga)\subset\CE_P^{\reg}(\Ga)$ the same holds for $\ga\in\CE_P^1(\Ga)$.  We fix $0<\ep<1$ and define for each $n\in\N$ the functions $\phi_{n,\ep}$ and $\psi_{n,\ep}$ to be the same sums as for $\phi_n$ and $\psi_n$ respectively but restricted to those classes $[\ga]\in\CE_P^1(\Ga)$ such that $1-\ep<\det(1-a_{\ga}b_{\ga}|\n)<1$.  It then follows that both
\begin{equation}
\label{eqn:PhiEp}
\frac{\phi_n(x)-\phi_{n,\ep}(x)}{x}\ra 0
\end{equation}
and
\begin{equation}
\label{eqn:PsiEp}
\frac{\psi_n(x)-\psi_{n,\ep}(x)}{x}\ra 0
\end{equation}
as $x\ra\infty$.  Now
$$
1-\ep<\det(1-a_{\ga}b_{\ga}|\n)<1
$$
immediately implies
$$
\frac{1-\ep}{\det(1-a_{\ga}b_{\ga}|\n)}<1<\rez{\det(1-a_{\ga}b_{\ga}|\n)},
$$
so summing up we get
$$
\frac{\phi_{n,\ep}(x)}{x}(1-\ep)<\frac{\psi_{n,\ep}(x)}{x}<\frac{\phi_{n,\ep}(x)}{x}.
$$
By Lemma \ref{lem:Phin} and (\ref{eqn:PhiEp}), it then follows that
$$
r_n(1-\ep)\leq\liminf_{x\ra\infty}\frac{\psi_{n,\ep}(x)}{x}\leq\limsup_{x\ra\infty}\frac{\psi_{n,\ep}(x)}{x}\leq r_n.
$$
It then follows from (\ref{eqn:PsiEp}) that
$$
r_n(1-\ep)\leq\liminf_{x\ra\infty}\frac{\psi_n(x)}{x}\leq\limsup_{x\ra\infty}\frac{\psi_n(x)}{x}\leq r_n.
$$
Since this holds for any value of $\ep$ in the interval $0<\ep<1$ the proposition is proven.
\qed

\section{Estimating $\psi^1(x)$}

\begin{lemma}
\label{lem:RnLimit}
The sequence $(r_n)$ is monotonically increasing and $r_n\ra 2$ as $n\ra\infty$.
\end{lemma}
\prf
That the sequence $(r_n)$ is monotonically increasing is clear since the sequence of functions $(g_n)$ is monotonically increasing.

By Weyl's integral formula (\cite{Knapp86}, Proposition 5.27), and since the functions $g_n$ are only non-zero on elliptic elements of $M$, we have
\begin{eqnarray*}
r_n & = & \int_M g_n(x)\ dx \\
  & = & \rez{|W(B:G)|}\int_B\int_{M/B}g_n(xbx^{-1})|D(b)|^2\ dx\,db \\
  & = & \rez{|W(B:G)|}\int_B\CO^M_b(g_n)|D(b)|^2\ db,
\end{eqnarray*}
where $W(B:G)$ is the Weyl group and $D(b)$ is the Weyl denominator.  Since the sequence of functions $g_n$ is monotonically increasing and the functions are all supported within a given compact subset of $M$ we can interchange integral and limit to get
\begin{eqnarray*}
\lim_{n\ra\infty}r_n & = & \lim_{n\ra\infty}\rez{|W(B:G)|}\int_B\CO^M_b(g_n)|D(b)|^2\ db \\
  & = & \rez{|W(B:G)|}\int_B\lim_{n\ra\infty}\CO^M_b(g_n)|D(b)|^2\ db.
\end{eqnarray*}
Furthermore, the orbital integrals $\CO^M_b(g_n)$ tend to one as $n\ra\infty$, except for $b$ in a set of measure zero, so it follows that
$$
\lim_{n\ra\infty}r_n=\rez{|W(B:G)|}\int_B|D(b)|^2\ db.
$$

The Weyl group $W(B:G)$ consists of the identity element and the element
$$
\matrixfour{1}{-1}{1}{-1}.
$$
We can also compute the Weyl denominator $D(b)$ for $b\in B$.  Let
$$
R(\th)=\matrix{\cos\th}{-\sin\th}{\sin\th}{\cos\th}\in\SO(2)
$$
and let
$$
R(\th,\phi)=\matrixtwo{R(\th)}{R(\phi)}\in B.
$$
Then
\begin{eqnarray*}
D(R(\th,\phi)) & = & e^{i(\th+\phi)}(1-e^{-2i\th})(1-e^{-2i\phi}) \\
  & = & (e^{i(\th+\phi)}+e^{-i(\th+\phi)})-(e^{i(\th-\phi)}+e^{-i(\th-\phi)}) \\
  & = & 2\cos(\th+\phi)-2\cos(\th-\phi) \\
  & = & 2(\cos\th\cos\phi-\sin\th\sin\phi-\cos\th\cos\phi-\sin\th\sin\phi) \\
  & = & -4\,\sin\th\,\sin\phi
\end{eqnarray*}
We have normalised the Haar measure on $B$ so that $\int_B db=1$, hence
\begin{eqnarray*}
\rez{|W(B:G)|}\int_B|D(b)|^2\ db & = & \rez{2}\int_0^{2\pi}\int_0^{2\pi}16\,\sin\!^2\th\,\sin\!^2\phi\ \frac{d\th\,d\phi}{4\pi^2} \\
  & = & \frac{2}{\pi^2}\int_0^{2\pi}\int_0^{2\pi}\sin\!^2\th\,\sin\!^2\phi\ d\th\,d\phi \\
  & = & \frac{2}{\pi^2}.\pi^2 \\
  & = & 2.
\end{eqnarray*}
This proves the lemma.
\qed

Let
$$
\psi^1(x)=\sum_{{[\ga]\in\CE_P^1(\Ga)}\atop{N(\ga)\leq x}}\chi_1(\Ga_{\ga_0})l_{\ga_0}.\index{$\psi(x)$}
$$
\begin{proposition}
\label{pro:PsiEst}
$\psi^1(x)\sim\frac{2x}{\log x}$.
\end{proposition}
\prf
For all $n\in\N$ and all $x>0$ we have
$$
\psi_n(x)\leq\psi^1(x)\leq\psi(x).
$$
It then follows, using Proposition \ref{pro:PsiHatEst} and Proposition \ref{pro:PsiNEst}, that for all $n\in\N$
$$
r_n=\liminf_{x\ra\infty}\frac{\psi_n(x)}{x}\leq\liminf_{x\ra\infty}\frac{\psi^1(x)}{x}\leq\limsup_{x\ra\infty}\frac{\psi^1(x)}{x}\leq\limsup_{x\ra\infty}\frac{\psi(x)}{x}=2.
$$
We can then deduce from Lemma \ref{lem:RnLimit} that
$$
\lim_{x\ra\infty}\frac{\psi^1(x)}{x}=2.
$$
\qed

\section{Estimating $\pi(x)$, $\tilde{\pi}(x)$ and $\pi^1(x)$}

To keep the notation less cluttered, in the sums over conjugacy classes in $\Ga$ which appear in this section we shall not specify which set of classes the sum is being taken over, since this will be clear from the context.  We shall always use $\ga_0$ to denote primitive elements and where $\ga$ and $\ga_0$ appear together in the same formula we shall mean that $\ga_0$ is the primitive element underlying $\ga$.

\begin{proposition}
\label{pro:PiEst}
$$
\displaystyle\lim_{x\ra\infty}\frac{\pi(x)}{2x/\log x} = \lim_{x\ra\infty}\frac{\psi(x)}{2x} = 1.
$$
\end{proposition}
\prf
We can write
$$
\psi(x)=\sum_{N(\ga_0)\leq x} n_{\ga_0}\chi_1(\Ga_{\ga_0})l_{\ga_0},
$$
where $n_{\ga_0}\in\N$ is maximal such that $N(\ga_0)^{n_{\ga_0}}\leq x$.  By definition $N(\ga_0)=e^{l_{\ga_0}}$, so $N(\ga_0)^n\leq x$ implies $nl_{\ga_0}\leq \log x$ and we can see that
\begin{equation}
\label{eqn:PsiPiIneq}
\psi(x)\leq\log(x)\pi(x).
\end{equation}

Next we fix a real number $0<a<1$.  By Theorem \ref{thm:ECharPos} the Euler characteristic $\chi_1(\Ga_{\ga_0})>0$ for all $\ga_0\in\CE_P^p(\Ga)$, so for $x>1$,
$$
\psi(x) \geq \sum_{x^a <N(\ga_0)\leq x} \chi_1(\Ga_{\ga_0})l_{\ga_0}.
$$
As above, $N(\ga_0)>x^a$ implies $l_{\ga_0}>\log x^a$, hence
$$
\psi(x)\geq a\log x\sum_{x^a <N(\ga_0)\leq x}\chi_1(\Ga_{\ga_0}) = a\log x(\pi(x)-\pi(x^a)).
$$
Since $\Ga\subset G$ is discrete, $\pi(x^a)<Cx^a$ for some constant $C$ so
$$
\psi(x)>a\pi(x)\log x - aCx^a\log x,
$$
which gives
$$
\frac{\psi(x)}{2x}>a\pi(x)\frac{\log x}{2x} - aC\frac{\log x}{2x^{1-a}}.
$$
Since $0<a<1$, it follows that $(\log x)/x^{1-a}\ra 0$ as $x\ra\infty$ and
$$
\lim_{x\ra\infty}\frac{\psi(x)}{2x} \geq a\lim_{x\ra\infty}\frac{\pi(x)}{2x/\log x}
$$
for all $0<a<1$.  Hence
$$
\lim_{x\ra\infty}\frac{\psi(x)}{2x} \geq \lim_{x\ra\infty}\frac{\pi(x)}{2x/\log x}.
$$
Together with (\ref{eqn:PsiPiIneq}) and Proposition \ref{pro:PsiHatEst} this proves the proposition.
\qed

\begin{proposition}
\label{pro:PiTildeEst}
$$
\displaystyle\lim_{x\ra\infty}\frac{\tilde{\pi}(x)}{8x/\log x} = \lim_{x\ra\infty}\frac{\tilde{\psi}(x)}{8x} = 1.
$$
\end{proposition}
\prf
Exactly as for the previous proposition, making use of Proposition \ref{pro:PsiTildeEst} and Lemma \ref{lem:SigmaTilde}.
\qed

\begin{proposition}
\label{pro:PiEst2}
$$
\pi(x)=2\li(x)+O\left(\frac{x^{3/4}}{\log x}\right).
$$
\end{proposition}
\prf
We consider the function
\begin{eqnarray*}
S(x) & = & \sum_{N(\ga)\leq x}\chi_1(\Ga_{\ga_0})\frac{l_{\ga_0}}{l_{\ga}} \\
  & = & \sum_{N(\ga_0)\leq x} \chi_1(\Ga_{\ga_0}) + \sum_{k\geq 2}\sum_{N(\ga_0)\leq x^{1/k}} \chi_1(\Ga_{\ga_0})\rez{k}.
\end{eqnarray*}
We consider the double sum on the right.  Since $\Ga\subset G$ is discrete there is a geodesic $\ga_{\rm min}$ of minimum length.  For a given $x$ the inner sum contains at least one summand only for $k\leq \log x/l_{\ga_{\rm min}}$.  For each such $k\geq 2$, by Proposition \ref{pro:PiEst}, the inner sum is equal to $O(\sqrt{x}/\log x)$.  Therefore we have
\begin{equation}
\label{eqn:SPiEst}
S(x) = \pi(x)+O\left(\sqrt{x}\right).
\end{equation}
Now
\begin{eqnarray*}
\int_2^x \frac{\psi(t)}{t\log^2\! t}\ dt & = & \int_2^x \sum_{N(\ga)\leq t} \chi_1(\Ga_{\ga_0})\frac{l_{\ga_0}}{t\log^2\! t}\ dt \\
  & = & \sum_{N(\ga)\leq x} \int_{N(\ga)}^x \chi_1(\Ga_{\ga_0})\frac{l_{\ga_0}}{t\log^2 t}\ dt \\
  & = & \sum_{N(\ga)\leq x} \chi_1(\Ga_{\ga_0})l_{\ga_0}\left(\rez{l_{\ga}}-\rez{\log x}\right) \\
  & = & S(x) - \frac{\psi(x)}{\log x}.
\end{eqnarray*}
Hence
$$
S(x)=\int_2^x \frac{\psi(t)}{t\log^2\! t}\ dt + \frac{\psi(x)}{\log x}.
$$
By Proposition~\ref{pro:PsiHatEst} we have
\begin{eqnarray*}
S(x) & = & \int_2^x \frac{2}{\log^2 t}\ dt + \frac{2x}{\log x} + O\left(\int_2^x \rez{t^{1/4}\log^2 t}\ dt\right) + O\left(\frac{x^{3/4}}{\log x}\right) \\
     & = & \left[-\frac{2t}{\log t}\right]_2^x + \int_2^x\frac{2}{\log t}\ dt + \frac{2x}{\log x} + O\left(\frac{x^{3/4}}{\log x}\right) \\
     & = & \int_2^x\frac{2}{\log t}\ dt + O\left(\frac{x^{3/4}}{\log x}\right).
\end{eqnarray*}
Together with (\ref{eqn:SPiEst}) this proves the proposition.
\qed

\begin{proposition}
\label{pro:PiTildeEst2}
$$
\tilde{\pi}(x)=8\li(x)+O\left(\frac{x^{3/4}}{\log x}\right).
$$
\end{proposition}

\prf
Exactly as for the previous proposition, making use of Proposition \ref{pro:PsiTildeEst} and Proposition \ref{pro:PiTildeEst}.
\qed

The proof of Theorem \ref{thm:PGT2} proceeds exactly as for the proof of Proposition \ref{pro:PiEst}, making use of Proposition \ref{pro:PsiEst}.

  \chapter{Division Algebras of Degree Four}
    \label{ch:DivAlg}
    \markright{\textnormal{\thechapter{. Division Algebras of Degree Four}}}
    \section{Central simple algebras and orders}
Let $F$ be a field.  An \emph{algebra} \index{algebra} $A$ over $F$ is a (not necessarily commutative) ring with unity, which is also a vector space over $F$, such that $x(ab)=(xa)b=a(xb)$ for all $x\in F$ and $a,b\in A$.  The field $F$ can be embedded into the centre of an $F$-algebra $A$ under the map $x\mapsto x.1$, where $1$ is the unity element of $A$.  A ring is said to be \emph{simple} \index{algebra!simple} if it has no non-trivial, proper, two-sided ideals.  If $A$ is simple and its centre is equal to $F$ then it is known as a \emph{central simple} \index{algebra!central simple} $F$-algebra.

If for every non-zero element of $A$ there exists a (two-sided) inverse in $A$, then $A$ is called a \emph{division algebra} \index{algebra!division}.  The centre of a division algebra is a field.  Every division algebra is simple and hence central simple over its centre.  Let $\CCC(F)$ \index{$\CCC(F)$} denote the set of isomorphism classes of finite dimensional central simple algebras over $F$ and $\CCD(F)$ \index{$\CCD(F)$} the set of isomorphism classes of finite dimensional division algebras with centre $F$, then $\CCD(F)\subset\CCC(F)$.

For convenience we collect together here a number of facts about central simple algebras and their orders, which will be used in the sequel.  Let $F$ be a field, $K$ an extension of $F$ and $A$ a finite dimensional algebra over $F$.

\begin{proposition}
\label{pro:algebras}

(a) If $n=\dim_F A$, then there exists an injective homomorphism $A\hra\Mat_n(F)$.

(b) $\dim_K(A\ox K)=\dim_F A$.

(c) If $A\in\CCC(F)$, then $A\ox K\in\CCC(K)$.

(d) If $B$ is a simple subalgebra of $A\in\CCC(F)$ with centraliser $C_A(B)$ in A, then $(\dim_F B)(\dim_F C_A(B))=\dim_F A$.
\end{proposition}
\prf
All the statements of the proposition are proved in \cite{Pierce82}.  Statement (a) is Corollary 5.5b; (b) is Lemma 9.4; (c) is Proposition 12.4b(ii) and (d) is Theorem 12.7.
\qed

The next proposition gives some results about subfields of algebras.

\begin{proposition}
\label{pro:subfields}
(a) If $A\in\CCD(F)$, then for every $x\in A$ the set $F[x]=\left\{f(x):f\in F[X]\right\}$ is a subfield of $A$.

(b) If $A\in\CCC(F)$, then $\dim_F A =m^2$ for some $m\in\N$.  The natural number $m$ is called the degree of $A$.  If $E$ is a subfield of $A$ then $[E:F]$ divides $m$.

If $K$ is a subfield of $A$ containing $F$ then it is said to be \emph{strictly maximal} \index{strictly maximal subfield} if $[K:F]=\deg A$.

(c) A subfield $K$ of $A\in\CCC(F)$ is strictly maximal if and only if $C_A(K)=K$.  If $A\in\CCD(F)$ then every maximal subfield of $A$ is strictly maximal.

(d) If $A\in\CCC(F)$ and $[K:F]=\deg A=n$, then $A\ox K\cong\Mat_n(K)$ if and only if $K$ is isomorphic as an $F$-algebra to a strictly maximal subfield of $A$.
\end{proposition}
\prf
All the statements of the proposition are proved in \cite{Pierce82}.  Statement (a) is Lemma 13.1b; (b) is Corollary 13.1a; (c) is Corollary 13.1b and (d) is Corollary 13.3.
\qed

Now let $F$ be a number field and $\CO_F$ the ring of integers in $F$.  For any finite dimensional vector space $V$ over $F$, a \emph{full $\CO_F$-lattice} \index{full $\CO_F$-lattice} in $V$ is a finitely generated $\CO_F$-submodule $M$ in $V$ which contains a basis of $V$.  A subring $\CO$ of $A$ which is also a full $\CO_F$-lattice in $A$ is called an \emph{$\CO_F$-order}\index{order!$\CO_F$-order}, or simply an \emph{order}\index{order}, in $A$.  A \emph{maximal} \index{order!maximal} $\CO_F$-order in $A$ is an $\CO_F$-order which is not properly contained in any other $\CO_F$-order in $A$.  An element $a\in A$ is said to be \emph{integral} \index{integral element} over $\CO_F$ if it is a root of a monic polynomial with coefficients in $\CO_F$.  The \emph{integral closure} \index{integral closure} of $\CO_F$ in $A$ is the set of all elements of $A$ which are integral over $\CO_F$.

\begin{proposition}
\label{pro:orders}
(a)(Skolem-Noether Theorem) \index{Skolem-Noether Theorem} Let $A\in\CCC(F)$, and $B$ a simple subring of $A$ such that $F\subset B\subset A$.  Then every $F$-isomorphism of $B$ onto a subalgebra of $A$ extends to an inner automorphism of $A$.

(b) Every element of an $\CO_F$-order $\CO$ is integral in $A$ over $\CO_F$.

(c) The $\CO_F$-order $\CO$ is maximal in $A$ if and only if for each prime ideal $\p$ of $\CO_F$ the $\p$-adic completion $\CO_{\p}$ is a maximal $\CO_{F,\p}$-order in $A_{\p}$.

(d) If $A\in\CCD(F_{\p})$ for some prime ideal $\p$ of $\CO_F$, then the integral closure of $\CO_{F,\p}$ in $A$ is the unique maximal $\CO_{F,\p}$-order in $A$.
\end{proposition}
\prf
All the statements of the proposition are proved in \cite{Reiner75}.  Statement (a) is Theorem 7.21; (b) is Theorem 8.6; (c) is Corollary 11.6 and (d) is Theorem 12.8.
\qed

\section{Division algebras of degree four}

Let $A$ be a central simple algebra over a field $F$.  By the Wedderburn Structure Theorem \index{Wedderburn Structure Theorem} (\cite{Pierce82}, Theorem 3.5) the algebra $A$ is isomorphic to $\Mat_n(D)$, where $n\in\N$ and $D$ is a finite dimensional division algebra over $F$, which by \cite{Pierce82}, Proposition 12.5b is unique up to isomorphism.  If $B$ is another central simple algebra over $F$ isomorphic to $\Mat_m(E)$, where $m\in\N$ and $E$ is a finite dimensional division algebra over $F$, then $A$ and $B$ are called \emph{Morita equivalent} \index{Morita equivalent} if and only if $D\cong E$.  By \cite{Pierce82}, Proposition 12.5a, the tensor product over $F$ induces a group structure on the set of equivalence classes.  The group defined in this way is called the \emph{Brauer group}.\index{Brauer group}

By \cite{Pierce82}, Theorem 17.10, if $F$ is a local, non-archimedean field, then there is a canonical isomorphism between the Brauer group of $F$ and the group $\Q/\Z$.  The \emph{Brauer invariant} \index{Brauer invariant} $\Inv A$ \index{$\Inv A$} of $A$ is defined to be the image in $\Q/\Z$ of the Morita equivalence class of $A$ under this isomorphism.  The \emph{Schur index} \index{Schur index} $\Ind A$ \index{$\Ind A$} of $A$ is defined to be the degree of $D$.  By \cite{Pierce82}, Corollary 17.10a(iii), the order of the Brauer invariant of $A$ in $\Q/\Z$ equals $\Ind A$.

The only (associative) division algebras over $\R$ are $\R$ itself, $\C$ and the quaternions $\H$ (see \cite{Pierce82}, Corollary 13.1c).  Since $\C$ is itself a field, any algebra which is central simple over $\R$ cannot be isomorphic to a matrix algebra over $\C$.  The Brauer group of $\R$ is therefore isomorphic to $\Z/2\Z$.  Let $A$ be a central simple algebra over $\R$.  If $A\cong\Mat_n(\R)$ for some $n\in\N$, then the Brauer invariant of $A$ is defined to be zero and the Schur index of $A$ is defined to be one.  If $A\cong\Mat_n(\H)$ for some $n\in\N$, then the Brauer invariant of $A$ is one half and the Schur index of $A$ is two.

Let $M$ be a division algebra of degree 4 over $\Q$.  Fix a maximal order $M(\Z)$ \index{$M(\Z)$} in $M$.  If $R$ is a commutative ring with unit then we denote by $M(R)$ \index{$M(R)$} the $R$-algebra $M(\Z)\ox_{\Z} R$.  In particular we have $M\cong M(\Q)$ \index{$M(\Q)$}.  For any ring $R$ the reduced norm induces a map $\det\!\!:M(R)\ra R$.  (We note that in the case that $M(R)\cong\Mat_4(R)$ the reduced norm of an element of $M(R)$ equals its determinant, justifying the choice of notation; see \cite{Pierce82}, Chapter 16.)

Let $p$ be a prime.  From Proposition \ref{pro:algebras}(b) we know that $\deg M(\Q_p)=\deg M(\Q)=4$.  It then follows that if $\Inv M(\Q_p)$ is equal to $\rez{4}$ or $\frac{3}{4}$ then $M(\Q_p)$ is a division algebra of degree four over $\Q_p$, and if $\Inv M(\Q_p)=0$ then $M(\Q_p)\cong\Mat_4(\Q_p)$.  The other possibility, that $\Inv M(\Q_p)=\rez{2}$, does not interest us here.  Proposition \ref{pro:algebras}(b) also tells us that $\deg M(\R)=4$.  We know that $\Inv M(\R)$ equals zero or one half.  It follows that $M(\R)$ is isomorphic to $\Mat_4(\R)$ or $\Mat_2(\H)$ respectively.

Let $S$ \index{$S$} be a finite, non-empty set of prime numbers with an even number of elements.  We say that $M(\Q)$ \emph{splits} \index{splitting!of division algebra} over a prime $p$ if $M(\Q_p)\cong\Mat_4(\Q_p)$.  For all primes $p$ we define $i_p$ as follows.  If $p\in S$ define $i_p$ to be either $\rez{4}$ or $\frac{3}{4}$ in such a way that $\sum_{p\in S}i_p\in\Z$.  Note that we are able to do this since we have specified that $S$ must contain an even number of elements.  For all other $p$ define $i_p=0$.  Then \cite{Pierce82}, Theorem 18.5 tells us that we may choose $M$ so that $\Inv M(\Q_p)=i_p$ for all primes $p$ and $\Inv M(\R)=0$.  We can then see that the set of places at which $M(\Q)$ does not split coincides with the set $S$.  More particularly, for $p\in S$ we have that $M(\Q_p)$ is a division $\Q_p$-algebra, for $p\notin S$ we have that $M(\Q_p)\cong\Mat_4(\Q_p)$ and we also have $M(\R)\cong\Mat_4(\R)$.

For a commutative ring $R$ with unity, let \index{$\CG(R)$}
$$
\CG(R)=\left\{x\in M(R):\det(x)=1\right\}.
$$
Then $\CG$ \index{$\CG$} is a linear algebraic group, defined over $\Z$.  Let $\Ga=\CG(\Z)$, \index{$\Ga$} then $\Ga$ forms a discrete subgroup of $G=\CG(\R)\cong\SL_4(\R)$. Since $M(\Q)$ is a division algebra, it follows that $\CG$ is anisotropic over $\Q$ and so $\Ga$ is cocompact in $G$ (see \cite{BorelHarder78}, Theorem A).  Let $P$ be the parabolic subgroup \index{$P$}
$$
P=\matrix{*}{*}{\begin{array}{cc} 0 & 0 \\ 0 & 0 \end{array}}{*}
$$
of G.  Then $P=MAN$, where \index{$M$}
\begin{eqnarray*}
M & = & {\rm S}\matrixtwo{\SL_2^{\pm}(\R)}{\SL_2^{\pm}(\R)} \\
  & = & \left\{\matrixtwo{X}{Y}:X,Y\in\Mat_2(\R),\ \det X=\det Y=\pm 1\right\},
\end{eqnarray*}
$A=\left\{ \diag (a,a,a^{-1},a^{-1}):a>0\right\}$ \index{$A$} and the elements of $N$ \index{$N$} have ones on the diagonal and the only other non-zero entries in the top right two by two square.  Let \index{$B$}
$$
B=\matrixtwo{\SO(2)}{\SO(2)}.
$$
Then $B$ is a compact subgroup of $M$.

Let $A^-=\left\{ \diag(a,a,a^{-1},a^{-1}):a\in (0,1)\right\}$ \index{$A^-$} be the negative Weyl chamber of $A$.  Let $\CE_P(\Ga)$ \index{$\CE_P(\Ga)$} be the set of conjugacy classes $[\ga]$ in $\Ga$ such that $\ga$ is conjugate in $G$ to an element $a_{\ga}b_{\ga}$ of $A^- B$ and let $\CE_P^p(\Ga)$ \index{$\CE_P^p(\Ga)$} be the subset of primitive conjugacy classes.

We say $g\in G$ is \emph{regular} \index{regular} if its centraliser is a torus and \emph{non-regular} \index{non-regular} (or \emph{singular}) \index{singular} otherwise.  Clearly, for $\ga\in\Ga$ regularity is a property of the $\Ga$-conjugacy class $[\ga]$, we denote by $\CE_P^{p,\reg}(\Ga)$ \index{$\CE_P^{p,\reg}(\Ga)$} the regular elements of $\CE_P^p(\Ga)$.  By abuse of notation we sometimes write $\ga\in\CE_P^p(\Ga)$ or $\ga\in\CE_P^{p,\reg}(\Ga)$ when we mean the conjugacy class of $\ga$.

We will call a quartic field extension $F/\Q$ \emph{totally complex} \index{totally complex field} if it has two pairs of conjugate complex embeddings and no real embeddings into $\C$.

\section{Subfields of $M(\Q)$ generated by $\Ga$}
\begin{lemma}
\label{lem:centralisers}
Let $[\ga]\in\CE_P(\Ga)$.

The centraliser $M(\Q)_{\ga}$ of $\ga$ in $M(\Q)$ is a totally complex quartic field if and only if $\ga$ is regular.

The centraliser $M(\Q)_{\ga}$ is a quaternion algebra over the real quadratic field $\Q[\ga]$ if and only if $\ga$ is non-regular.
\end{lemma}

\prf
The centraliser $M(\Q)_{\ga}$ of $\ga$ in $M(\Q)$ is a division subalgebra of $M(\Q)$.  Suppose that $M(\Q)_{\ga}=\Q$.  Then since $\ga\in M(\Q)_{\ga}$ we must have $\ga\in\Q$, but then $\ga$ is central so $M(\Q)_{\ga}=M(\Q)$, a contradiction.  Suppose instead that $M(\Q)_{\ga}=M(\Q)$, ie. $\ga$ is central.  Then $\ga\in\Q$, so $\det\ga=\ga^4=1$ and it follows that $\ga=\pm 1$, which possibility is excluded since $\ga=\pm 1$ is not primitive.

By Proposition \ref{pro:algebras}(d) we have that $\dim_{\Q}M(\Q)_{\ga}$ divides $\dim_{\Q}M(\Q)$, which is 16, so $\dim_{\Q}M(\Q)_{\ga}$ is equal to 2, 4 or 8.  Also
$$
\dim_{\Q}M(\Q)_{\ga}=\dim_{\R}M(\R)_{\ga}=\dim_{\R}M(\R)_{a_{\ga}b_{\ga}}.
$$
Depending on whether neither, one or both of the $\SO(2)$ components of $b_{\ga}$ are $\pm I_2$, the dimension $\dim_{\R}M(\R)_{a_{\ga}b_{\ga}}$ is equal to 4, 6 or 8 respectively.  Hence $\dim_{\Q}M(\Q)_{\ga}$ is either 4 or 8.  If $\dim_{\Q}M(\Q)_{\ga}=4$ then $b_{\ga}=(R(\th),R(\phi))$ with $\th,\phi\notin\pi\Z$, where
$$
R(\th)=\matrix{\cos\th}{-\sin\th}{\sin\th}{\cos\th}.
$$
If $\dim_{\Q}M(\Q)_{\ga}=8$ then $b_{\ga}$ is diagonal.

In the first case $G_{\ga}\cong AB\cong\R^+\x\SO(2)\x\SO(2)$ so $\ga$ is regular.  From Proposition \ref{pro:algebras}(d) and Proposition \ref{pro:subfields}(a) it follows that $\Q[\ga]=\left\{f(\ga):f(x)\in \Q[x]\right\}$ is a subfield of $M(\Q)$ of degree 4.  Let $f_{\ga}(x)$ be the minimal polynomial of $\ga$ over $\Q$.  Then $a_{\ga}b_{\ga}$ also satisfies $f_{\ga}(x)=0$.  Now $a_{\ga}=\diag(a,a,a^{-1},a^{-1})$ for some $a\in(0,1)$ and $b_{\ga}=(R(\th),R(\phi))$ for some $\th, \phi\notin\pi\Z$, so the complex numbers $z_1=ae^{i\th}$ and $z_2=a^{-1}e^{i\phi}$ also satisfy $f_{\ga}(x)=0$.  Neither $z_1$ nor $z_2$ are real, nor do we have $z_1=\bar{z_2}$, so $\Q[\ga]\cong\Q[z_1,z_2]$, which is a totally complex field.  The division subalgebra $M(\Q)_{\ga}$ of $M(\Q)$ has dimension four over $\Q$ and contains $\Q[\ga]$, hence is equal to $\Q[\ga]$.

Suppose instead that $b_{\ga}$ is diagonal.  Then $G_{\ga}\cong AM$, so $\ga$ is not regular.  From Proposition \ref{pro:algebras}(d) and Proposition \ref{pro:subfields} (a) it follows that $\Q[\ga]=\left\{f(\ga):f(x)\in \Q[x]\right\}$ is a subfield of $M(\Q)$ of degree 2.  Let $f_{\ga}(x)$ be the minimal polynomial of $\ga$ over $\Q$.  Then $a_{\ga}b_{\ga}$ also satisfies $f_{\ga}(x)=0$.  Now $a_{\ga}=\diag(a,a,a^{-1},a^{-1})$ for some $a\in(0,1)$ and $b_{\ga}=(\pm I_2,\pm I_2)$.  Let us suppose that $b_{\ga}=I_4$, the other cases are similar.  Then $f_{\ga}(a)=f_{\ga}(a^{-1})=0$, so $f_{\ga}=(x-a)(x-a^{-1})$ and $\Q[\ga]$ is a real quadratic field.

The division algebra $M(\Q)_{\ga}$ is central simple over its centre $Z(M(\Q)_{\ga})$.  Clearly, $\Q[\ga]$ is contained in $Z(M(\Q)_{\ga})$.  By Proposition \ref{pro:subfields}(b), the dimension of $M(\Q)_{\ga}$ over its centre is $m^2$ for some $m\in\N$.  Hence either $M(\Q)_{\ga}$ is a field or $Z(M(\Q)_{\ga})=\Q[\ga]$.  However, by Proposition \ref{pro:subfields}(b) again, if $F$ is a subfield of $M(\Q)$ then $[F:\Q]$ = 1, 2 or 4.  This rules out the possibility of $M(\Q)_{\ga}$ being a field, since $\dim_{\Q}M(\Q)_{\ga}=8$.  We conclude that $Z(M(\Q)_{\ga})=\Q[\ga]$ and by \cite{Pierce82}, Theorem 13.1, $M(\Q)_{\ga}$ is a quaternion algebra over $\Q[\ga]$.
\qed

\section{Field and order embeddings}

In this section we prove a number of lemmas which will be needed later.

\begin{lemma}
\label{lem:unit}
Let $u\in M(\Z)$.  Then $u$ is a unit in $M(\Z)$ iff $\det(u)=\pm 1$.  If $\Q[u]$ is quadratic or is totally complex quartic, or if $u=\pm 1$, then $\det(u)=1$.
\end{lemma}

\prf
By Proposition \ref{pro:subfields}(a) the set $F=\Q[u]$ is a subfield of $M(\Q)$ containing $u$.  The set $F\cap M(\Z)$ is an order in $F$.  The element $u$ is in $F\cap M(\Z)$ so by Proposition \ref{pro:orders}(b) it is in the integral closure $\CO_F$ of $\Z$ in $F$.  

The field norm $N_{F|\Q}:F\ra\Q$ can be written as a product over all embeddings of $F$ into $\C$:
$$
N_{F|\Q}(x)=\prod_{\si}\si x.
$$
Hence we can see that it restricts to a group homomorphism $F^{\x}\ra\Q^{\x}$ and maps integral elements in $F$ to elements of $\Z$ (see \cite{Neukirch99}).  Suppose there exists $v\in\CO_F$ such that $vu=1$.  Then $N_{F|\Q}(v),N_{F|\Q}(u)\in\Z$ and
$$
N_{F|\Q}(v)N_{F|\Q}(u)=N_{F|\Q}(uv)=N_{F|\Q}(1)=1,
$$
so $N_{F|\Q}(u)=\pm 1$.  On the other hand, suppose there exists $n\in\Z$ such that $nN_{F|\Q}(u)=1$.  Then
$$
1=n\prod_{\si}\si u=yu
$$
for some $y\in\CO_F^{\x}$.

We have shown that $u$ is a unit in $\CO_F$ if and only if $N_{F|\Q}(u)=\pm 1$.  Proposition 16.2a of \cite{Pierce82} tells us that
\begin{equation}
\label{eqn:detNormReln}
\det(u)=N_{F|\Q}(u)^k,
\end{equation}
where $\deg M(\Q)=k[F:\Q]$.  The first statement of the lemma follows.

If $u=\pm 1$ (so that $F=\Q$) or $F$ is quadratic then by (\ref{eqn:detNormReln}) we have $\det(u)=1$.  Suppose $F$ is a totally complex quartic extension.  Let $f(X)=X^4+a_3 X^3+a_2 X^2+a_1 X+a_0$ be the minimal polynomial of $u$ over $\Q$.  Then $\det(u)=a_0$.  If $a_0=\det(u)=-1$ then $f$ must have a real root, which contradicts the supposition on $F$.  Hence $\det(u)=1$.
\qed

The following lemma holds for any algebraic extension $F$ of $\Q$.
\begin{lemma}
\label{lem:maxorders}
Let $\CO$ be an order in the number field $F/\Q$ and $p$ a prime number.  Let $F_p=F\ox\Q_p$ \index{$F_p$} and $\CO_p=\CO\ox\Z_p$\index{$\CO_p$}.  Then $\CO_p$ is the maximal order of $F_p$ for all but finitely many primes $p$.
\end{lemma}

\prf
Let $\CO_M$ be the maximal order of $F$ and let $\CO_{M,p}=\CO_M\ox\Z_p$.  We define the \emph{conductor} \index{conductor} of $\CO$ to be the ideal $\CF=\{\al\in\CO_M|\al\CO_M\subset\CO\}$ of $\CO_M$.  Then for all primes $p$, if $\CF\nsubseteq p\CO_M$, then $\CO_p=\CO_{M,p}$.  In fact, if $\CF\nsubseteq p\CO_M$ then there is an element $\al\in\CF$ such that $\al\notin p\CO_M$.  Every element of $\CO_{M,p}$ can be written as a unit of $\CO_{M,p}$ times a power of $p$.  Since $\al\notin p\CO_M$ we have $\al\in\CO_{M,p}^{\x}$, which implies $\al\CO_{M,p}=\CO_{M,p}$.  Then, as $\al\in\CF$ it follows that
$$
\CO_{M,p}=\al\CO_{M,p}\subset\CO_p.
$$
It is clear that $\CO_p\subset\CO_{M,p}$.  Since $\CF\subset p\CO_M$ for only finitely many primes $p$, the lemma follows.
\qed

Let $F$ be a number field.  A prime number $p$ is called {\it non-decomposed} \index{prime!non-decomposed} \index{non-decomposed (prime)} in $F$ if there is only one place in $F$ lying above $p$.

\begin{lemma}
\label{lem:fieldembedding}
A number field $F$ embeds into $M(\Q)$ if and only if $[F:\Q]=1, 2$ or $4$ and $p$ is non-decomposed in $F$ for all $p\in S$.
\end{lemma}

\prf
The statement about the degree of $F$ is Proposition \ref{pro:algebras} (c).  The case $F=\Q$ is trivial.

Let $F$ be a quartic field and suppose that $F$ embeds into $M(\Q)$, then $F$ is a strictly maximal subfield of $M(\Q)$.  We say that a subfield $K$ of $M(\Q)$ is a \emph{splitting field} \index{splitting!field} for $M(\Q)$ (or that $K$ \emph{splits} $M(\Q)$) if $M(\Q)\ox K\cong \Mat_4(K)$.  Proposition \ref{pro:subfields}(d) tells us that a quartic field embeds into $M(\Q)$ if and only if it splits $M(\Q)$.

For $p$ a rational prime, the field $\Q_p$ is an extension of $\Q$.  Hence, by Proposition \ref{pro:algebras}(c), the $\Q_p$-algebra $M(\Q_p)$ is central simple over $\Q_p$.  Let $p\in S$.  The Schur index, $\Ind M(\Q_p)$, of $M(\Q_p)$ is, by \cite{Pierce82}, Corollary 17.10a(iii), equal to the order of $[M(\Q_p)]$ in the Brauer group of $\Q_p$, where the square brackets denote the Morita equivalence class.  By choice of $M(\Q)$ we have $\Ind M(\Q_p)=4$.  Let $\p_i$ be the (finitely many) primes in $F$ above $p$.  Then Theorem 32.15 of \cite{Reiner75} tells us that, since $F$ splits $M(\Q)$,
$$
4=\Ind M(\Q_p)\ |\ [F_{\p_i}:\Q_p]
$$
for all $i$.  Corollary 8.4 of Chapter II of \cite{Neukirch99} says that
\begin{equation}
\label{eqn:padicSum}
[F:\Q]=\sum_i[F_{\p_i}:\Q_p],
\end{equation}
which implies that there is only one prime $\p$ in $F$ above $p$ and that $[F_{\p}:\Q_p]=4$.

Conversely, suppose that for each $p\in S$ there is only one prime $\p$ in $F$ above $p$.  From (\ref{eqn:padicSum}) we get that $[F_{\p}:\Q_p]=4$, so
\begin{equation}
\label{eqn:indexDivides}
\Ind M(\Q_p)\ |\ [F_{\p}:\Q_p].
\end{equation}
For $p\notin S$, by choice of $M(\Q)$, the Schur index, $\ind M(\Q_p)$, is equal to 1, so in this case (\ref{eqn:indexDivides}) also holds for each $\p$ in $F$ above $p$.  Then, by \cite{Reiner75}, Theorem 32.15, the field $F$ splits $M(\Q)$.  By Proposition \ref{pro:subfields}(d) we then deduce that $F$ embeds into $M(\Q)$.

Now let $F$ be a quadratic field and suppose that $F$ embeds into $M(\Q)$, then we consider $F$ as a subfield of $M(\Q)$.  By Proposition \ref{pro:subfields}(c), the subfield $F$ of $M(\Q)$ is not maximal so there exists a maximal subfield $K$ of $M(\Q)$ properly containing $F$, which by the same proposition is quartic.  It was shown above that every prime in $S$ is non-decomposed in $K$ and hence also in $F$.

Suppose $p$ is non-decomposed in $F$ for all $p\in S$.  We choose a quaternion algebra $A$ over $F$ by specifying that the local Brauer invariant at $\p$ be $\rez{2}$ for all places $\p$ in $F$ over some $p\in S$, and that the local Brauer invariant be 0 at all other places of $F$.  Proposition \ref{pro:subfields}(c) tells us that there exists a subfield $K$ of $A$ containing $F$ with $[K:F]=2$.  Then by the same argument as above, if $\p$ is a place of $F$ over $p$ for some $p\in S$, then $\p$ is non-decomposed in $K$.  It follows that $p$ is non-decomposed in $K$ for all $p\in S$.  It was shown above that the quartic extension $K$ of $\Q$, and thus also the subfield $F$ of $K$, embeds into $M(\Q)$.  The lemma is proven.
\qed

Let $\A_{\fin}$ \index{$\A_{\fin}$} denote the ring of finite adeles over $\Q$, ie.
$$
\A_{\fin}=\left.\prod_p\right.'\Q_p,
$$
where $p$ ranges over the rational primes and $\prod'$ denotes the restricted product, that is, almost all components of a given element are integral.  Let $\hat{\Z}=\prod_p\Z_p\subset\A_{\fin}$\index{$\hat{\Z}$}.  There is a natural embedding of $\Q$ into $\A_{\fin}$, which maps $q\in\Q$ to the adele every one of whose components is $q$.  From this embedding we derive an embedding of $M(\Q)$ in $M(\A_{\fin})$ and an embedding of $\CG(\Q)$ in $\CG(\A_{\fin})$.

\begin{lemma}
\label{lem:Adeles}
$M(\A_{\fin})^{\x}=M(\hat{\Z})^{\x}\,M(\Q)^{\x}$.
\end{lemma}
\prf
We have the following commutative diagram with exact rows:
\begin{displaymath}
\xymatrix{1 \ar[r] & \CG(\Q) \ar[r] \ar[d] & M(\Q)^{\x} \ar[r]^{\det} \ar[d] & \Q^{\x} \ar[r] \ar[d] & 1 \\
   1 \ar[r] & \CG(\A_{\fin}) \ar[r] & M(\A_{\fin})^{\x} \ar[r]^{\det} & \A_{\fin}^{\x} \ar[r] & 1 \\
   1 \ar[r] & \CG(\hat{\Z}) \ar[r] \ar[u] & M(\hat{\Z})^{\x} \ar[r]^{\det} \ar[u] & \hat{\Z}^{\x} \ar[r] \ar[u] & 1}
\end{displaymath}
where the vertical arrows denote the natural embeddings.  The commutativity of the diagram and the exactness of the rows are clear, except that the surjectivity of the three determinant maps requires justification.  The surjectivity of the map $\det\!\!:M(\Q)^{\x}\ra\Q^{\x}$ follows from \cite{Reiner75}, Theorem 33.15.  For all primes $p$ we have that $M(\Q_p)$ is a central simple $\Q_p$-algebra so it follows from \cite{Reiner75}, Theorem 33.4 that the map $\det\!\!:M(\Q_p)^{\x}\ra\Q_p^{\x}$ is surjective.  For $p\in S$, we have that $M(\Q_p)$ is a division $\Q_p$-algebra and so an element $x$ of $M(\Q_p)$ is in $M(\Z_p)$ if and only if $\det x\in\Z_p$ (\cite{Reiner75}, Theorem 12.5).  For $p\notin S$, we have that $M(\Q_p)\cong\Mat_4(\Q_p)$ and
$$
\det\matrixfour{x}{1}{1}{1}=x
$$
for all $x\in\Z_p$.  It follows that the map $\det\!\!:M(\hat{\Z})^{\x}\ra\hat{Z}^{\x}$ is surjective.  From this we can then see that map $\det\!\!:M(\A_{\fin})^{\x}\ra\A_{\fin}^{\x}$ is also surjective.

Let $v\in M(\A_{\fin})^{\x}$ and let $w=\det v\in\A_{\fin}^{\x}$.  Since $\A_{\fin}^{\x}=\hat{\Z}^{\x}\Q^{\x}$, there exist $r\in\hat{\Z}^{\x}$ and $q\in\Q^{\x}$ such that $w=rq$.  By the surjectivity of $\det$, there exist $\bar{r}\in M(\hat{\Z})^{\x}$ and $\bar{q}\in M(\Q)^{\x}$ such that $\det\bar{r}=r$ and $\det\bar{q}=q$.  Let $\bar{v}=\bar{r}^{-1}v\bar{q}^{-1}$.  Then $\det\bar{v}=1$, so $\bar{v}\in\CG(\A_{\fin})$.  By the Strong Approximation Theorem (see \cite{Kneser66}), we have that $\CG(\Q)$ is dense in $\CG(\A_{\fin})$.  Hence there exist $\hat{r}\in\CG(\hat{\Z})$ and $\hat{q}\in\CG(\Q)$ such that $\bar{v}=\hat{r}\hat{q}$.  Finally, we have
$$
v=\bar{r}\bar{v}\bar{q}=\bar{r}(\hat{r}\hat{q})\bar{q}=(\bar{r}\hat{r})(\hat{q}\bar{q}),
$$
where $\bar{r}\hat{r}\in M(\hat{\Z})^{\x}$ and $\hat{q}\bar{q}\in M(\Q)^{\x}$.  The lemma follows.
\qed

\begin{lemma}
\label{lem:orderembedding}
Let $F$ be a field that embeds into $M(\Q)$.  For each embedding $\si:F\ra M(\Q)$ the order $\CO_{\si}=\si^{-1}\left(\si(F)\cap M(\Z)\right)$ \index{$\CO_{\si}$} is maximal at all $p\in S$.  Conversely, if $\CO$ is an order in the field $F$ which is maximal at all $p\in S$, then there exists an embedding $\si:F\ra M(\Q)$ such that $\CO=\CO_{\si}$.
\end{lemma}

Let $\si:F\ra M(\Q)$ be an embedding of the field $F$ into $M(\Q)$.  By Lemma \ref{lem:fieldembedding} each prime $p$ in the set $S$ is non-split in $F$.  Since for $p\in S$ there is only one place in $F$ above $p$, we can note that the $p$-adic completion $F_p=F\ox\Q_p$ is again a field (see \cite{Neukirch99}, Chapter 1, Proposition 8.3).  Firstly we need to show that for all $p\in S$, the completion $\CO_{\si,p}=\CO_{\si}\ox\Z_p$ is the maximal order in $F_p$.  By Proposition \ref{pro:orders}(c), $M(\Z_p)$ is the maximal $\Z_p$-order in $M(\Q_p)$.  Since $M(\Q_p)$ is a division algebra, $M(\Z_p)$ coincides with the integral closure of $\Z_p$ in $M(\Q_p)$, by Proposition \ref{pro:orders}(d).  We can extend $\si$ to an embedding of $F_p$ in $M(\Q_p)$ and then $\CO_{\si,p}=\si^{-1}\left(\si(F_p)\cap M(\Z_p)\right)$ is the integral closure of $\Z_p$ in $F_p$, which is the maximal order of $F_p$.

For the converse, let $\CO$ be an order of $F$ such that the completion $\CO_{p}$ is maximal for each $p\in S$.  Via $\si$ we consider $F$ to be a subfield of $M(\Q)$.  For any $u\in M(\Q)^{\x}$ let $\CO_u =F\cap u^{-1}M(\Z)u$.  Let $^u\si=u\si u^{-1}$, that is: $^u\si(x)=u\si(x)u^{-1}$ for all $x\in F$.  We will show that there is a $u\in M(\Q)^{\x}$ such that $\CO=\CO_u$.  The embedding $^u\si$ is then the one required by the lemma.

Let $\CO_1=F\cap M(\Z)$.  Since $\CO$ and $\CO_1$ are orders, by Lemma \ref{lem:maxorders} they are both maximal at all but finitely many places.  Let $T$ be the set of primes $p$ such that either $\CO$ or $\CO_1$ is not maximal at $p$.  Then $T$ is finite and $T\cap S=\varnothing$.  Furthermore, if $T_F$ denotes the set of places of $F$ lying over $T$, we have that for any place $\p$ of $F$ with $\p\notin T_F$ the completions $\CO_{\p}$ and $\CO_{1,\p}$ are maximal in $F_{\p}$ and thus equal, by uniqueness of maximal orders.  Hence, for $p\notin T$ we have
$$
\CO_p=F_p\cap M(\Z_p).
$$

Let $p\in T$, then $p\notin S$ so $M(\Q_p)\cong\Mat_4(\Q_p)$.  Considering $F_p$ as a $\Q_p$-vector space we see that we can embed $F_p$ into $\Lin_{\Q_p}(F_p,F_p)$ by sending an element $x\in F_p$ to the linear map $a\mapsto ax$.  We choose a $\Z_p$-basis of $\CO_p$, which is then also a $\Q_p$-basis of $F_p$.  This gives us an isomorphism
$$
\Lin_{\Q_p}(F_p,F_p)\cong\Mat_4(\Q_p)\cong M(\Q_p).
$$
This isomorphism and the above embedding $F_p\hra\Lin_{\Q_p}(F_p,F_p)$ then give us an embedding $\si_p:F_p\hra M(\Q_p)$ such that $\CO_p=\si_p^{-1}\left(\si_p(F_p)\cap M(\Z_p)\right)$.  The map $\si_p$ is a $\Q_p$-isomorphism of $F_p$ (considered as a subfield of $M(\Q_p)$ via a suitable extension of $\si$) onto its image in $M(\Q_p)$, so by the Skolem-Noether Theorem (Proposition \ref{pro:orders}(a)), there exists $u_p\in M(\Q_p)^{\x}$ such that $u_pxu_p^{-1}=\si_p(x)$ for all $x\in F_p$.  Hence
$$
\CO_p=F_p\cap u_p^{-1}M(\Z_p)u_p.
$$

For $p\notin T$ set $u_p=1$ and let $\ut=\left(u_p\right)\in\A_{\fin}$.  By Lemma \ref{lem:Adeles}
$$
M(\A_{\fin})^{\x}=M(\hat{\Z})^{\x}\,M(\Q)^{\x}
$$
so there exists an element $u\in M(\Q)^{\x}$ such that for all primes $p$ we have
$$
u^{-1}M(\Z_p)u=u_p^{-1}M(\Z_p)u_p.
$$
This implies that
\begin{equation}
\label{eqn:padicid}
\CO_p=F_p\cap u^{-1}M(\Z_p)u
\end{equation}
for all primes $p$ and we deduce that
$$
\CO=F\cap u^{-1}M(\Z)u=\CO_u.
$$
This completes the proof of the lemma.
\qed

\section{Counting order embeddings}

Let $F/\Q$ be an algebraic number field and let $\CO$ be an order in $F$.  Let $I(\CO)$ \index{$I(\CO)$} be the set of all finitely generated $\CO$-submodules of $F$.  According to the Jordan-Zassenhaus Theorem \index{Jordan-Zassenhaus Theorem} (\cite{Reiner75}, Theorem 26.4), the set $[I(\CO)]$ of isomorphism classes of elements of $I(\CO)$ is finite.  Let $h(\CO)$ \index{$h(\CO)$} be the cardinality of the set $[I(\CO)]$, called the \emph{class number} \index{class number} of $\CO$.

Let $I$ be a non-trivial, finitely generated $\CO$-submodule of $F$.  Then, by \cite{Reiner75}, Theorem 10.6, $I$ is isomorphic as an $\CO$-module to an ideal of $\CO$ so also has the property that $I\ox\Q=F$.  If $J$ is another finitely generated $\CO$-submodule of $F$, such that $J\cong I$, then the isomorphism extends to an automorphism of $F$ considered as a module over itself, and hence $J=Ix$ for some $x\in F^{\x}$.  So $h(\CO)=\card([I(\CO)])=|I(\CO)/F^{\x}|$.

For a prime $p\in\Z$ let $F_p=F\ox\Q_p$\index{$F_p$}.  Let $\p_1,...\p_n$ be the prime ideals of $F$ above $p$.  By \cite{Neukirch99}, Chapter 1, Proposition 8.3 there exists a canonical $\Q_p$-algebra isomorphism defined by:
\begin{eqnarray}
\label{eqn:localproduct}
F_p=F\ox\Q_p & \cong   & \prod_{i=1}^n F_{\p_i} \\
    \al\ox x   & \mapsto & \left(\tau_i(\al).x\right)_{i=1}^n, \nonumber
\end{eqnarray}
where $F_{\p_i}$ denotes the completion of $F$ with respect to the $\p_i$-adic absolute value and $\tau_i$ the embedding of $F$ into $F_{\p_1}$.  If $\CO$ is an order in $F$ and $p\in\Z$ a prime then let $\CO_p=\CO\ox\Z_p$\index{$\CO_p$}.

Recall that $S$ is a finite, non-empty set of prime numbers with an even number of elements.  Let $C(S)$ \index{$C(S)$} be the set of all totally complex quartic fields $F$ such that all primes $p$ in $S$ are non-decomposed in $F$.  Note in particular that by the isomorphism (\ref{eqn:localproduct}) this means that for each field $F$ in $C(S)$ and for each $p\in S$ we have that $F_p=F\ox\Q_p$ is once again a field, namely the completion of $F$ at the unique place of $F$ above $p$.  For $F\in C(S)$ let $O_F(S)$ \index{$O_F(S)$} be the set of isomorphism classes of orders $\CO$ in $F$ which are maximal at all $p\in S$, ie. are such that the completion $\CO_p=\CO\ox\Z_p$ is the maximal order of the field $F_p$ for all $p\in S$.  Let $O(S)$ \index{$O(S)$} be the union of all $O_F(S)$, where $F$ ranges over $C(S)$.

Let $F\in C(S)$ and $\CO\in O_F(S)$, then, by Lemma \ref{lem:fieldembedding}, $F$ embeds into $M(\Q)$.  By Lemma \ref{lem:unit}, the group $M(\Z)^{\x}$ of units in $M(\Z)$ contains $\Ga$ as a subgroup.  By Lemma \ref{lem:orderembedding} we know that there is an embedding $\si$ of $F$ into $M(\Q)$ such that $\CO = \CO_{\si} = \si^{-1}\left(\si(F)\cap M(\Z)\right)$.  Let $u\in M(\Z)^{\x}$ and, as above, let $^u\si=u\si u^{-1}$ for embeddings $\si:F\ra M(\Q)$.  Then $\CO_{^u\si}=\CO_{\si}$, so the group $M(\Z)^{\x}$ acts on the set $\Si(\CO)$ \index{$\Si(\CO)$} of all embeddings $\si$ with $\CO_{\si}=\CO$.  This $M(\Z)^{\x}$ action also restricts to an action of $\Ga$ on $\Si(\CO)$.

We define \index{$\la_S(\CO)$}
$$
\la_S(\CO)=\la_S(F)=\prod_{p\in S} f_p(F),
$$
where $f_p(F)$ is the inertia degree of $p$ in $F$.  Our aim is to prove the following proposition.

\begin{proposition}
\label{pro:embeddings}
The quotient $M(\Z)^{\x}\bs\Si(\CO)$ is finite and has cardinality equal to the product $h(\CO)\la_S(\CO)$.
\end{proposition}

We prepare for the proof of the propostition with the following discussion and then prove some lemmas, from which the proposition will follow.

Fix an embedding $\si:F\hra M(\Q)$ with $\CO=\CO_{\si}$ and consider $F$ as a subfield of $M(\Q)$ such that $\CO=F\cap M(\Z)$.  For $u\in M(\Q)^{\x}$ let
$$
\CO_u =F\cap u^{-1}M(\Z)u.
$$
Let $U$ be the set of all $u\in M(\Q)^{\x}$ such that
$$
\CO=F\cap M(\Z)=F\cap u^{-1}M(\Z)u=\CO_u.
$$
Let $a\in F^{\x}$, then
$$
F\cap a^{-1}u^{-1}M(\Z)ua=a^{-1}\left(F\cap u^{-1}M(\Z)u\right)a=F\cap u^{-1}M(\Z)u,
$$
so $F^{\x}$ acts on $U$ by multiplication from the right.  Let $v\in M(\Z)^{\x}$, then
$$
F\cap u^{-1}v^{-1}M(\Z)vu=F\cap u^{-1}M(\Z)u,
$$
so $M(\Z)^{\x}$ acts on $U$ by multiplication from the left.  Let $\tau\in\Si(\CO)$, so that $\CO_{\tau}=\CO$.  Then by the Skolem-Noether Theorem (Proposition \ref{pro:orders}(a)) there exists $u\in M(\Q)$ such that $\tau=\,^u\si=u\si u^{-1}$, ie.$\tau(F)=uFu^{-1}$.  Then 
\begin{eqnarray*}
\CO_u = F\cap u^{-1}M(\Z)u & = & u^{-1}\left(uFu^{-1}\cap M(\Z)\right)u \\
  & = & u^{-1}\left(\tau(F)\cap M(\Z)\right)u \\
  & = & \CO_{\tau} \\
  & = & \CO.
\end{eqnarray*}
Conversely, it is clear that for all $u\in U$ we have $^u\si\in\Si(\CO)$, and that $^u\si=\ ^v\si$ if and only if $u=vx$ for some $x\in F^{\x}$.  Hence
$$
|M(\Z)^{\x}\bs U/F^{\x}|=|M(\Z)^{\x}\bs \Si(\CO)|,
$$
and the proposition will be proved if we can show that the left hand side equals $h(\CO)\la_S(\CO)$.

For $u\in U$ let
$$
I_u=F\cap M(\Z)u.
$$
Then $I_u$ is a finitely generated $\CO$-module in $F$.  We shall prove that the map \index{$\Psi$}
\begin{eqnarray*}
\Psi:M(\Z)^{\x}\bs U/F^{\x} & \ra & I(\CO)/F^{\x} \\
  u & \mapsto & I_u,
\end{eqnarray*}
is surjective and $\la_S(\CO)$ to one.  This will be done in the following lemmas through localisation and strong approximation.

Let $p\in\Z$ be a prime and let $I(\CO_p)$ be the set of all finitely generated $\CO_p$-submodules of $F_p$ and $[I(\CO_p)]$ the set of isomorphism classes.  Let $I_p$ be a non-trivial, finitely generated $\CO_p$-submodule of $F_p$.  Then, by \cite{Reiner75}, Theorem 10.6, the submodule $I_p$ is isomorphic as an $\CO_p$-module to an ideal of $\CO_p$ so also has the property that $I_p\ox\Q_p=F_p$.  If $J_p$ is another finitely generated $\CO_p$-submodule of $F_p$, such that $J_p\cong I_p$, then the isomorphism extends to an automorphism of $F_p$ considered as a module over itself, and hence $J_p=I_p x$ for some $x\in F_p^{\x}$.  So $\card([I(\CO_p)])=|I(\CO_p)/F_p^{\x}|$.

For a prime $p$ and $u_p\in M(\Q_p)^{\x}$ let 
$$
\CO_{p,u_p}=F_p\cap u_p^{-1}M(\Z_p)u_p
$$
and recall that
$$
\CO_p=F_p\cap M(\Z_p).
$$
Let $U_p$ be the set of all $u_p\in M(\Q_p)^{\x}$ such that
$$
\CO_{p,u_p}=\CO_p.
$$

As in the global case above, $F_p^{\x}$ acts on $U_p$ from the right and $M(\Z_p)^{\x}$ acts on $U_p$ from the left.  For $u_p\in U_p$ let
$$
I_{u_p}=F_p\cap M(\Z_p)u_p.
$$
Then $I_{u_p}$ is a finitely generated $\CO_p$-module in $F_p$.

\begin{lemma}
\label{lem:LocalIdeals}
Let $p$ be a prime.  Then $|I(\CO_p)/F_p^{\x}|=1$.
\end{lemma}
\prf
Suppose $p\in S$.  Then $p$ is non-decomposed in $F$ so there exists a unique prime ideal $\p$ of $F$ above $p$, so by the isomorphism (\ref{eqn:localproduct}) we can see that $F_p\cong F_{\p}$ and hence $F_p$ is itself a field.  By \cite{Neukirch99}, Chapter II, Proposition 3.9, every ideal of $\CO_p$ is principal and hence all non-zero ideals of $\CO_p$ are isomorphic.  The claim holds for $p\in S$.

Suppose $p\notin S$.  Let $\p_1,...,\p_n$ be the prime ideals of $F$ over $p$.  Under the isomorphism (\ref{eqn:localproduct}) an ideal in $\CO_p$ gives ideals in $\CO_{\p_1},...,\CO_{\p_n}$, where $\CO_{\p_1}$ denotes the ring of integers in the field $F_{\p_i}$.  As we saw above, these ideals are each unique up to isomorphism and hence the original ideal in $\CO_p$ is also unique up to isomorphism and the claim holds for $p\notin S$.
\qed

We now prove the surjectivity of $\Psi$.

\begin{lemma}
\label{lem:3}
The map $\Psi$ is surjective.
\end{lemma}
\prf
Let $I\subset\CO$ be an ideal and for a prime $p$ let $I_p=I\ox\Z_p$.  By Lemma \ref{lem:LocalIdeals}, for all primes $p$ the $\CO_p$-module $F_p\cap M(\Z_p)$ is isomorphic to $I_p$, so there exists $u_p\in F_p^{\x}$ such that
$$
F_p\cap M(\Z_p)u_p=I_p.
$$
We have further that
\begin{eqnarray*}
F_p\cap u_p^{-1}M(\Z_p)u_p & = & u_p^{-1}(F_p\cap M(\Z_p)u_p \\
  & = & F_p\cap M(\Z_p).
\end{eqnarray*}

Let $\ut=(u_p)\in\A_{\fin}^{\x}$.  By Lemma \ref{lem:Adeles} there exist $\tilde{z}\in M(\hat{\Z})^{\x}$ and $u\in M(\Q)^{\x}$ such that $\ut=\tilde{z}u$.  Then for all primes $p$ we have that
$$
F_p\cap u^{-1}M(\Z_p)u=F_p\cap M(\Z_p)
$$
and
$$
F_p\cap M(\Z_p)u=I_p.
$$
We deduce from this that
$$
F\cap u^{-1}M(\Z)u=F\cap M(\Z)
$$
and
$$
F\cap M(\Z)u=I.
$$
This is what was required to complete the proof of the lemma.
\qed

In the following series of lemmas we prove that $\Psi$ is $\la_S(\CO)$ to one.

\begin{lemma}
\label{lem:Diag}
Let $K$ be a non-archimedian local field and $\CO_K$ its ring of integers.  For $n\in\N$ let $\D_n\subset\Mat_n(K)$ be the subspace of diagonal matrices and suppose that $u\in\GL_n(K)$ is such that
$$
\D_n\cap\Mat_n(\CO_K)\subset\D_n\cap u^{-1}\Mat_n(\CO_K)u.
$$
Then there exist $z\in\Mat_n(\CO_K)^{\x}$ and $x\in\D_n^{\x}$ such that $zu=x$.

If we suppose further that
$$
\D_n\cap\Mat_n(\CO_K)=\D_n\cap\Mat_n(\CO_K)u
$$
then we have $x\in(\D_n\cap\Mat_n(\CO_K))^{\x}$.
\end{lemma}
\prf
The first claim of the lemma will be proved by induction on $n$.  For $n=1$ we take $z=1$ and $x=u$ and the claim holds.  Assume now that the claim holds for some $n-1\in\N$.  For any ring $R$ let $\Upp_n(R)$ denote the subalgebra of upper triangular matrices in $\Mat_n(R)$.  We can choose an element $y\in\Mat_n(\CO_k)^{\x}$ such that $u'=yu\in\Upp_n(K)$ is upper triangular.  Then $u'$ has the form
$$
u'=\matrix{a}{b}{}{c},
$$
where $a\in K$, $b\in\Mat_{1\x(n-1)}(K)$ and $c\in\Upp_{n-1}(K)$.  Then
$$
(u')^{-1}=\matrix{a^{-1}}{-a^{-1}bc^{-1}}{}{c^{-1}}.
$$

Since $u$ satisfies
$$
\D_n\cap\Mat_n(\CO_K)\subset\D_n\cap u^{-1}\Mat_n(\CO_K)u
$$
it follows that $u'$ satisfies
$$
\D_n\cap\Upp_n(\CO_K)\subset\D_n\cap (u')^{-1}\Upp_n(\CO_K)u'.
$$
It follows that there exist $\al\in\CO_K$, $\be\in\Mat_{1\x(n-1)}(\CO_K)$ and $\ga\in\Upp_{n-1}(\CO_K)$ such that
$$
\matrix{\al}{a^{-1}(\al b+\be c-bc^{-1}\ga c)}{}{c^{-1}\ga c}=(u')^{-1}\matrix{\al}{\be}{}{\ga}u'=\matrix{1}{}{}{0},
$$
where the matrix on the right is the $n$ by $n$ matrix with a one in the top left corner and all other entries zero.  It then follows that $\al=1$, $c^{-1}\ga c=0$ and $b=-\be c$.  Let
$$
w=\matrix{1}{\be}{}{I_{n-1}}\in\Mat_n(\CO_K)^{\x}.
$$
Then
$$
wu'=\matrix{1}{\be}{}{I_{n-1}}\matrix{a}{-\be c}{}{c}=\matrixtwo{a}{c}.
$$
so we have that
$$
wyu=\matrixtwo{a}{c}.
$$
It is then clear that $c$ satisfies
$$
\D_{n-1}\cap\Mat_{n-1}(\CO_K)\subset\D_{n-1}\cap c^{-1}\Mat_{n-1}(\CO_K)c
$$
so by the inductive hypothesis there exist $z'\in\Mat_{n-1}(\CO_K)^{\x}$ and $x'\in\D_{n-1}^{\x}$ such that $z'c=x'$.  Setting
$$
z=\matrixtwo{1}{z'}wy\in\Mat_n(\CO_K)^{\x}\ \ \ \textrm{and}\ \ \ x=\matrixtwo{a}{x'}\in\D_n^{\x}
$$
we get that $zu=x$ and the proof of the first claim is complete.

If we also assume that $u$ satisfies
$$
\D_n\cap\Mat_n(\CO_K)=\D_n\cap\Mat_n(\CO_K)u
$$
then since $z\in\Mat_n(\CO_K)^{\x}$ and $zu=x\in\D_n^{\x}$ it follows that
$$
\D_n\cap\Mat_n(\CO_K)=\D_n\cap\Mat_n(\CO_K)zu=(\D_n\cap\Mat_n(\CO_K))x.
$$
This can be the case only if $x\in(\D_n\cap\Mat_n(\CO_K))^{\x}$ so the lemma is proved.
\qed

\begin{lemma}
\label{lem:1}
For $p\notin S$ we have $|M(\Z_p)^{\x}\bs U_p/F_p^{\x}|=1$.
\end{lemma}
\prf
Let $p\notin S$ and let $u_p\in U_p$ so that
\begin{equation}
\label{eqn:Fp1}
F_p\cap M(\Z_p)=F_p\cap u_p^{-1}M(\Z_p)u_p.
\end{equation}
By Lemma \ref{lem:LocalIdeals}, there exists $\la_p\in F_p^{\x}$ such that
$$
F_p\cap M(\Z_p)=F_p\cap M(\Z_p)u_p\la_p.
$$
Replacing $u_p$ with $u_p\la_p$ we may assume that $u_p$ satisfies
\begin{equation}
\label{eqn:Fp2}
F_p\cap M(\Z_p)=F_p\cap M(\Z_p)u_p.
\end{equation}
We shall show that there exist $z_p\in M(\Z_p)^{\x}$ and $x_p\in F_p^{\x}$ such that $z_pu_px_p$ is the identity element in $M(\Q_p)$.

Let $\ep$ be an integral element of $F$ generating $F$ over $\Q$, and let $L$ be an extension of $\Q$ containing all the zeros of the minimal polynomial of $\ep$.  Let $L_p=L\ox\Q_p$, note that neither this nor $F_p$ are necessarily fields, rather a finite direct product of fields.  The embedding $F\hra M(\Q)$ gives us an embedding $F_p\hra M(\Q_p)$.  Since $p\notin S$ we have $M(\Q_p)\cong\Mat_4(\Q_p)$ so we get an embedding
$$
F_p\ox L\hra M(\Q_p)\ox L\cong\Mat_4(L_p).
$$

Let $\CO_L$ be the ring of integers in $L$ and let $\CO_{L_p}=\CO_L\ox\Z_p$.  Let $\D$ denote the subspace of diagonal matrices in $\Mat_4(L_p)$.  It follows from our choice of $L$ that $\ep$ is diagonalisable in $M(L_p)$ by an element $A\in M(\CO_{L_p})^{\x}$.  Since $F_p=\Q_p(\ep)$, the whole of $F_p\ox L$ is simultaneously diagonalisable in $M(L_p)$, that is, we have that $A^{-1}(F_p\ox L)A$ is a subspace of $\D$, and by dimensional considerations we must have that $A^{-1}(F_p\ox L)A=\D$.  Then, tensoring (\ref{eqn:Fp1}) and (\ref{eqn:Fp2}) with $\CO_L$ and conjugating by $A$ we get
\begin{equation}
\label{eqn:Lp1}
\D\cap M(\CO_{L_p})=\D\cap\bar{u}_p^{-1}M(\CO_{L_p})\bar{u}_p
\end{equation}
and
\begin{equation}
\label{eqn:Lp2}
\D\cap M(\CO_{L_p})=\D\cap M(\CO_{L_p})\bar{u}_p,
\end{equation}
where $\bar{u}_p=u_p A$.

We are currently working in the matrix algebra $M(L_p)\cong\Mat_4(L_p)$.  Let $\p_1,...,\p_n$ be the prime ideals of $L$ above $p$.  By virtue of the isomorphism (\ref{eqn:localproduct}) we can consider separately each $\p_i$-adic component of the entries of the matrices.  It then follows from Lemma \ref{lem:Diag} and equations (\ref{eqn:Lp1}) and (\ref{eqn:Lp2}) that there exist $\bar{z}\in M(\CO_{L_p})^{\x}$ and $\bar{x}\in(\D\cap M(\CO_{L_p}))^{\x}$ such that $\bar{z}\bar{u}_p=\bar{x}$.  Setting
$$
z=A\bar{z}\in M(\CO_{L_p})^{\x}
$$
and
$$
x=A\bar{x}A^{-1}\in ((F_p\ox L)\cap M(\CO_{L_p}))^{\x}
$$
we get that $zu_p=x$.

Consider now the trace map $\tr_{L/\Q}$ of the field extension $L/\Q$.  The image of $\CO_L$ under this map is an ideal in $\Z$.  Let $\nu\in\Z$ be the greatest such that $\tr_{L/\Q}(\CO_L)\subset p^{\nu}\Z$.  We define a linear map
\begin{eqnarray*}
T:M(L_p)=M(\Q_p)\ox L & \ra &     M(\Q_p)\ox\Q\cong M(\Q_p) \\
           x\ox y     & \mapsto & x\ox p^{-\nu}(\tr_{L/\Q}\,y).
\end{eqnarray*}
Note that $T(M(\CO_{L_p}))=T(M(\Z_p)\ox\CO_L)\subset M(\Z_p)\ox\Z\cong M(\Z_p)$.

We denote by $\CO_{L_1}$ the ring of integers of the field $L_{\p_1}$.  Let $\xi\in(F\ox L_{\p_1})\cap M(\CO_{L_1})$ and let $x_1$ and $z_1$ be the $\p_1$ components of $x$ and $z$.  Then $z_1u_p=x_1$.  Now set $\al=x_1^{-1}\xi\in (F\ox L_{\p_1})\cap M(\CO_{L_1})$.  From (\ref{eqn:Fp1}) we deduce that
$$
(F\ox L_{\p_1})\cap M(\CO_{L_1})=(F\ox L_{\p_1})\cap u_p^{-1}M(\CO_{L_1})u_p
$$
so there exists $\be\in M(\CO_{L_1})$ such that $\al=u_p^{-1}\be u_p$, or equivalently $u_p\al=\be u_p$.  Setting $\ze=z_1\be\in M(\CO_{L_1})$, we then get
$$
\ze u_p=z_1\be u_p=z_1u_p\al=x_1\al=\xi.
$$
Hence we can see that $\xi u_p^{-1}\in M(\CO_{L_1})$.

We now consider the linear map defined by
\begin{eqnarray*}
T':M(\Q)\ox L_{\p_1} & \ra &     M(\Q)\ox\Q_p=M(\Q_p) \\
      x \ox y        & \mapsto & x\ox p^{-\nu}(\tr_{L_{\p_1}/\Q_p}\,y),
\end{eqnarray*}
where $\tr_{L_{\p_1}/\Q_p}$ is the trace map of the field extension $L_{\p_1}/\Q_p$.  We note that $T'$ maps $(F\ox L_{\p_1})\cap M(\CO_{L_{\p}})$ to $F_p\cap M(\Z_p)$.  The map $T'$ induces a  linear map
$$
T'':\frac{M(\CO_{L_1})}{\p_1M(\CO_{L_1})}\ra\frac{M(\Z_p)}{pM(\Z_p)}.
$$
Note that $M(\CO_{L_1})/\p_1M(\CO_{L_1})$ is a vector space over the field $\CO_{L_1}/\p_1\CO_{L_1}$ and $M(\Z_p)/pM(\Z_p)$ is a vector space over the field $\Z_p/p\Z_p\cong\F_p$.  By the choice of the integer $\nu$ we can see that the map $T''$ is not the zero map, that is, its kernel is a proper subspace of $M(\CO_{L_1})/\p_1M(\CO_{L_1})$.  An element of $M(\CO_{L_1})/\p_1M(\CO_{L_1})$ not in the kernel of $T''$ corresponds to an element of $M(\CO_{L_1})$ whose image under $T'$ is in $M(\Z_p)^{\x}$.  We choose an element $\xi\in(F\ox L_{\p_1})\cap M(\CO_{L_1})$ such that $T'(\xi)\in (F_p\cap M(\Z_p))^{\x}$ and $T'(\xi u_p^{-1})\in M(\Z_p)^{\x}$.

We shall write $\xi$ also for the element of $(F_p\ox L)\cap M(\CO_{L_p})$ with $\p_1$ component equal to $\xi$ and all other components zero.  Similarly we denote by $\ze$ the element of $M(\CO_{L_p})$ with $\p_1$ component equal to $\xi u_p^{-1}$ and all other components zero.  Then $\ze u_p=\xi$.  From the definitions of the maps $T$ and $T'$ and the isomorphism (\ref{eqn:localproduct}) it follows that
$$
T(\xi)=T'(\xi')\in(F_p\cap M(\Z_p))^{\x}
$$
and
$$
T(\ze)=T'(\ze)\in M(\Z_p)^{\x}.
$$

Now set $x_p=T(\xi)\in(F_p\cap M(\Z_p))^{\x}$ and $z_p=T(\ze)\in M(\Z_p)^{\x}$.  We then have that
$$
z_pu_p=x_p
$$
and the lemma is proved.
\qed

\begin{lemma}
\label{lem:2}
For $p\in S$ we have $|M(\Z_p)^{\x}\bs U_p/F_p^{\x}|=f_p(F)$.
\end{lemma}
\prf
Let $p\in S$.  We then have that $M(\Q_p)$ is a division algebra so in particular $M(\Q_p)^{\x}=M(\Q_p)\smallsetminus\{ 0\}$.  Further, since $\Q_p$ is a local field, $M(\Z_p)$ is the unique maximal order of $M(\Q_p)$ (cf.Proposition \ref{pro:orders}(d)) and hence $u_p^{-1}M(\Z_p)u_p=M(\Z_p)$ for all $u_p\in M(\Q_p)^{\x}$.  Thus the claim is equivalent to
$$
|M(\Z_p)^{\x}\bs M(\Q_p)^{\x}/F_p^{\x}|=f_p(F).
$$
Let $v_M$ be the $p$-adic valuation on $M(\Q_p)$ and let
$$
e_M=e(M(\Q_p)/\Q_p)=|v_M(M(\Q_p)^{\x})/v_M(\Q_p^{\x})|
$$
be the ramification index of $\Q_p$ in $M(\Q_p)$.  The valuation $v_M$ is a surjective group homomorphism of $M(\Q_p)^{\x}$ onto $\Z$ with kernel $M(\Z_p)^{\x}$.  It follows that
$$
e_M=|M(\Z_p)^{\x}\bs M(\Q_p)^{\x}/\Q_p^{\x}|.
$$
By Proposition 17.7 of \cite{Pierce82} we have
$$
|M(\Z_p)^{\x}\bs M(\Q_p)^{\x}/\Q_p^{\x}|=e_M=\deg M(\Q_p)=4.
$$

Let $e_p(F)$ be the ramification index of $p$ in $F$.  Then by \cite{Neukirch99}, Chapter I, Proposition 8.2, and since $p$ does not split in $F$, we have
\begin{equation}
\label{eqn:fundid}
e_p(F)f_p(F)=[F:\Q]=4.
\end{equation}
Let $\CCO$ be the valuation ring of $F_p$ and let $\Pi=\pi\CCO$ be the unique maximal ideal of $\CCO$, where $\pi\in\CCO$ is a generator of the principal ideal $\Pi$.  Let $\p=p\Z_p$ be the unique maximal ideal of $\Z_p$.  Then $\p=\Pi\cap\Z_p$ and $p=\pi^{e_p(F)}$.  Thus $v_M(\pi)=v_M(p)/e_p(F)$, which implies that
$$
|M(\Z_p)^{\x}\bs M(\Q_p)^{\x}/F_p^{\x}|=\frac{|M(\Z_p)^{\x}\bs M(\Q_p)^{\x}/\Q_p^{\x}|}{e_p(F)}=\frac{4}{e_p(F)},
$$
and hence from (\ref{eqn:fundid}) we get
$$
|M(\Z_p)^{\x}\bs M(\Q_p)^{\x}/F_p^{\x}|=f_p(F),
$$
as required.
\qed

\begin{lemma}
$\Psi$ is $\la_S(\CO)$ to one.
\end{lemma}
\prf
Recall that $\CO\in O_F(S)$ for some $F\in C(S)$.  Let $\hat{F}=\prod'_p F_p$ be the restricted product over all primes $p$, that is, all but finitely many components are integral, and let $\hat{U}$ be the set of all $\ut=(u_p)\in M(\A_{\fin})^{\x}$ such that $\CO_p=F_p\cap u_p^{-1}M(\Z_p)u_p=\CO_{p,u_p}$ for all $p$.

For $I\in I(\CO)$ let $U_I=\Psi^{-1}(IF^{\x})$.  Let $u\in U_I$ for some $I\in I(\CO)$.  We denote by $\hat{u}$ the element of $M(\A_{\fin})$ each of whose components is $u$.  Setting $u_p=u$ for all primes $p$, we see that for all $p$
$$
\CO_p=\CO\ox\Z_p=\CO_u\ox\Z_p=\CO_{p,u_p}
$$
so $\hat{u}\in\hat{U}$.

For each $I\in I(\CO)$ we define the map
$$
\Phi_I:M(\Z)^{\x}\bs U_I/F^{\x}\ra M(\hat{\Z})^{\x}\bs\hat{U}/\hat{F}^{\x}
$$
by setting $\Phi_I(M(\Z)^{\x}uF^{\x})=M(\hat{\Z})^{\x}\hat{u}\hat{F}^{\x}$.  We shall show that $\Phi_I$ is a bijection.  Since from Lemmas \ref{lem:1} and \ref{lem:2} we know that $|M(\hat{\Z})^{\x}\bs\hat{U}/\hat{F}^{\x}|=\la_S(\CO)$, the lemma will then follow.

First we show surjectivity.  Let $\ut=(u_p)\in\hat{U}$, so that $\CO_{p,u_p}=\CO_p$ for all primes $p$.  For all primes $p$, by Lemma \ref{lem:LocalIdeals} we have $I\ox\Q_p\cong F_p\cap M(\Z_p)u_p$ so there exists $x_p\in F_p^{\x}$ such that $I\ox\Q_p=F_p\cap M(\Z_p)u_px_p$.  We replace $u_p$ with $u_px_p$.  By Lemma \ref{lem:Adeles} there exist $\tilde{z}=(z_p)\in M(\hat{\Z})^{\x}$ and $u\in M(\Q)^{\x}$ such that $\ut=\tilde{z}\hat{u}$.  It follows that $\CO_{p,u}=\CO_p$ for all primes $p$ and hence we get that $\CO_u=\CO$, so that $u\in U$.  Furthermore, since
$$
I\ox\Q_p=F_p\cap M(\Z_p)u_p=F_p\cap M(\Z_p)u
$$
for all primes $p$, it follows that $I=F\cap M(\Z)u$, so $u\in U_I$.  We have proven that $\Phi_I$ is surjective.

For the proof of injectivity, suppose that $u,v\in U_I$ are such that $\Phi_I(u)=\Phi_I(v)$.  The assumption that $u,v\in U_I$ means that there is an isomorphism of $\CO$-modules $F\cap M(\Z)u\cong F\cap M(\Z)v$.  As we saw above this implies that there exists $x\in F^{\x}$ such that $F\cap M(\Z)u=F\cap M(\Z)vx$.  Replacing $v$ with $vx$ it follows that
\begin{equation}
\label{eqn:IdealEq}
\hat{F}\cap M(\hat{\Z})\hat{u}=\hat{F}\cap M(\hat{\Z})\hat{v}.
\end{equation}
The assumption that $\Phi_I(u)=\Phi_I(v)$ means that there exist $\tilde{z}\in M(\hat{\Z})^{\x}$ and $\tilde{x}\in\hat{F}^{\x}$ such that $\hat{u}=\tilde{z}\hat{v}\tilde{x}$.  Together with (\ref{eqn:IdealEq}) this implies that
$$
\hat{F}\cap M(\hat{\Z})\hat{v}=\hat{F}\cap M(\hat{\Z})\hat{v}\tilde{x}
$$
and hence $\tilde{x}\in\hat{v}^{-1}M(\hat{\Z})^{\x}\hat{v}$.  In other words, there exists $\tilde{y}\in M(\hat{\Z})^{\x}$ such that $\hat{v}\tilde{x}=\tilde{y}\hat{v}$.  We deduce that
$$
\hat{u}=\tilde{z}\hat{v}\tilde{x}=(\tilde{z}\tilde{y})\hat{v}.
$$
Replacing $\tilde{z}$ with $\tilde{z}\tilde{y}\in M(\hat{\Z})^{\x}$, we now have that
$$
\tilde{z}=\hat{u}\hat{v}^{-1}\in M(\hat{\Z})^{\x}\cap M(\Q)^{\x}=M(\Z)^{\x}.
$$
Writing $z=\tilde{z}\in M(\Z)^{\x}$ we get that $u=zv$ and injectivity is proven.
\qed

  \chapter{Comparing Geodesics and Orders}
    \label{ch:Compare}
    \markright{\textnormal{\thechapter{. Comparing Geodesics and Orders}}}
    In this final chapter we complete the proof of the Main Theorem.  We do this by using our results on orders in subfields of division algebras of degree four to interpret the Prime Geodesic Theorem algebraically.  This will then yield information about the class numbers of those orders.  All notation is as defined in previous chapters.  In particular we have $G=\SL_4(\R)$ and take $\Ga=\CG(\Z)$ as in Chapter \ref{ch:DivAlg}.

\section{Regular geodesics}
Let $F$ be a number field with $r$ real embeddings and $s$ pairs of complex conjugate embeddings and let $t=r+s-1$.  Let $\CO_F$ denote the ring of integers of $F$ and $\CO_F^{\x}$ the unit group therein.  By Dirichlet's unit theorem (\cite{Neukirch99}, I Theorem 7.4), there exist units $\ep_1,...,\ep_t\in\CO_F^{\x}$ such that every $\ep\in\CO_F^{\x}$ can be written uniquely as a product
$$
\ep=\zeta\ep_1^{\nu_1}\cdots\ep_t^{\nu_t},
$$
where $\zeta$ is a root of unity and $\nu_i\in\Z$.  The set $\{\ep_1,...,\ep_t\}$ is not uniquely determined.  We shall refer to such a set of units as a \emph{system of fundamental units}\index{fundamental unit}.  If $\CO\subset\CO_F$ is an order of $F$ and $\{\ep_1,...,\ep_t\}$ a system of fundamental units then there exist $\mu_i\in\N$ and units $\bar{\ep}_i=\ep_i^{\mu_i}$ such that every $\ep\in\CO^{\x}$ can be written uniquely as a product
$$
\ep=\zeta\bar{\ep}_1^{\nu_1}\cdots\bar{\ep}_t^{\nu_t}
$$
where $\zeta$ is a root of unity and $\nu_i\in\Z$.  Once again, the set $\{\bar{\ep}_1,...,\bar{\ep}_t\}$ is not uniquely determined.  We shall call such a set of units \emph{a system of fundamental units of $\CO$}.

If $\tau_1,...,\tau_{t+1}$ are distinct embeddings of $F$ into $\C$ which are pairwise non-conjugate, then the \emph{regulator} \index{regulator} $R(\CO)$ of an order $\CO\subset\CO_F$ is defined as the absolute value of the determinant of an arbitrary minor of rank $t$ of the matrix
$$
\left(\log|\tau_i(\bar{\ep}_j)|\right)_{ij},
$$
where $1\leq i\leq t+1$ and $1\leq j\leq t$.  The value of the regulator does not depend on the choice of a system of fundamental units.

We shall sometimes refer, somewhat imprecisely, to the elements of a system of fundamental units simply as \emph{fundamental units}.  In the case $t=1$, let $\ep_1$ be a fundamental unit of an order $\CO$ in $F$.  We shall refer to the elements of the set $\{\zeta\ep_1^{\pm 1}:\zeta\textrm{ root of unity in }F\}$ as \emph{the fundamental units of $\CO$}, or in the case that $\CO=\CO_F$ \emph{the fundamental units of $F$}.

In the case that $F$ is a totally complex quartic field we have $t=1$, so we can choose a fundamental unit $\ep_1$ such that for all $\ep\in\CO_F^{\x}$ we have $\ep=\zeta\ep_1^{\nu_1}$, where $\zeta$ is a root of unity and $\nu_1\in\Z$.  Then $\ep_1$ is determined up to inversion and multiplication by a root of unity.  In this case the regulator of an order $\CO\subset\CO_F$ of $F$ with fundamental unit $\bar{\ep}_1$ satisfies
$$
R(\CO)=\left|\,\log|\tau(\bar{\ep}_1)|\,\right|,
$$
where $\tau$ is any embedding of $F$ into $\C$.

We consider subfields of the field $F$.  Apart from $\Q$ these are all quadratic.  By Dirichlet's unit theorem, real quadratic fields contain a single fundamental unit (up to inversion and multiplication by $\pm 1$) and complex quadratic fields contain no fundamental unit.  We can see then that $F$ can have at most one real subfield.  Let $C^r(S)$ \index{$C^r(S)$} and $C^c(S)$ \index{$C^c(S)$} be respectively the subsets of $C(S)$ consisting of fields with and without a real quadratic subfield.  Let $O^r(S)$ \index{$O^r(S)$} and $O^c(S)$ \index{$O^c(S)$} be respectively the subsets of $O(S)$ consisting of orders in fields in $C^r(S)$ and $C^c(S)$.  Then $C(S)=C^r(S)\cup C^c(S)$ and $O(S)=O^r(S)\cup O^c(S)$.

Let $[\ga]\in\CE_P^{p,\reg}(\Ga)$.  Then by Lemma \ref{lem:centralisers} the centraliser $M(\Q)_{\ga}$ is a totally complex quartic field which embeds into $M(\Q)$ and we denote this subfield of $M(\Q)$ by $F_{\ga}$.  We say that an element $\ga\in\Ga$ is weakly neat if the adjoint $\Ad(\ga)$ has no non-trivial roots of unity as eigenvalues.  Let $\CE_P^{p,\rw}(\Ga)$ \index{$\CE_P^{p,\rw}(\Ga)$} be the set of primitive, regular, weakly neat elements in $\CE_P(\Ga)$.

\begin{lemma}
\label{lem:WeaklyNeat}
Let $[\ga]\in\CE_P^{p,\reg}(\Ga)$.  Then $F_{\ga}\in C^c(S)$ if and only if $[\ga]\in\CE_P^{p,\rw}(\Ga)$.
\end{lemma}
\prf
Let $\CO$ be the maximal order in $F$ and let $\CO^{\x}$ denote the group of units.  By Dirichlet's unit theorem we can see that $F_{\ga}$ contains a real quadratic subfield if and only if $\CO^{\x}\cap\R\neq\{\pm 1\}$.  Let $\ga_1\in\CO^{\x}$ be a fundamental unit.  The element $\ga$ is a unit $F_{\ga}$ so there exists a root of unity $\ze$ and $n\in\Z$ such that $\ga=\ze\ga_1^n$.  If there also exist a root of unity $\xi$ and $m\in\Z\smallsetminus\{0\}$ such that $\xi\ga_1^m\in\R$ then $\ga_1$ is not weakly neat and so neither is $\ga$.

Conversely, suppose $\ga$ is not weakly neat.  We know that $\ga$ is conjugate in $G$ to an element $a_{\ga}b_{\ga}\in A^-B$.  From the assumption that $\ga$ is not weakly neat it follows that there exists $m\in\N$ such that one of the components of $b_{\ga}^m$ is diagonal.  Since the group $AB$ is abelian we have that $\ga^m$ is conjugate in $G$ to $a_{\ga}^m b_{\ga}^m$.  Since $\ga^m\in\CE_P(\Ga)$ it follows from Lemma \ref{lem:centralisers} that $b_{\ga}^m$ is diagonal and hence $\ga^m$ generates a real quadratic subfield of $F_{\ga}$.  The lemma follows.
\qed

Define the map
$$
\th:\CE_P^{p,\rw}(\Ga)\ra O^c(S)
$$
by
$$
\th:[\ga]\mapsto F_{\ga}\cap M(\Z).
$$
Let $\ga,\d\in\Ga$.  Then $F_{\d\ga\d^{-1}}=\d F_{\ga}\d^{-1}$, and since by Lemma \ref{lem:unit} the group $\Ga$ is a subgroup of the multiplicative group $M(\Z)^{\x}$, we can conclude that
$$
F_{\d\ga\d^{-1}}\cap M(\Z)=\d\left(F_{\ga}\cap M(\Z)\right)\d^{-1}\cong F_{\ga}\cap M(\Z).
$$
Hence the map is well defined.

\begin{proposition}
The map $\th$ is surjective.
\end{proposition}

\prf
Let $\CO\in O^c(S)$ be an order in the field $F\in C^c(S)$.  Since $\deg F=4$, Proposition \ref{pro:subfields}(b) implies that $F$ is a maximal subfield of $M(\Q)$.

Lemma \ref{lem:orderembedding} tells us that there exists an embedding $\si:F\ra M(\Q)$ such that
$$
\CO=\CO_{\si}=\si^{-1}\left(\si(F)\cap M(\Z)\right).
$$
Let $\ep$ be a fundamental unit in $\CO^{\x}$.  Let $\ga$ be the image of $\ep$ under the embedding $\si$.  The unit $\ep$ is integral in $F$, so $\ga\in M(\Z)$.

First we show that $\ga\sim_G a_{\ga}b_{\ga}$ for some $a_{\ga}\in A^{-}$, $b_{\ga}\in B$.  The field $F$ has two pairs of conjugate embeddings into $\C$.  Let $\si_1$, $\si_2$ be two distinct, non-conjugate embeddings of $F$ into $\C$.  The $\R$-algebra $F\ox_{\Q}\R$ is isomorphic to the commutative $\R$-algebra $\C\oplus\C$ via the map 
$$
\al\ox x \mapsto \left(x\si_1(\al),x\si_2(\al)\right).
$$
We note that $\CO\ox_{\Z}\R =F\ox_{\Q}\R$ and we get the series of inclusions
$$
\CO \subset \CO\ox_{\Z}\R = F\ox_{\Q}\R \subset M(\Q)\ox_{\Q}\R = \Mat_4(\R) \subset \Mat_4(\C).
$$
Thus we see that $\ga$ may be considered as an element, which we will denote $\bar{\ga}$, of a commutative subalgebra of $\Mat_4(\C)$ isomorphic to $\C\oplus\C$.  The matrix $\bar{\ga}$ has real entries and can be diagonalised in $\Mat_4(\C)$ to the matrix
$$
X=\matrixfour{\si_1(\ga)}{\overline{\si_1(\ga)}}{\si_2(\ga)}{\overline{\si_2(\ga)}},
$$
since the eigenvalues of $\bar{\ga}$ are the roots of its minimal polynomial, ie. the values of $\ga$ under its different embeddings into $\C$.  Let $a,b\in\R^+$; $\th,\phi\in [-\pi,\pi]$ be such that $\si_1(\ga)=ae^{i\th}$ and $\si_2(\ga)=be^{i\phi}$.  Since $\ga\in M(\Z)$, we have $\prod_{\si}\si(\ga)=N_{F|\Q}(\ga)=\det\ga=\pm1$, where $N_{F|\Q}$ is the field norm and the product is over all embeddings $\si$ of $F$ into $\C$ (\cite{Neukirch99}, I Proposition 2.6(iii)).  Hence $a^2 b^2=\pm 1$.  But $a$ and $b$ are both real so we must have $a^2 b^2=1$, and so $b=a^{-1}$ and we have
$$
X=\matrixfour{ae^{i\th}}{ae^{-i\th}}{\rez{a}e^{i\phi}}{\rez{a}e^{-i\phi}}.
$$
Without loss of generality we assume $a<b$ so that $a\in(0,1)$.
Let $R$ be the $2\x 2$ matrix
$$
\rez{\sqrt{2}}\matrix{i}{-1}{1}{-i}
$$
and $Y$ the $4\x 4$ matrix
$$
\matrixtwo{R}{R}\in\Mat_4(\C)
$$
Then $X$ is conjugate in $\Mat_4(\C)$ via $Y$ to the real matrix
$$
X'=\matrixtwo{aR(\th)}{\rez{a}R(\phi)}\in A^-B.
$$
So $\bar{\ga}$ is conjugate in $\Mat_4(\C)$ to an element of $A^-B$, hence by \cite{Lang02}, XIV Corollary 2.3, the matrix $\bar{\ga}$ is conjugate in $\Mat_4(\R)$ to $X'$, that is $(Z')^{-1}\bar{\ga}Z'=X'$ for some $Z'\in\Mat_4(\R)$.  Since $Z'$ is invertible we have that $n=\det\,Z'\neq 0$.  Let $Z=|n|^{-\frac{1}{4}}Z'$.  Then $\det\,Z=\pm 1$ and $Z^{-1}\bar{\ga}Z=X'$, so that $\ga$ is conjugate in $\SL_4^{\pm}(\R)$ to an element of $A^-B$.  Suppose that $\det\,Z=-1$.  Let $W$ be the matrix
$$
W=\left( \begin{array}{cccc} 
	       0 & 1 & \ & \ \\ 
	       1 & 0 & \ & \ \\  
	       \ & \ & 1 & 0 \\  
	       \ & \ & 0 & 1 \end{array}\right).
$$
Then $ZW\in G$ and $(ZW)^{-1}\bar{\ga}ZW=X''$, where
$$
X''=\matrixtwo{aR(-\th)}{\rez{a}R(\phi)}\in A^-B.
$$
Hence $\ga$ is conjugate in $G$ to an element of $A^-B$.

We also saw above that $\det\ga=1$ and so $\ga\in\Ga$ and hence $\ga\in\CE_P(\Ga)$.  It then follows from Lemma \ref{lem:centralisers} that $\ga$, and hence also $\ep$, generates either a real quadratic or a totally complex quartic field over $\Q$.  The field generated by $\ep$ over $\Q$ is a subfield of $F$.  However $F\in C^c(S)$ so has no real quadratic subfields.  Hence $\Q(\ga)\cong\Q(\ep)$ is a totally complex quartic field and Lemma \ref{lem:centralisers} tells us that $\ga$ is regular.  It then follows from Lemma \ref{lem:WeaklyNeat} that $\ga$ is weakly neat.

It remains to show that $\ga$ is primitive in $\Ga$.  Note that $\Q(\ep)$ is a quartic field contained in and hence equal to the field $F$.  Recall that we write $F_{\ga}=M(\Q)_{\ga}$ for the centraliser of $\ga$ in $M(\Q)$.  Then we have that $F_{\ga}=\Q(\ga)\cong F$ and $\CO^{\x}\cong F_{\ga}\cap M(\Z)^{\x}$.  The fundamental unit $\ep$ is primitive in $\CO^{\x}$ and so $\ga=\si(\ep)$ is primitive in $F_{\ga}\cap M(\Z)^{\x}$.  Suppose $\ga$ is not primitive in $\Ga$.  Then there exists $\ga_0\in\Ga\subset M(\Z)^{\x}$ such that $\ga=\ga_0^{\kappa}$ for some $\kappa\in\N$.  However we also have $\ga_0\in F_{\ga}$ and so $\ga_0\in F_{\ga}\cap M(\Z)^{\x}$, which is a contradiction.  Hence $\ga$ is primitive in $\Ga$ and we have $\ga\in\CE_P^{p,\reg}(\Ga)$.

This concludes the proof of the surjectivity of $\th$.
\qed

\begin{lemma}
\label{lem:rootsOfUnity}
Let $\ze$ be a torsion element of $M(\Z)$.  Then $\det(\ze)=1$.
\end{lemma}
\prf
By Proposition \ref{pro:subfields}(b), the element $\ze$ generates a subfield $\Q(\ze)$ of $M(\Q)$ of degree 1, 2 or 4 over $\Q$.  If $[\Q(\ze):\Q]=1$ or 2 then by Lemma \ref{lem:unit} we have $\det(\ze)=1$.

For a primitive $n^{\rm th}$ root of unity $\xi$, the degree $[\Q(\xi):\Q]$ is equal to $\varphi(n)=\#\{1\leq m\leq n:(m,n)=1\}$.  From \cite{Chandrasekharan68}, Chapter II, \S 2 we see that
$$
\varphi(n)=\prod_{i=1}^r p_i^{a_i-1}(p_i-1),
$$
where $n=p_1^{a_1}...p_r^{a_r}$ is the prime decomposition of $n$.  It is then straightfoward to show that $[\Q(\xi):\Q]=4$ if and only if $n=5,8,10$ or 12.  In these cases the minimal polynomial of $\xi$ is respectively $x^4+x^3+x^2+x+1$, $x^4+1$, $x^4-x^3+x^2-x+1$ or $x^4-x^2+1$.  In each of these cases we read off from the constant term of the minimal polynomial that $N_{\Q(\xi)|\Q}(\xi)=1$.  The lemma follows.
\qed

For $\CO\in O(S)$ let $\mu_{\CO}$ be the number of roots of unity in the order $\CO$.  For $[\ga]\in\CE_P^{p,\reg}(\Ga)$ let $\mu_{\ga}$ denote the number of roots of unity in the order $F_{\ga}\cap M(\Z)$.  Note that $F_{\ga}\cap M(\Z)$ is equal to the centraliser $M(\Z)_{\ga}$ of $\ga$ in $M(\Z)$, and by Lemma \ref{lem:rootsOfUnity} the roots of unity in $M(\Z)_{\ga}$ are all in $\Ga$ and hence in $\Ga_{\ga}$.  It follows that $\mu_{\ga}$ is equal to the cardinality of the torsion part of $\Ga_{\ga}$.  It is also immediate from the definitions that, writing $\CO_{\ga}=F_{\ga}\cap M(\Z)$, we have $\mu_{\ga}=\mu_{\CO_{\ga}}$.

\begin{lemma}
\label{lem:SplitFactors}
Let $\CO\in O(S)$ be an order in the field $F\in C(S)$ and let $f(x)$ be the minimal polynomial over $\Q$ of a fundamental unit in $\CO$.  Then the number of roots of $f(x)$ in $F$ which are also fundamental units in $\CO$ is independent of the choice of fundamental unit.  We call this number $\ka(\CO)$\index{$\ka(\CO)$}.  Note that $\ka(\CO)$ is equal to 1,2 or 4.
\end{lemma}

\prf
Let $\ep$ be a fundamental unit in $\CO$ with minimal polynomial $f(x)$ over $\Q$.  Suppose there exists another fundamental unit $\d\in\CO^{\x}$ with $\d\neq\ep$ which is also a root of $f(x)$.  Then there exist embeddings $\al$ and $\be$ of $F$ into $\C$ such that $\al(\ep)=\be(\d)$.

Let $\eta\in\CO^{\x}$ also be a fundamental unit in $\CO$ with minimal polynomial $g(x)$ over $\Q$.  Then $\eta=\ze\ep^{\pm 1}$ for some root of unity $\ze\in\CO^{\x}$.  We have
$$
\al(\eta)=\al(\ze\ep^{\pm 1})=\al(\ze)\al(\ep)^{\pm 1}=\al(\ze)\be(\d)^{\pm 1}.
$$
Let $n$ be the order of the root of unity $\ze$.  The order $\CO$ is a ring so contains all roots of unity of order $n$ and hence there exists an $n^{th}$ root of unity $\xi\in\CO^{\x}$ such that $\be(\xi)=\al(\ze)$.  Let $\th$ be the fundamental unit $\xi\d^{\pm 1}\in\CO^{\x}$.  Then $\al(\eta)=\be(\th)$ and so $\th$ is also a root of $g(x)$.  The lemma follows.
\qed

\begin{proposition}
\label{pro:ThetahlaToOne}
The map $\th$ is $\frac{4h(\CO)\la_S(\CO)\mu_{\CO}}{\ka(\CO)}$ to one.
\end{proposition}
\prf
The question is, given $\CO\in O(S)$, how many conjugacy classes $[\ga]\in\CE_P^{p,\reg}(\Ga)$ are there such that $F_{\ga}\cap M(\Z)\cong\CO$.  Suppose $\CO$ is an order in the field $F$.  We saw in the proof of the previous proposition that given an embedding $\si:F\hra M(\Q)$ such that $\CO_{\si}=\CO$ and a fundamental unit $\ep\in\CO^{\x}$ the element $\ga=\si(\ep)$ satisfies both $F_{\ga}\cap M(\Z)\cong\CO$ and $[\ga]\in\CE_P^{p,\reg}(\Ga)$.  Conversely, it is clear that for $[\ga]\in\CE_P^{p,\reg}(\Ga)$, the element $\ga$ is a fundamental unit in the order $F_{\ga}\cap M(\Z)$.

Recall that $M(\Z)^{\x}$ acts on the set $\Si(\CO)$ from the left and that, by Proposition \ref{pro:embeddings}, the cardinality of the set $M(\Z)^{\x}\bs\Si(\CO)$ is $h(\CO)\la_S(\CO)$.  If $\ep,\d\in\CO^{\x}$ are fundamental units and $\si,\tau$ are embeddings of $F$ into $M(\Q)$ such that $\CO_{\si}=\CO_{\tau}=\CO$, then $\si(\ep)$ and $\tau(\d)$ are in the same $M(\Z)^{\x}$-conjugacy class if and only if there exists an automorphism $\al$ of $F$ such that $\al(\d)=\ep$ and the embeddings $\si$ and $\tau\circ\al$ are in the same class modulo $M(\Z)^{\x}$.

We will show that for each $\ga\in\CE_P^{p,\reg}(\Ga)$ the $M(\Z)^{\x}$-conjugacy class of $\ga$ decomposes into two $\Ga$-conjugacy classes.  If $\ga,\d\in\Ga$ are in the same $M(\Z)^{\x}$-conjugacy class we shall write $[[\ga]]=[[\d]]$.  For $\ga\in\Ga$, if $\ga$ is central in $M(\Q)$ then clearly the conjugacy class $[[\ga]]=\{\ga\}$.  Conversely, if $[[\ga]]=\{\ga\}$ then $\ga$ commutes with every element of $M(\Z)$ and, since $M(\Z)$ is an order in $M(\Q)$, it follows that $\ga$ commutes with every element of $M(\Q)$.  Hence the conjugacy class $[[\ga]]$ consists of a single element if and only if $\ga$ is central in $M(\Q)$.  Note that, by \cite{Reiner75}, Theorem 34.9, the image of $M(\Z)$ under the map $\det\!$ is equal to $\Z$, so there exist elements of $M(\Z)$ whose image under $\det\!$ is $-1$.

Let $\ga\in\CE_P^{p,\reg}(\Ga)$.  Then $\ga$ is not central in $M(\Q)$ so there exists $\d\neq\ga$ in $[[\ga]]$ such that $\ga=\eta^{-1}\d\eta$ for some $\eta\in M(\Z)^{\x}$ with $\det(\eta)=-1$.  Suppose further that $[\ga]=[\d]$, so that $\ga=\th^{-1}\d\th$ for some $\th\in\Ga$.  Then $\ga=\eta^{-1}\th\ga\th^{-1}\eta$, so $\xi=\eta^{-1}\th$ is an element of $F_\ga$ with $\det(\xi)=-1$.  By Lemma \ref{lem:unit} we then have $\xi\in F_{\ga}\cap M(\Z)^{\x}$, where $\CO=F_{\ga}\cap M(\Z)$ is an order in the totally complex quartic field $F_{\ga}$ and $\CO^{\x}=F_{\ga}\cap M(\Z)^{\x}$ is the group of units.  But since $F_{\ga}$ is totally complex quartic, Lemma \ref{lem:unit} tells us that $\det(\xi)=1$ and hence, by contradiction, $[\ga]\neq[\d]$.  This shows that each $M(\Z)^{\x}$-conjugacy class of elements in $\CE_P^{p,\reg}(\Ga)$ decomposes into at least two $\Ga$-conjugacy classes.

Suppose that $[[\ga]]=[[\d]]$ but $[\ga]\neq[\d]$.  Then there exists $\eta\in M(\Z)^{\x}$ such that $\ga=\eta^{-1}\d\eta$ and $\det(\eta)=-1$.  Let $\be\in\Ga$ also be such that $[[\ga]]=[[\be]]$ but $[\ga]\neq[\be]$.  Then there exists $\th\in M(\Z)^{\x}$ such that $\ga=\th^{-1}\be\th$ and $\det(\th)=-1$.  Then $\d=\eta\th^{-1}\be\th\eta^{-1}$, where $\det(\th\eta^{-1})=1$ so $[\d]=[\be]$.  Hence each $M(\Z)^{\x}$-conjugacy class of elements in $\CE_P^{p,\reg}(\Ga)$ decomposes into precisely two $\Ga$-conjugacy classes.

We conclude that the choice of a fundamental unit $\ep\in\CO^{\x}$ gives us $2h(\CO)\la_S(\CO)$ $\Ga$-conjugacy classes $[\ga]\in\CE_P^{p,\reg}(\Ga)$ such that $F_{\ga}\cap M(\Z)\cong\CO$.  We have seen that the choice of fundamental unit is determined up to inversion and multiplication by a root of unity so there are $2\mu_{\CO}$ choices for $\ep$.  We saw in the previous proposition that for a given embedding $\si$ of $F$ into $M(\Q)$ and choice of a fundamental unit $\ep\in\CO$ there is a conjugacy class $[\ga_{\si,\ep}]\in\CE_P^{p,\reg}(\Ga)$ with $\si(\ep)=\ga_{\si,\ep}$.  However, we saw above that if there exists an automorphism $\al$ of $F$ such that $\d=\al(\ep)$, then setting $\tau=\si\circ\al^{-1}$ we have $[\ga_{\si,\ep}]=[\ga_{\tau,\d}]$.  Such an automorphism exists if and only if $\d$ is a root of the minimal polynomial of $\ep$ over $\Q$.  The claim of the proposition folows.
\qed

\begin{lemma}
\label{lem:KaCondition}
Let $\ga\in\CE_P^{p,\rw}(\Ga)$ so that $\ga$ is conjugate in $G$ to $a_{\ga}b_{\ga}\in A^-B$ where
$$
a_{\ga}=\matrixfour{a}{a}{a^{-1}}{a^{-1}}\ \ \textrm{ and }\ \ b_{\ga}=\matrixtwo{R(\th)}{R(\phi)}
$$
for some $a\in(0,1)$, $\th,\phi\in\R$.  If $\ka(\CO_{\ga})>1$ then either $\th+\phi$ or $\th-\phi$ is in $\frac{\pi}{2}\Z\cup\frac{\pi}{3}\Z$.
\end{lemma}
\prf
The element $\ga$ is a fundamental unit in $\CO_{\ga}$.  We shall write $f(x)$ for its minimal polynomial.  The roots of $f(x)$ are the eigenvalues of $\ga$, which are $\al_1=ae^{i\th}$, $\al_2=ae^{-i\th}$, $\al_3=a^{-1}e^{i\phi}$, $\al_4=a^{-1}e^{-i\phi}$.

Suppose $\ka(\CO_{\ga})>1$.  Then there exist $i,j\in\{1,2,3,4\}$, $i\neq j$ and a root of unity $\ze$ such that $\al_i=\ze\al_j^{\pm 1}$.  If $\al_i=\ze\al_j$ then $\al_i\al_j^{-1}=\ze$ and either $\{i,j\}=\{1,2\}$ or $\{i,j\}=\{3,4\}$.  Without loss of generality take $i=1$ and $j=2$, then $\ze=\al_1\al_2^{-1}=e^{i2\th}$.  It follows that if $n\in\N$ is such that $\ze^n=1$ then $\al_1^{2n}=a^n\in(0,1)$ and so $F_{\ga}$ has a real quadratic subfield.  However, by Lemma \ref{lem:WeaklyNeat} this contradicts the assumption that $\ga\in\CE_P^{p,\rw}(\Ga)$, so we must have $\al_i=\ze\al_j^{-1}$.

By considering the possible values of $i$ and $j$ which would satisfy $\al_i=\ze\al_j^{-1}$ we see that either $\th+\phi$ or $\th+\phi$ is equal to $q\pi$, where $q\in\Q$ is such that $\ze=e^{iq\pi}$.  It remains to show what the possible values for $q$ are, that is, which roots of unity may occur in $\CO_{\ga}$.  We saw in the proof of Lemma \ref{lem:rootsOfUnity} that the only roots of unity which may occur in a quartic field are $\pm 1$ and roots of order $3,4,5,6,8,10$ or $12$.  If $\ze$ is a root of order $5,8,10$ or $12$ then $\ze$ generates a totally complex quartic field, however, $\Q(\ze)$ has the real quadratic subfield $\Q(\ze+\ze^{-1})$, which possibility is excluded in our case.  So we are left with the possibility that $\ze=\pm 1$ or $\ze$ is a root of order $3,4$ or $6$, which implies the claim of the lemma.
\qed

\begin{lemma}
Let $[\ga]\in\CE_P^{p,\reg}(\Ga)$.  Then
$$
\chi_1(\Ga_{\ga})=\rez{\mu_{\ga}}.
$$
\end{lemma}
\prf
By Theorem \ref{thm:ECharPos},
$$
\chi_1(\Ga_{\ga})=\frac{\left[\Ga_{\ga,A}:\Ga_{\ga,A}'\right]}{\left[\Ga_{\ga}:\Ga_{\ga}'\right]},
$$
where $\Ga'\subset\Ga$ is a torsion free subgroup of finite index.  In the case under consideration we have that $\Ga_{\ga}$ is isomorphic as a group to the group of norm one elements in the order $\CO_{\ga}=F_{\ga}\cap M(\Z)$.  By Lemma \ref{lem:unit} this is equal to the group of units $\CO_{\ga}^{\x}$.  By Dirichlet's unit theorem and Lemma \ref{lem:rootsOfUnity} it follows that
$$
\Ga_{\ga}\cong\ep^{\Z}\x\mu(\ga),
$$
for some generator $\ep$ and where $\mu(\ga)$ is the finite cyclic group of roots of unity in $\CO_{\ga}^{\x}$.  It then follows that
$$
\Ga_{\ga}'\cong\ep^{k\Z},
$$
for some $k\in\N$.  Since every torsion element in $\Ga_{\ga}$ is sent to the identity under the projection onto $\Ga_{\ga,A}$ we also have the isomorphisms
$$
\Ga_{\ga,A}\cong\ep^{\Z}
$$
and
$$
\Ga_{\ga,A}'\cong\ep^{k\Z}.
$$
The lemma follows.
\qed

Let $[\ga]\in\CE_P^{p,\reg}(\Ga)$.  Then $\ga$ is conjugate in $G$ to a matrix $a_{\ga}b_{\ga}$, where $a_{\ga}=\diag(a,a,a^{-1},a^{-1})$ for some $0<a<1$ and
$$
b_{\ga}=\matrixtwo{R(\th)}{R(\phi)}
$$
for some $\th$, $\phi$.  Recall that the length $l_{\ga}$ of $\ga$ is defined to be $8|\log a|$ and $N(\ga)=e^{l_{\ga}}$.  For $\CO\in O(S)$ let $\ep_{\CO}$ be a fundamental unit of $\CO$ and $R(\CO)=2|\log|\ep_{\CO}||$ the regulator of $\CO$ and let $r(\CO)=e^{4R(\CO)}$.  Under the map $\th$, the element $\ga$ corresponds to a fundamental unit of the order $\th([\ga])$ and it is clear that
$$
r\left(\th\left([\ga]\right)\right)=N(\ga).
$$
We recall from the proof of Lemma \ref{lem:centralisers} that $\th([\ga])$ is an order in the field $\Q(ae^{i\th},a^{-1}e^{i\phi})$.  By Lemma \ref{lem:SigmaTilde} we have that $\tr\tilde{\si}(\ga)=4(1-\cos 2\th)(1-\cos 2\phi)$.  Let
$$
\nu(\ga)=\prod_{\al}\left(1-\frac{\al(\ga)}{|\al(\ga)|}\right),
$$
where the product is over the embeddings of $F_{\ga}$ into $\C$.  A simple calculation then shows that
$$
\tr\tilde{\si}(\ga)=(1-e^{i\th})(1-e^{-i\th})(1-e^{i\phi})(1-e^{-i\phi})=\nu(\ga).
$$

We summarise the results of this section in the following:
\begin{proposition}
\label{pro:regOrders}
The map $\th$ is surjective and $\frac{4h(\CO)\la_S(\CO)\mu_{\CO}}{\ka(\CO)}$ to one.  For $[\ga]\in\CE_P^{p,\reg}(\Ga)$ we have that $\chi_1(\Ga_{\ga})=1/\mu_{\ga}$, that $\tr\tilde{\si}(\ga)=\nu(\ga)$ and that $N(\ga)=r\left(\th\left([\ga]\right)\right)$.
\end{proposition}

\section{Class numbers of orders in totally complex quartic fields}
We are now in a position to prove our main theorem.
\begin{theorem}\textnormal{(Main Theorem)}
\label{thm:Main}
\index{Main Theorem}
Let $S$ be a finite, non-empty set of prime numbers with an even number of elements.  For $x>0$ let
$$
\pi_S(x)=\sum_{{\CO\in O^c(S)}\atop{R(\CO)\leq x}}\la_S(\CO)h(\CO),\index{$\pi_S(x)$}
$$
Then, as $x\ra\infty$ we have
$$
\pi_S(x)\sim\frac{e^{4x}}{8x}.
$$
\end{theorem}
\prf
By \cite{Borel69}, Proposition 17.6 we know that $\Ga$ has a weakly neat subgroup $\Ga'$ which is of finite index.  Hence there exists $n_{\Ga}\in\N$ such that all roots of unity which are eigenvalues of $\Ad(\ga)$ for some $\ga\in\Ga$ have order less than or equal to $n_{\Ga}$.  Define the sets
$$
B^0_1=\left\{\matrixtwo{R(\th)}{R(\phi)}\in B:\th,\phi\notin\pi\Z\cup\frac{\pi}{2}\Z\cup\cdots\cup\frac{\pi}{n_{\Ga}}\Z\right\}
$$
and
$$
B^0_2=B^0_1\cup\left\{\matrixtwo{R(\th)}{R(\phi)}\in B:(\th+\phi),(\th-\phi)\notin\frac{\pi}{2}\Z\cup\frac{\pi}{3}\Z\right\}.
$$
For $x>0$ define the functions
$$
\pi^{\rw}_1(x)=\sum_{{[\ga]\in\CE^{p,\rw}_P(\Ga)}\atop{N(\ga)\leq x}}\chi_1(\Ga_{\ga})\ \textrm{ and }\ \ 
\pi^{\rw}_2(x)=\sum_{{[\ga]\in\CE^{p,\rw}_P(\Ga);\ka(\CO_{\ga})=1}\atop{N(\ga)\leq x}}\chi_1(\Ga_{\ga}).
$$
Setting $B^0=B^0_1$ in Theorem \ref{thm:PGT2} we get as $x\ra\infty$
\begin{equation}
\label{eqn:Pi1Asymptotic}
\pi^{\rw}_1(x)\sim\frac{2x}{\log x}.
\end{equation}
Setting $B^0=B^0_2$ in Theorem \ref{thm:PGT2} we get, using Lemma \ref{lem:KaCondition}, as $x\ra\infty$
\begin{equation}
\label{eqn:Pi2Asymptotic}
\pi^{\rw}_2(x)\sim\frac{2x}{\log x}.
\end{equation}
We also define the following functions for $x>0$:
$$
\pi_{S,1}(x)=\sum_{{\CO\in O^c(S)}\atop{R(\CO)\leq x}}\frac{\la_S(\CO)}{\ka(\CO)}h(\CO)\ \textrm{ and }\ \ 
\pi_{S,2}(x)=\sum_{{\CO\in O^c(S);\ka(\CO_{\ga})=1}\atop{R(\CO)\leq x}}\la_S(\CO)h(\CO).
$$
By Proposition \ref{pro:regOrders} and the fact that $N(\ga)=e^{4R(\CO_{\ga})}$ we deduce from (\ref{eqn:Pi1Asymptotic}) that
\begin{equation}
\label{eqn:PiS1Asymptotic}
\pi_{S,1}(x)\sim\frac{e^{4x}}{8x}
\end{equation}
and from (\ref{eqn:Pi2Asymptotic}) that
\begin{equation}
\label{eqn:PiS2Asymptotic}
\pi_{S,2}(x)\sim\frac{e^{4x}}{8x}.
\end{equation}
For $x>0$ define
$$
\pi'_{S,2}(x)=\pi_{S,1}(x)-\pi_{S,2}(x)=\sum_{{\CO\in O^c(S);\ka(\CO_{\ga})>1}\atop{R(\CO)\leq x}}\frac{\la_S(\CO)}{\ka(\CO)}h(\CO).
$$
Then from (\ref{eqn:PiS1Asymptotic}) and (\ref{eqn:PiS2Asymptotic}) it follows that
$$
\pi'_{S,2}(x)=o\left(\frac{e^{4x}}{x}\right).
$$
Since $\ka(\CO)\leq 4$ for all $\CO\in O^c(S)$ we have
$$
\pi_{S,2}(x)\leq\pi_S(x)\leq\pi_{S,2}(x)+4\pi'_{S,2}(x)
$$
and the theorem follows.
\qed

We can also prove the following:

\begin{theorem}
\label{thm:Main2}
\index{Main Theorem}
Let $S$ be a finite, non-empty set of prime numbers with an even number of elements.  For $x>0$ let
$$
\tilde{\pi}_S(x)=\sum_{{\CO\in O^c(S)}\atop{R(\CO)\leq x}}\nu(\CO)\la_S(\CO)h(\CO),\index{$\tilde{\pi}_S(x)$}
$$
where
$$
\nu(\CO)=\rez{2\mu_{\CO}}\sum_{\ep}\prod_{\al}\left(1-\frac{\al(\ep)}{|\al(\ep)|}\right),\index{$\nu(\CO)$}
$$
where the sum is over the $2\mu_{\CO}$ different fundamental units of $\CO$ and the product is over the embeddings of $\CO$ into $\C$.

Then, as $x\ra\infty$ we have
$$
\pi_S(x)\sim\frac{e^{4x}}{2x}.
$$
\end{theorem}
\prf
Let $\CE_P^{p,\rnw}(\Ga)=\CE_P^{p,\reg}(\Ga)\smallsetminus\CE_P^{p,\rw}(\Ga)$ be the set of regular, non weakly neat elements in $\CE_P^p(\Ga)$.  Setting $B^0=B^{\reg}$ in Theorem \ref{thm:PGT2} we get
$$
\sum_{{[\ga]\in\CE^{p,\reg}_P(\Ga)}\atop{N(\ga)\leq x}}\chi_1(\Ga_{\ga})\sim\frac{2x}{\log x}
$$
and with (\ref{eqn:Pi1Asymptotic}) above it follows that
$$
\sum_{{[\ga]\in\CE^{p,\rnw}_P(\Ga)}\atop{N(\ga)\leq x}}\chi_1(\Ga_{\ga})=o\left(\frac{x}{\log x}\right).
$$
It follows from Lemma \ref{lem:SigmaTilde} that $0\leq\tr\tilde{\si}(b_{\ga})\leq 16$ for all $\ga\in\CE_P(\Ga)$, and hence we have
$$
\sum_{{[\ga]\in\CE^{p,\rnw}_P(\Ga)}\atop{N(\ga)\leq x}}\chi_1(\Ga_{\ga})\tr\tilde{\si}(b_{\ga})=o\left(\frac{x}{\log x}\right).
$$
It then follows from Theorem \ref{thm:PGT} that
$$
\sum_{{[\ga]\in\CE^{p,\rw}_P(\Ga)}\atop{N(\ga)\leq x}}\chi_1(\Ga_{\ga})\tr\tilde{\si}(b_{\ga})\sim\frac{8x}{\log x}.
$$
We can also deduce from (\ref{eqn:Pi1Asymptotic}) and (\ref{eqn:Pi2Asymptotic}) that
$$
\sum_{{[\ga]\in\CE^{p,\rw}_P(\Ga);\ka(\CO_{\ga})>1}\atop{N(\ga)\leq x}}\chi_1(\Ga_{\ga})\tr\tilde{\si}(b_{\ga})=o\left(\frac{x}{\log x}\right).
$$
The theorem then follows from Proposition \ref{pro:regOrders} using the same arguments as in the proof of the previous theorem.
\qed

The methods we have used to deduce the asymptotic result for class numbers from the Prime Geodesic Theorem were not sharp enough to preserve the error term which was proven there.  However, we make the following conjecture:

\begin{conjecture}
Under the conditions of Theorem \ref{thm:Main}, as $x\ra\infty$ we have
$$
\pi_S(x)=\rez{2}L(4x)+O\left(\frac{e^{3x}}{x}\right),
$$
where
$$
L(x)=\int_1^x\frac{e^t}{t}\,dt.
$$
\end{conjecture}

\backmatter

\nocite{Deitmar96}
\nocite{Deitmar96a}
\nocite{Serre73}
\nocite{Humphreys78}
\nocite{Bourbaki89}
\nocite{Serre69}
\nocite{Helgason01}
\nocite{Juhl01}
\nocite{AlperinBell95}
\nocite{Dixmier77}
\nocite{BunkeOlbrich95}
\nocite{Karatsuba93}
\nocite{BrockerTomDieck95}
\nocite{Weil82}
\nocite{IwaniecKowalski04}
\nocite{Huber59}
\nocite{McKean72}
\nocite{Sarnak82}
\nocite{Hejhal76}
\nocite{Venkov82}
\nocite{Venkov79}
\nocite{Rademacher73}
\nocite{Gangolli77}
\nocite{Jost02}
\nocite{Reiter68}
\nocite{MoscoviciStanton91}

\newpage
\markright{\textnormal{Bibliography}}
\bibliography{/home/pavey/bibliography/biblio}
\bibliographystyle{acm}

\newpage
\markright{\textnormal{Index}}
\printindex

\newpage
\thispagestyle{empty}
\begin{center}\Large \textbf{Mark Pavey}\end{center}

Mathematisches Institut der Universit\"at T\"ubingen,

Auf der Morgenstelle 10,

72076,

T\"ubingen.

\bigskip

Geburtsdatum: 1/11/78

Geburtsort: Cheltenham, UK

\bigskip
\bigskip

\textbf{Lebenslauf und Bildungsgang:}

\begin{tabular}{ll}
Okt.2004-Gegenwart & Universit\"at T\"ubingen \\
 & Wissenschaftlicher Angestellter/Doktorand \\
\\

Okt.2002-Sept.2004 & Universit\"at Exeter, UK \\
 & Doktorand \\
\\

Sept.2001-Sept.2002 & RHS Harntec Ltd., Exeter, UK \\
 & Programmierer \\
\\

Okt.1998-Jun.2001 & St Edmund Hall, Universit\"at Oxford \\
 & BA(Hons) (Klasse I) \\
 & Mathematische Wissenschaften \\
\\

1990-1997 & Pates Grammar School, Cheltenham, UK \\

\end{tabular}

\end{document}